\documentclass[11pt]{article}

\usepackage[margin=1in]{geometry}
\usepackage{graphicx}
\usepackage{multirow}
\usepackage{amsmath,amssymb,amsfonts,amsthm}
\usepackage{mathrsfs}
\usepackage[title]{appendix}
\usepackage{xcolor}
\usepackage{textcomp}
\usepackage{booktabs}
\usepackage{algorithm}
\usepackage{algorithmicx}
\usepackage{algpseudocode}
\usepackage{enumitem}
\usepackage{bbm}
\usepackage{parskip}
\usepackage{makecell}
\usepackage{array}
\usepackage{caption}
\usepackage{tikz}
\usepackage{subcaption}
\usepackage{float}
\usepackage{hyperref}
\usepackage{placeins}
\usepackage[numbers]{natbib}

\usetikzlibrary{shapes.geometric, arrows.meta, positioning, calc, shapes.multipart, decorations.pathreplacing, matrix, shapes.misc, backgrounds, fit}

\theoremstyle{plain}

\newtheorem{proposition}{Proposition}[section]
\newtheorem{lemma}{Lemma}[section]

\theoremstyle{remark}

\numberwithin{equation}{section}

\begin{document}

\title{Two-Product Make-to-Stock System: \\ Strategic Joining and Optimal Inventory Levels}

\author{
Odysseas Kanavetas\thanks{Mathematical Institute, Leiden University, Einsteinweg 55, Leiden 2333CC, Netherlands. Email: \texttt{o.kanavetas@math.leidenuniv.nl}}
\and
Ekaterina Kosarevskaia\thanks{Mathematical Institute, Leiden University, Einsteinweg 55, Leiden 2333CC, Netherlands. Email: \texttt{e.kosarevskaia@math.leidenuniv.nl}}
}

\date{}

\maketitle

\begin{abstract}
This paper analyzes a two-product make-to-stock queueing system where a single production facility serves two customer classes with independent Poisson arrivals. Customers make strategic join-or-balk decisions without observing current inventory levels. The analysis establishes the existence and uniqueness of Nash equilibria in customer joining strategies for various inventory scenarios. Optimal base-stock levels are characterized from both profit-maximizing and welfare-maximizing perspectives, with closed-form expressions for key performance measures. 
\end{abstract}


\newcommand*{\lambar}[1]{{\lambda}_{#1}}
\newcommand*{\SWe}{\mathrm{SW}}
\newcommand*{\BRf}[1]{\mathrm{BR}_{#1}}

\section{Introduction}
Make-to-stock systems aim to hold inventory in anticipation of uncertain demand, yet the moment customers arrive, their individual choices can reshape the entire operation. In multi-product settings, where a shared production facility must manufacture different product types, these dynamics become particularly complex. Producers must determine appropriate target inventory levels for each distinct product while balancing holding costs against stockout risks. When customers strategically decide whether to purchase specific products, their collective behavior significantly impacts overall system performance. This resource-sharing challenge appears across industries --- from manufacturing facilities producing both premium and economy variants to cloud computing systems allocating capacity between interactive and batch processing services.

\textit{--~Classical inventory models.} Inventory management research dates to \citet{Harris1913} economic order quantity model. Subsequent work incorporated stochastic demand and dynamic decisions: \citet{Arrow1951} on optimal inventory policies, \citet{Scarf1960} on $(S,s)$ policies, and \citet{Karlin1958} on dynamic inventory control. Multi-product and nonstationary extensions followed \citep{Veinott1965}, with comprehensive treatments in \citet{Zipkin2000}, \citet{Silver1998}, and \citet{Graves2000}. Multi-item models with shared capacity include \citet{Federgruen1996}. 

\textit{--~Inventory models with queues and make-to-stock systems.} Queueing formulations of make-to-stock systems model the inventory shortfall as a single-server queue. \citet{Buzacott1993} develop this representation systematically; the non-strategic version of our two-product system corresponds to one of the standard make-to-stock queues in their framework. \citet{Veatch1996} use the same representation to derive optimal dynamic control in a tandem production--inventory system. Closer to our setting, \citet{Ha1997} and \citet{DeVericourt2002} study single-item make-to-stock systems with multiple demand classes sharing a production resource, showing that optimal policies take the form of inventory-rationing or stock-allocation rules. \citet{Benjaafar2012}, \citet{Abouee2012}, and \citet{Zhao2011} likewise analyze multiclass or multistage queueing--inventory systems under target inventory and rationing policies.

\textit{--~Strategic queues without inventory.} In parallel, another stream models customers as strategic agents deciding whether to join a queue. \citet{Naor1969} analyzes joining behavior in an M/M/1 queue and shows how tolls align individual and social objectives. \citet{Edelson1975} study congestion pricing for Poisson queues. Extensions include \citet{Hassin1985}, \citet{Mendelson1985} on computer service pricing, \citet{Burnetas2007}, and the monographs of \citet{HassinHaviv2003} and \citet{Hassin2016}. \citet{HassinRoetGreen2017} consider an unobservable queue where customers may pay to inspect queue length before joining, characterizing equilibrium inspection and joining behavior.

\textit{--~Make-to-stock queues with demand responsiveness.} Within the make-to-stock framework, several papers incorporate demand responsiveness to prices or congestion at an aggregate level. \citet{Caldentey2006} analyze a single-product system with two selling channels, using dynamic pricing and admission control. \citet{Ata2006} study dynamic control and quantify the value of demand information. \citet{Cachon2009} consider strategic consumers who may delay purchases in anticipation of markdowns.

Most closely related to our work is the literature on single-product make-to-stock queues with strategic joining decisions. \citet{Oz2015} study a single-product queue with unobservable inventory, where customers choose equilibrium joining probabilities based on stationary performance measures; their comparison of decentralized and centralized objectives makes this the closest single-product analogue to our work. Closely related extensions retain strategic joining but differ mainly in the producer's control lever: \citet{Li2016} consider an unobservable finite-capacity system with setup times, characterizing optimal two-threshold production policies. \citet{Zhang2019} analyze dynamic-service make-to-stock systems, identifying optimal inventory thresholds in a Stackelberg game. Other variants modify the customer-flow dynamics or add additional instruments: \citet{WangCai2023} extend the framework to retrial customers. \citet{Cai2020} study pricing and information disclosure under different observability structures. Finally, \citet{Zare2021} incorporate strategic customers in a two-echelon supply chain, focusing on vertical manufacturer--distributor interactions rather than horizontal product competition.

Building on the above make-to-stock queueing-inventory and strategic-joining literature, we study a two-product setting where both product classes share a single production server under target inventory (base-stock) policies, and customers make joining decisions without observing inventory or backlog levels. The analysis explores how target inventory levels and price/waiting trade-offs influence customer joining strategies and system performance. Different inventory targets for each product affect equilibrium joining probabilities, which in turn determine stockouts and waiting times. We examine both a profit-maximizing producer's perspective and a social planner's view of overall welfare, highlighting the externalities that arise when products share production capacity but face distinct customer segments.
These findings also apply to settings beyond conventional manufacturing and retail. Shared resources in healthcare (e.g., diagnostic equipment), cloud computing (e.g., server allocation), or transportation (e.g., electric vehicle charging) all involve users who make decisions under imperfect information. Modern computing environments must balance capacity across multiple service tiers while accommodating strategic users; see \citet{Gong2018, Hong2017}. Accounting for strategic behavior is thus critical for coordinating resource allocation and managing costs.

The remainder of the paper is organized as follows. Section~\ref{sec:model} introduces our two-product make-to-stock system, defining the key parameters and interactions. Section~\ref{sec:fundamentals} establishes the analytical foundation by deriving closed-form expressions for critical performance measures --- expected waiting times, inventory levels, and backlogs. Sections~\ref{sec:customer} and~\ref{sec:producer} analyze the system from a decentralized perspective, where decisions are made independently by participants. Section~\ref{sec:customer} explores how customers strategically decide whether to join the system, characterizing the structure of Nash equilibria under different inventory configurations. Section~\ref{sec:producer} shifts to the producer's viewpoint in a Stackelberg game framework, where the producer acts as a leader by setting optimal inventory targets while anticipating customer responses. In contrast, Section~\ref{sec:social} examines a centralized approach where a social planner coordinates both inventory and customer participation decisions simultaneously to maximize overall welfare rather than individual objectives. Section~\ref{sec:numerics} presents our numerical findings on how asymmetry in customer waiting costs affects system performance when products compete for limited production capacity.

\section{Model Description}\label{sec:model}
\subsection{Single-Product Make-to-Stock System}\label{subsec:single-product}
We begin with the single-product make-to-stock system, following \citet{Buzacott1993} and \citet{Oz2015}, which serves as the foundation for our two-product extension.

A single production facility maintains a target inventory level (base-stock level) $S\in\mathbb{Z}_+$. Customers arrive as a Poisson process with rate $\Lambda$ and, without observing the current inventory or queue length, decide whether to join or balk. A joining customer pays $p$ and either receives a unit immediately or, if inventory is zero, waits in the backlog. In either case, their arrival triggers exactly one production job added to the queue.

The facility processes jobs first-come-first-served (FCFS) with exponential service times at rate $\mu$. Crucially, the job created by an arrival is not tied to that customer's unit. For example, let $S=2$ and suppose inventory is zero, no customer is backlogged, and two jobs are in the queue (both currently destined to rebuild inventory). A customer who joins is backlogged; the next completion serves this customer, while the job their arrival created later rebuilds inventory (provided no additional customers arrive before it completes).

\subsection{Two-Product Extension}\label{subsec:two-product}

We now extend to two products sharing a single production line, with target inventory levels $S_1$ and $S_2$. Customers requesting product $i$ arrive as independent Poisson processes with rate $\Lambda_i$, and do not observe current inventories or queue lengths when deciding whether to join.

The production mechanics extend directly: a joining type-$i$ customer pays $p_i$ and triggers one type-$i$ job appended to a single FCFS queue shared by both products. The customer receives the product immediately if in stock, or waits in the backlog otherwise. The production line processes all jobs at rate $\mu$; a completed type-$i$ job serves the oldest waiting type-$i$ customer if one exists, otherwise replenishes on-hand inventory. Holding cost is $h_i$ per unit of on-hand inventory per unit time.

While FCFS is natural in the single-product model (all jobs identical), it need not be optimal with two products sharing capacity: production capacity may be spent on type-1 replenishment even while type-2 customers are waiting. Studying alternative service disciplines is an interesting direction for future research; we keep FCFS as a tractable baseline that still captures competition for shared capacity.

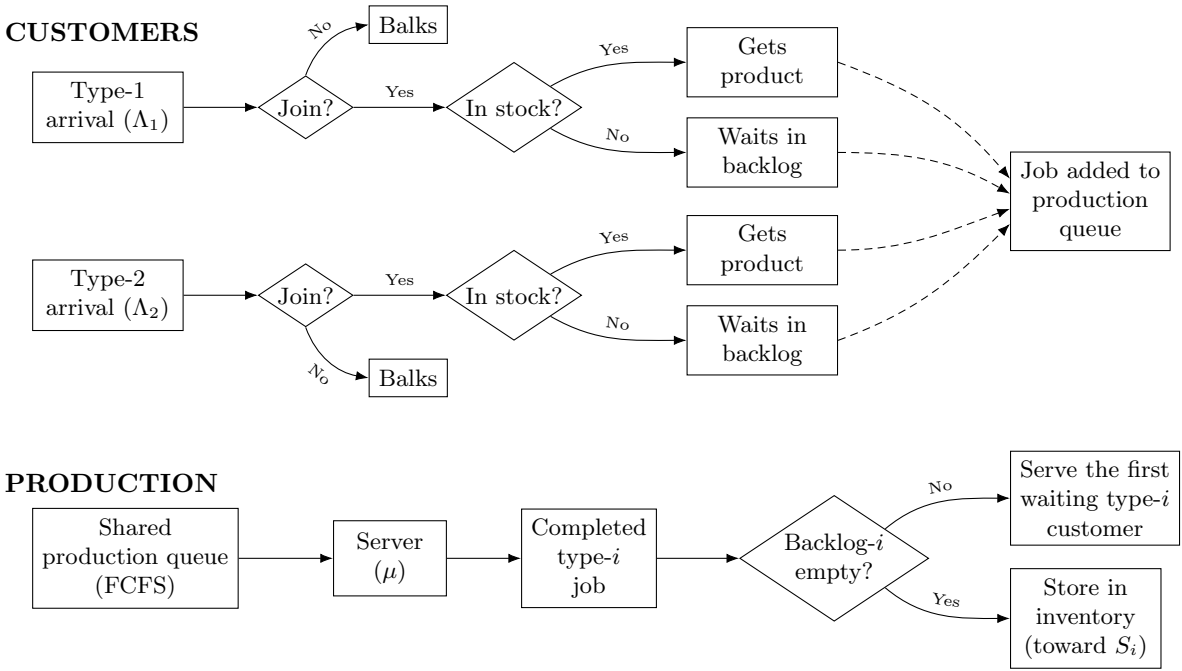
\begin{figure}[!htbp]
\begin{tikzpicture}[>=Latex]

\tikzset{
    rectangle_node/.style={
        rectangle,
        minimum width=2cm,
        minimum height=0.8cm,
        draw,
        align=center,
        font=\footnotesize,
        anchor=west
    },
    small_rectangle_node/.style={
        rectangle,
        draw,
        align=center,
        font=\footnotesize,
        minimum height=0.5cm,
        anchor=west
    },
    floor2_rectangle_node/.style={
        rectangle,
        draw,
        align=center,
        font=\footnotesize,
        minimum width=1.5cm,
        minimum height=1cm,
        anchor=west
    },
    diamond_node/.style={
        diamond,
        aspect=1.5,
        draw,
        align=center,
        font=\footnotesize,
        inner sep=1pt,
        anchor=west
    },
    arrow/.style={->},
    dashed_arrow/.style={->, densely dashed},
    floor_title/.style={font=\small\bfseries, anchor=west}
}


\def\yTypeOne{4.5}
\def\yTypeTwo{2}
\def\yProd{-1.5}

\def\xArrival{0}
\def\xJoin{3}
\def\xStock{5.5}
\def\xOutcome{8.7}          
\def\xJobLabel{13}      

\def\xProdQueue{0}
\def\xServer{4}
\def\xCompleted{6.5}
\def\xBacklogCheck{9.4}
\def\xFinal{13}           

\def\dyBalk{1.1}          
\def\dyBranch{0.6}        
\def\dyFinal{0.8}         

\def\xTitle{-0.5}
\def\yCustomersTitle{5.5}
\def\yProdTitle{-0.5}

\node[floor_title] at (\xTitle, \yCustomersTitle) {CUSTOMERS};
\node[floor_title] at (\xTitle, \yProdTitle) {PRODUCTION};

\node[rectangle_node] (arr1) at (\xArrival, \yTypeOne) {Type-1\\arrival ($\Lambda_1$)};
\node[diamond_node]   (join1) at (\xJoin, \yTypeOne) {Join?};
\node[diamond_node]   (stk1) at (\xStock, \yTypeOne) {In stock?};
\node[rectangle_node] (get1) at (\xOutcome, \yTypeOne + \dyBranch) {Gets\\product};
\node[rectangle_node] (bl1) at (\xOutcome, \yTypeOne - \dyBranch) {Waits in\\backlog};

\node[small_rectangle_node] (balk1) at ($(join1.east)+(0.2, \dyBalk)$) {Balks};

\draw[arrow] (arr1.east) -- (join1.west);
\draw[arrow] (join1.east) -- node[above, font=\tiny] {Yes} (stk1.west);
\draw[arrow] (join1.north) to[out=70,in=190] node[above, sloped, font=\tiny] {No} (balk1.west);
\draw[arrow] (stk1.30) to[out=30,in=180] node[above, sloped, font=\tiny] {Yes} (get1.west);
\draw[arrow] (stk1.-30) to[out=-30,in=180] node[above, sloped, font=\tiny] {No} (bl1.west);

\node[rectangle_node] (arr2) at (\xArrival, \yTypeTwo) {Type-2\\arrival ($\Lambda_2$)};
\node[diamond_node]   (join2) at (\xJoin, \yTypeTwo) {Join?};
\node[diamond_node]   (stk2) at (\xStock, \yTypeTwo) {In stock?};
\node[rectangle_node] (get2) at (\xOutcome, \yTypeTwo + \dyBranch) {Gets\\product};
\node[rectangle_node] (bl2) at (\xOutcome, \yTypeTwo - \dyBranch) {Waits in\\backlog};

\node[small_rectangle_node] (balk2) at ($(join2.east)+(0.2, -\dyBalk)$) {Balks};

\draw[arrow] (arr2.east) -- (join2.west);
\draw[arrow] (join2.east) -- node[above, font=\tiny] {Yes} (stk2.west);
\draw[arrow] (join2.south) to[out=-70,in=170] node[below, sloped, font=\tiny] {No} (balk2.west);
\draw[arrow] (stk2.30) to[out=30,in=180] node[above, sloped, font=\tiny] {Yes} (get2.west);
\draw[arrow] (stk2.-30) to[out=-30,in=180] node[above, sloped, font=\tiny] {No} (bl2.west);

\pgfmathsetmacro{\yJobLabel}{(\yTypeOne + \yTypeTwo) / 2}
\node[rectangle_node] (joblabel) at (\xJobLabel, \yJobLabel) {Job added to\\production \\queue};

\draw[dashed_arrow] (get1.east) to[out=-20,in=130] ($(joblabel.west)+(0,0.3)$);
\draw[dashed_arrow] (bl1.east)  to[out=0,in=150] ($(joblabel.west)+(0,0.1)$);
\draw[dashed_arrow] (get2.east) to[out=0,in=200] ($(joblabel.west)+(0,-0.1)$);
\draw[dashed_arrow] (bl2.east)  to[out=20,in=225] ($(joblabel.west)+(0,-0.3)$);

\node[floor2_rectangle_node] (prodqueue) at (\xProdQueue, \yProd) {Shared \\production queue \\(FCFS)};
\node[floor2_rectangle_node] (prodserver) at (\xServer, \yProd) {Server\\($\mu$)};
\node[floor2_rectangle_node] (completed) at (\xCompleted, \yProd) {Completed\\type-$i$\\ job};
\node[diamond_node]   (backlogempty) at (\xBacklogCheck, \yProd) {Backlog-$i$\\empty?};
\node[rectangle_node] (towaiting) at (\xFinal, \yProd + \dyFinal) {Serve the first\\waiting type-$i$\\customer};
\node[rectangle_node] (toinventory) at (\xFinal, \yProd - \dyFinal) {Store in\\inventory\\(toward $S_i$)};

\draw[arrow] (prodqueue.east) -- (prodserver.west);
\draw[arrow] (prodserver.east) -- (completed.west);
\draw[arrow] (completed.east) -- (backlogempty.west);
\draw[arrow] (backlogempty.30) to[out=40,in=180] node[above, sloped, font=\tiny] {No} (towaiting.west);
\draw[arrow] (backlogempty.-30) to[out=-40,in=180] node[above, sloped, font=\tiny] {Yes} (toinventory.west);

\end{tikzpicture}
\caption{Two-Product Make-to-Stock System: Customer Flow and Production Process}
\label{fig:system-flow}
\end{figure}

Figure~\ref{fig:system-flow} illustrates this system. Key parameters include prices $p_i$, rewards $R_i > p_i$, waiting costs $c_i$ per unit time, and holding costs $h_i$. While $\Lambda_i$ denotes the arrival rate, $\lambar{i}$ denotes the effective arrival rate of customers who actually join. In the decentralized setting, $\lambar{i} = \Lambda_i q_i$, where $q_i$ is the joining probability; under the social planner, $\lambar{i}$ is directly chosen as the optimal participation rate.

This system shows a key tradeoff: higher inventory reduces customer waiting times but increases holding costs. Because customers cannot see current inventory levels, they must make joining decisions based on expected rather than actual waiting times, similar to how retail shoppers decide whether to visit a store without knowing its current stock levels.

\section{Queueing Analysis and Performance Measures}
\label{sec:fundamentals}
We now derive closed-form expressions for waiting times, inventory levels, and backlogs. These formulas support our later strategic analysis while also providing practical insights for non-strategic settings (where $\lambar{i} = \Lambda_i$). The expressions reveal how waiting times, inventory levels, and backlogs respond to system parameters, demonstrating the inherent tradeoffs when two products share production capacity.

\subsection{Stationary Distributions}
\subsubsection{Baseline Results}
\label{subsec:known-distr}
This subsection summarizes stationary distributions established by \citet[Section~4.4]{Buzacott1993} for two-product make-to-stock systems, which serve as the foundation for our analysis. In this setting, a single production server processes both product types at rate $\mu$, and the system maintains a target inventory level of $S_i$ units for each product~$i$.  If $N_i(t)$, the number of product-$i$ jobs in the production queue, exceeds $S_i$, the difference $B_i(t) = N_i(t)-S_i$ is interpreted as a backlog of customers who arrived to find no finished inventory available.  

Conversely, $I_i(t) = [S_i - N_i(t)]^+$ gives the on-hand inventory of product~$i$ whenever $N_i(t)$ falls below $S_i$.  Hence, $N_i$, $I_i$, and $B_i$ capture the state of the system from different angles:
\begin{itemize}
  \item $N_i(t)$: the number of type-$i$ jobs in the production queue (either waiting or in service),
  \item $I_i(t)$: the on-hand inventory: units of product~$i$ available to immediately serve arrivals,
  \item $B_i(t)$: the backlog of type-$i$ customers still waiting for a unit to be produced.
\end{itemize}
These three variables provide a complete picture of the system's state. When inventory level is positive, customers receive products immediately upon arrival. As inventory level depletes and jobs accumulate in the queue, $N_i(t)$ increases. Once number of jobs exceeds inventory level, the backlog grows, indicating customers who must wait for production to complete. The balance between these measures directly impacts customer satisfaction and operational costs --- high inventory reduces waiting but increases holding costs, while high backlog indicates poor service quality.

Across products, $N(t) = N_1(t) + N_2(t)$ is simply the total number of jobs in the queue. Although the processes $\{N_1(t)\}$ and $\{N_2(t)\}$ are not Markov individually, their sum $N(t)$ behaves as an M/M/1 queue with total arrival rate $\lambar{1} + \lambar{2}$ and service rate $\mu$. In the long run, the system reaches a steady state where the probability distribution of queue length no longer changes with time. This stationary distribution is:
\begin{equation}\label{eq:distr_jobs_total}
\mathbb{P}\{N = n\}
=
\left(1 - \frac{\lambar{1} + \lambar{2}}{\mu}\right)
\left(\frac{\lambar{1} + \lambar{2}}{\mu}\right)^n,
\quad
n = 0,1,2,\dots.
\end{equation}
Since arrivals of types~1 and~2 are independent Poisson processes with rates $\lambar{1}$ and $\lambar{2}$, each job in $N$ is of type~$i$ with probability $\lambar{i} / (\lambar{1} + \lambar{2})$.  By multinomial thinning, conditional on $N = n$, the vector $(N_1, N_2)$ follows a multinomial distribution:
 \[
\mathbb{P}\{N_1 = n_1, N_2 = n_2 \mid N = n\}
=
\binom{n}{n_1} \left(\frac{\lambar{1}}{\lambar{1}+\lambar{2}}\right)^{n_1} \left(\frac{\lambar{2}}{\lambar{1}+\lambar{2}}\right)^{n_2}, \quad n_1 + n_2 = n.
\]
Combining this with \eqref{eq:distr_jobs_total} yields the joint distribution
\begin{align*}	
\mathbb{P}&\{N_1 = n_1, N_2 = n_2\}\\
&=
\left(1 - \frac{\lambar{1} + \lambar{2}}{\mu}\right)
\left(\frac{\lambar{1} + \lambar{2}}{\mu}\right)^{n_1 + n_2}
\binom{n_1 + n_2}{n_1} \left(\frac{\lambar{1}}{\lambar{1}+\lambar{2}}\right)^{n_1} \left(\frac{\lambar{2}}{\lambar{1}+\lambar{2}}\right)^{n_2}.
\end{align*}
Summing over $n_j$ and applying standard geometric series identities, the marginal distribution of $N_i$ takes a geometric form:
\[
\mathbb{P}\{N_i = n_i\}
=
\left(1 - \frac{\lambar{i}}{\mu - \lambar{j}}\right)
\left(\frac{\lambar{i}}{\mu - \lambar{j}}\right)^{n_i},
\quad
\text{for } i,j \in \{1,2\},\; j \neq i.
\]
The stationary distribution of on-hand inventory $I_i$ follows by applying the transformation $I_i = [S_i - N_i]^+$ to the geometric distribution of $N_i$. Since $I_i = 0$ corresponds to $N_i \ge S_i$ and $I_i = k \in \{1,\dots,S_i\}$ corresponds to $N_i = S_i - k$, we obtain:
\begin{equation}\label{eq:distr_inv_level}
\mathbb{P}\{I_i = k\}
=
\begin{cases}
\displaystyle
\left(\frac{\lambar{i}}{\mu - \lambar{j}}\right)^{S_i}, 
& k = 0,\\[6pt]
\displaystyle
\left(1 - \frac{\lambar{i}}{\mu - \lambar{j}}\right)
\left(\frac{\lambar{i}}{\mu - \lambar{j}}\right)^{S_i - k},
& k = 1,\dots,S_i,
\end{cases}
\quad
j \neq i.
\end{equation}

The probability that $I_i$ is zero captures a stockout situation, meaning all $S_i$ target units are in the production queue. When the queue has fewer than $S_i$ jobs, some finished inventory remains available.

The stationary distribution of backlog $B_i$ is obtained similarly via $B_i = [N_i - S_i]^+$. The event $B_i = 0$ corresponds to $N_i \le S_i$, whereas $B_i = k \ge 1$ corresponds to $N_i = S_i + k$, yielding the piecewise geometric form:
\begin{equation}\label{eq:distr_backlog}
\mathbb{P}\{B_i = k\}
=
\begin{cases}
1 - \displaystyle\left(\frac{\lambar{i}}{\mu - \lambar{j}}\right)^{S_i + 1}, 
& k = 0,\\[6pt]
\displaystyle
\left(1 - \frac{\lambar{i}}{\mu - \lambar{j}}\right)
\left(\frac{\lambar{i}}{\mu - \lambar{j}}\right)^{S_i + k},
& k = 1,2,\dots.
\end{cases}
\end{equation}

The event $B_i=0$ means the system has sufficient inventory to prevent backlog, while $B_i>0$ indicates customers are waiting for products to be manufactured.

The next parts of this section leverage these results to derive closed-form expressions for waiting times and expected on-hand inventories, then expand to backlog expectations as well.

\subsubsection{Queue Position of Arriving Customer}
\label{subsec:arriving-cust-position}
To derive the expected waiting time, we identify how many jobs must complete before a newly arriving type-$i$ customer receives a unit. Since the arrival triggers a new type-$i$ job at the tail of the queue, one might assume the customer must wait for that job. The situation is more subtle.

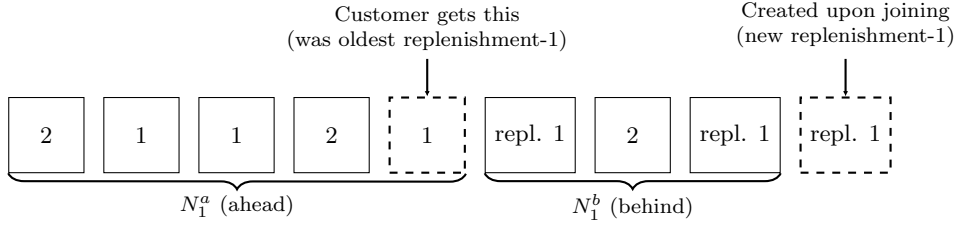
\begin{figure}[!htbp]
\centering
\begin{tikzpicture}[
  job/.style={rectangle, draw, minimum height=1cm, minimum width=1cm, align=center, font=\footnotesize, text=black},
  type1/.style={},
  type2/.style={},
  newjob/.style={thick, dashed},
  brace/.style={decorate, decoration={brace, mirror, amplitude=6pt}, thick},
  label/.style={font=\scriptsize},
  arrow/.style={-{Latex[width=2pt,length=2pt]}, thick}
]

\node[job, type2] (j1) at (0,0) {2};
\node[job, type1] (j2) [right=0.25cm of j1] {1};
\node[job, type1] (j3) [right=0.25cm of j2] {1};
\node[job, type2] (j4) [right=0.25cm of j3] {2};
\node[job, newjob] (j5) [right=0.25cm of j4] {1};
\node[job, type1] (j6) [right=0.25cm of j5] {repl. 1};
\node[job, type2] (j7) [right=0.25cm of j6] {2};
\node[job, type1] (j8) [right=0.25cm of j7] {repl. 1};
\node[job, newjob] (j9) [right=0.25cm of j8] {repl. 1};

\draw[brace] (j1.south west) -- node[below=4pt, label] {$N_1^a$ (ahead)} (j5.south east);
\draw[brace] (j6.south west) -- node[below=4pt, label] {$N_1^b$ (behind)} (j8.south east);

\draw[arrow] ($(j5.north)+(0,0.5)$) -- (j5.north);
\node[align=center, font=\scriptsize] at ($(j5.north)+(0,0.9)$) {
  Customer gets this\\ (was oldest replenishment-1)
};

\draw[arrow] ($(j9.north)+(0,0.5)$) -- (j9.north);
\node[align=center, font=\scriptsize] at ($(j9.north)+(0,0.95)$) {
  Created upon joining\\ (new replenishment-1)
};

\end{tikzpicture}
\caption{Queue Position for Backlogged Arrival}
\label{fig:queue-position}
\end{figure}

Throughout this subsection, we use two state-dependent labels for jobs in the production queue: backlog-serving (completions allocated to currently waiting customers) and replenishment (completions that rebuild inventory toward $S_i$ once backlog is cleared). These labels may shift with the next arriving customer; they do not affect production and are used only for the waiting-time derivation. Under a target inventory level $S_i>0$ and a stockout event ($I_i=0$), the queue contains $S_i$ type-$i$ replenishment jobs (possibly preceded by type-$i$ backlog-serving jobs). A new type-$i$ customer who joins is served by the completion of the job that, immediately before the arrival, was the oldest replenishment type-$i$ job already in the queue. The newly triggered job enters at the tail as the newest replenishment job. Figure~\ref{fig:queue-position} illustrates this mechanism for the case where $i=1$ and shows how type-2 jobs may be interspersed throughout.

Let $N_i^a(t)$ be the number of jobs ahead of, and including the serving type-$i$ job (i.e., the job whose completion supplies a newly arriving type-$i$ backorder at time $t$, after any earlier backorders). Let $N_i^b(t)$ be the number of jobs strictly behind it. Hence, at time $t$, before any potential arrival,
\[
N(t)  =  N_i^a(t)  +  N_i^b(t),
\]
where (see Section~\ref{subsec:known-distr}) $N(t)$ follows an M/M/1-type stationary distribution under effective arrival rates $\lambar{1}, \lambar{2}$, and service rate $\mu$.

The joint stationary distribution of $\left(N_i^a(t), N_i^b(t)\right)$ is:
\begin{equation}\label{eq:joint-position}
\mathbb{P}\{N_i^a = n^a,N_i^b = n^b\}
=
\binom{n^b}{S_i-1}
\frac{(\lambar{1}+\lambar{2})^{n^a-1}\lambar{i}^{S_i}\lambar{j}^{n^b - S_i+1}\left(\mu - \lambar{1} - \lambar{2}\right)}
     {\mu^{n^a + n^b + 1}}, \quad j\neq i,
\end{equation}
valid for $n^a \ge 1$ and $n^b \ge S_i-1$, and zero otherwise. The binomial coefficient arises because exactly $S_i-1$ of the $n^b$ jobs behind the customer's position must be type-$i$ replenishment jobs, with the remaining jobs assigned types through standard thinning of the Poisson process. For the detailed proof, see Appendix \ref{app:position-distribution}.

This distribution specifies the 'position in queue' $N_i^a$ for an arriving type-$i$ backlog under the condition of out-of-stock. The next subsection uses $\mathbb{E}[N_i^a]$ to derive an explicit formula for the expected waiting time of that newly arrived customer.
\subsection{Waiting Times}
\label{sec:waiting_time}
When a customer arrives to find no inventory available, their expected waiting time depends on the system's congestion level, the inventory policy, and the interaction between products sharing service capacity.

For a newly arrived type-$i$ customer, the expected waiting time $\mathbb{E}W_i(\lambar{1}, \lambar{2})$ equals the expected number of jobs ahead of it $N_i^a$ multiplied by the mean service time $1/\mu$. By analyzing the queue structure and applying appropriate summation techniques on the joint distribution (see Appendix \ref{app:waiting-time} for details), we obtain:

\begin{proposition}\label{prop:waiting_time}
The expected waiting time for a newly arriving type-$i$ customer is
\begin{equation}\label{eq:waiting_time1}
\mathbb{E}[W_i(\lambar{1}, \lambar{2})]
=
\left(\frac{\lambar{i}}{\mu - \lambar{j}}\right)^{S_i}
\cdot
\frac{1}{\mu - \lambar{1} - \lambar{2}}, \qquad j\neq i.
\end{equation}
\end{proposition}

This formula reveals several key operational insights. First, waiting time decreases as $S_i$ increases, since $\lambar{i} < \mu - \lambar{j}$ for stability, making $(\frac{\lambar{i}}{\mu - \lambar{j}})^{S_i}$ smaller with larger inventory targets. The fraction $\frac{\lambar{i}}{\mu - \lambar{j}}$ represents a product-specific congestion measure --- the ratio between product $i$'s demand and its effective service capacity (after accounting for capacity consumed by product $j$). The power $S_i$ shows how inventory buffers against this congestion. The term $(\mu - \lambar{1}-\lambar{2})$ in the denominator reflects overall system utilization --- as total arrival rates approach service capacity, waiting times increase for all products. 

Notably, this expression depends on the other product only through its arrival rate $\lambar{j}$, not its inventory level $S_j$. This happens because each arriving customer of type $j$ adds exactly one job to the production queue, regardless of whether it is a replenishment job or a backlog-serving job. From the perspective of a type-$i$ customer, the distinction between these job types for product $j$ is irrelevant --- what matters is only how many total jobs from product $j$ are ahead in the queue, which is determined solely by $\lambar{j}$.

\subsection{Inventory Levels}\label{sec:inv_lev}
Using the stationary distribution \eqref{eq:distr_inv_level} of on-hand inventory $I_i(\lambar{1}, \lambar{2})$, we compute the expected inventory for each product. Through appropriate summation and application of geometric series identities, we obtain:
\begin{equation}\label{eq:inv_level}
\mathbb{E}[I_i(\lambar{1}, \lambar{2})]
=
S_i 
- \frac{\lambar{i}}{\mu - \lambar{1} - \lambar{2}}
\left[
  1 
  - 
  \left(\tfrac{\lambar{i}}{\mu - \lambar{j}}\right)^{S_i}
\right],
\quad
j \neq i.
\end{equation}

This expression reveals how expected inventory balances target levels against system dynamics. The starting point $S_i$ represents the target inventory level, which is then reduced by the expected depletion term. This depletion depends on two congestion measures: $\frac{\lambar{i}}{\mu - \lambar{1} - \lambar{2}}$ captures overall system congestion, while $(\frac{\lambar{i}}{\mu - \lambar{j}})^{S_i}$ reflects product-specific stockout probability. As $S_i$ increases, this stockout probability decreases (since $\frac{\lambar{i}}{\mu - \lambar{j}} < 1$ for stability), meaning the bracketed term approaches 1. This creates a diminishing returns effect --- each additional unit of target inventory contributes less to the expected on-hand inventory. Under high arrival rates approaching capacity limits, this effect strengthens, as both congestion measures increase and depletion accelerates.

\subsection{Backlog Levels}\label{sec:backlog_lev}
For the backlog $B_i(\lambar{1}, \lambar{2})$, we derive the expected value using the distribution from \eqref{eq:distr_backlog}.

Summing over all possible backlog values and applying the identity for infinite geometric series, we obtain:
\[
\mathbb{E}[B_i(\lambar{1}, \lambar{2})]
= \left(\frac{\lambar{i}}{\mu - \lambar{j}}\right)^{S_i}
\cdot
\frac{\lambar{i}}
     {\bigl(\mu - \lambar{1} - \lambar{2}\bigr)}.
\]
This expression reveals how expected backlog decreases with higher inventory targets, sharing structural similarities with the waiting time formula. The term $\frac{\lambar{i}^{S_i+1}}{(\mu - \lambar{j})^{S_i}}$ incorporates the same product-specific congestion measure raised to power $S_i$, showing how inventory buffers against stockouts. The additional $\lambar{i}$ factor reflects that backlog only occurs for customers who actually join the system. The denominator term $(\mu - \lambar{1} - \lambar{2})$ again captures overall system utilization --- as total demand approaches capacity, backlog grows regardless of inventory levels. Higher $S_i$ values reduce expected backlog by decreasing the probability of stockout events, though this protection becomes less effective as congestion increases.

\section{Strategic Customer Behavior}\label{sec:customer}
Using the performance measures from Section~\ref{sec:fundamentals}, we now begin our examination of the decentralized decision-making approach. In this framework, customers and producers act independently to maximize their own objectives rather than coordinating decisions. This section focuses on the first stage of decentralization: how customers strategically decide whether to join the system based on their expected utilities, considering both the producer's inventory choices and the anticipated behavior of other customers. Section~\ref{sec:producer} will complete this decentralized perspective by examining how the producer optimizes inventory targets while anticipating these customer responses.

\subsection{Utility Functions and Best Response}
Consider a type-$i$ customer ($i=1,2$) who observes strategies $q_1,q_2$ being played by other customers in the system. These joining probabilities determine the effective arrival rates $\lambar{i} = q_i\Lambda_i$ used in our earlier performance measures. The customer's expected utility from joining is:
\[
U_i(\mathrm{join}; q_1, q_2) = R_i - p_i - c_i \mathbb{E}\left[W_i(q_1\Lambda_1,q_2\Lambda_2)\right],
\]
where $R_i$ is the reward, $p_i$ is the price, $c_i$ is the waiting-cost rate, and $\mathbb{E}[W_i(\Lambda_1 q_1,\Lambda_2 q_2)]$ is the expected waiting time given others' strategies (Section~\ref{sec:waiting_time}).

If the customer balks, the utility $U_i(\text{balk}; q_1, q_2)$ is zero.

If the customer joins with probability $p$, the expected utility becomes
\[
U_i(p; q_1,q_2) = p\cdot U_i(\mathrm{join}; q_1,q_2) + (1-p)\cdot 0 = p\cdot U_i(\mathrm{join}; q_1,q_2).
\]

This linear relationship between $p$ and expected utility means that the customer's optimal decision will depend only on the sign of $U_i(\mathrm{join}; q_1,q_2)$. For notational simplicity, throughout the remainder of this paper, we will use $U_i(q_1,q_2)$ to denote $U_i(\mathrm{join}; q_1, q_2)$.

A customer's decision process mirrors what shoppers intuitively do every day: weigh the value of a product against its price and the hassle of potentially waiting. When many customers adopt joining strategy $q_i$, they collectively influence the system's congestion level, creating a feedback loop where individual decisions affect everyone's experience.

The utility function has several important properties that shape the strategic interactions in this system:
\begin{lemma}[Continuity and Monotonicity of $U_i$]
\label{lem:monotonicity}
For each $i=1,2$, the function 
\[
U_i(q_1,q_2) 
 =  
R_i  -  p_i  -  c_i\mathbb{E}\left[W_i(q_1\Lambda_1, q_2\Lambda_2)\right]
\]
is continuous in $(q_1,q_2)$ over $[0,1]^2$. Moreover, $U_i$ is strictly decreasing in both $q_1$ and $q_2$.
\end{lemma}

The proof is provided in Appendix \ref{app:lem-monotonicity}. It follows by standard calculus arguments, examining partial derivatives of the utility function in both the $S_i=0$ and $S_i>0$ cases to establish strict monotonicity.

Then, we define the best response for type~$i$ as
\[
\mathrm{BR}_i(q_1,q_2)
 = 
\arg \max_{p \in [0,1]}  
p \cdot U_i(q_1,q_2).
\]
From linearity, the best response is:
\[
\mathrm{BR}_i(q_1,q_2) = \begin{cases}
1 & \text{if } U_i(q_1,q_2) > 0\\
0 & \text{if } U_i(q_1,q_2) < 0\\
[0,1] & \text{if } U_i(q_1,q_2) = 0
\end{cases},
\]
where the value $1$ corresponds to the pure action 'always join' and $0$ to 'always balk.' In particular, away from indifference points where $U_i(q_1,q_2)=0$, the unique best response is pure; mixing (i.e., any probability in $[0,1]$) is optimal only when the customer is indifferent between joining and balking. The next lemma studies the structure of these best responses when one type's joining probability is fixed.
\begin{lemma}[Unique Fixed Points in One Variable]\label{lem:unique-fixed}
Fix any $q_2\in[0,1]$.  As a function of $q_1\in[0,1]$, the mapping 
$\mathrm{BR}_1(\cdot,q_2)$ admits exactly one fixed point $q_1^*$ in $[0,1]$.  
Concretely:
\begin{enumerate}
\item 
$\displaystyle q_1^*=0$ if and only if $U_1(0,q_2)  \le  0.$

\item
$\displaystyle q_1^*=1$ if and only if $U_1(1,q_2)  \ge  0.$

\item 
Otherwise, exactly one point $q_1^*\in(0,1)$ satisfies $U_1(q_1^*,q_2)=0$, and that $q_1^*$ is the unique fixed point.
\end{enumerate}

The same holds for type~2 by fixing $q_1$ and considering $\mathrm{BR}_2(q_1,q_2)$ in $q_2$.

\end{lemma}

The proof is provided in Appendix \ref{app:lem-unique-fixed}. It relies on the strict monotonicity established in Lemma~\ref{lem:monotonicity} and applies the Intermediate Value Theorem to guarantee a unique crossing point. 

The strict decrease of $U_i$ in both $q_1$ and $q_2$ established in Lemma~\ref{lem:monotonicity} is a standard avoid-the-crowd property in congestion games: when more customers join, the incentive to join decreases. In one-dimensional queueing models, such avoid-the-crowd behavior directly yields existence and uniqueness of a symmetric equilibrium joining probability; see, for example, the discussion in \citet[][Section~1.1]{HassinRoetGreen2017} and the references therein.
A fixed point $q_1^*$ represents a self-consistent strategy where, given the strategy $q_2$ of type-2 customers, $q_1^*$ is the best response for type-1 customers. This means that when all type-1 customers adopt probability $q_1^*$ and type-2 customers use $q_2$, no individual type-1 customer has an incentive to deviate from $q_1^*$. Despite the set-valued nature of best responses at indifference points, the lemma guarantees that such a fixed point is always unique for any fixed strategy of the other customer type.

\subsection{Nash Equilibrium Characterization}

\subsubsection{Classification}

In the previous analysis, the best-response functions $\mathrm{BR}_1$ and $\mathrm{BR}_2$ have been studied for each customer type separately. When customers of both types make their decisions simultaneously, we search for stable strategies where no individual would benefit by unilaterally changing their decision. This concept of Nash equilibrium (NE) characterizes the likely long-term behavior in the system.

While best responses may be set-valued at indifference points (where $U_i=0$), in equilibrium each customer type adopts a specific joining probability. A pair $(q_1^*,q_2^*)\in[0,1]^2$ is a Nash equilibrium if each probability is consistent with the best response to the other:
\[
q_1^*
 \in 
\mathrm{BR}_1\left(q_1^*,q_2^*\right),
\quad
q_2^*
 \in 
\mathrm{BR}_2\left(q_1^*,q_2^*\right).
\]
Equivalently, if each best-response is single-valued at $\left(q_1^*,q_2^*\right)$, it is possible (by Lemma~\ref{lem:unique-fixed}) to write
\[
\left(q_1^*,q_2^*\right)
 = 
\left(\mathrm{BR}_1(q_1^*,q_2^*), \mathrm{BR}_2(q_1^*,q_2^*)\right),
\]
meaning $\left(q_1^*,q_2^*\right)$ is the fixed point of $\left(\mathrm{BR}_1,\mathrm{BR}_2\right)$.

These equations capture the self-enforcing nature of equilibrium --- given that all customers join with probability $(q_1^*,q_2^*)$, a tagged type-1 customer finds $q_1^*$ optimal, and a tagged type-2 customer finds $q_2^*$ optimal.

Because, by Lemma~\ref{lem:unique-fixed}, each equilibrium joining probability for type~$i$ must be either $0$, $1$, or a unique interior value in $(0,1)$, there are nine possible ways to form $(q_1^*, q_2^*)$.
The proposition below compiles these nine patterns and summarizes each one's utility-based conditions.

\begin{proposition}[Nash Equilibrium Classification]\label{prop:NE-summary}
Any Nash equilibrium $\left(q_1^*,q_2^*\right)$ must lie in exactly one of the nine cases shown in Table~\ref{tab:NE_cases}. Each case corresponds to a fixed point of $(\mathrm{BR}_1,\mathrm{BR}_2)$, with the relevant utility inequalities (or equalities) determining when that fixed point occurs.
\end{proposition}

\begin{table}[ht]
\centering
\caption{Possible Equilibrium Conditions}
\label{tab:NE_cases}
\small  
\begin{tabular}{|p{3.6cm}|p{4.2cm}|p{3.6cm}|}
\hline
\textbf{(0,0) - NE:} & 
\textbf{(0, $\boldsymbol{q_2^*}$), $q_2^*\in(0,1)$ - NE:} & 
\textbf{(0,1) - NE:} \\ 
$U_1(0,0)\le 0$ & 
$U_1(0,q_2^*) \le 0$ & 
$U_1(0,1) \le 0$ \\
$U_2(0,0)\le 0$ & 
$U_2(0,q_2^*) = 0$ & 
$U_2(0,1) \ge 0$ \\
\hline
\textbf{($\boldsymbol{q_1^*}$,0), $q_1^*\in(0,1)$ - NE:} & 
\textbf{($\boldsymbol{q_1^*}$,$\boldsymbol{q_2^*}$), $q_1^*,q_2^*\in(0,1)$ - NE:} & 
\textbf{($\boldsymbol{q_1^*}$,1), $q_1^*\in(0,1)$ - NE:} \\ 
$U_1(q_1^*,0) = 0$ & 
$U_1(q_1^*,q_2^*) = 0$ & 
$U_1(q_1^*,1) = 0$ \\
$U_2(q_1^*,0) \le 0$ & 
$U_2(q_1^*,q_2^*) = 0$ & 
$U_2(q_1^*,1) \ge 0$ \\
\hline
\textbf{(1,0) - NE:} & 
\textbf{(1,$\boldsymbol{q_2^*}$), $q_2^*\in(0,1)$ - NE:} & 
\textbf{(1,1) - NE:} \\ 
$U_1(1,0) \ge 0$ & 
$U_1(1,q_2^*) \ge 0$ & 
$U_1(1,1) \ge 0$ \\
$U_2(1,0) \le 0$ & 
$U_2(1,q_2^*) = 0$ & 
$U_2(1,1) \ge 0$ \\
\hline
\end{tabular}
\end{table}

The proof is provided in Appendix \ref{app:ne-classification}. It proceeds by systematic applying Lemma~\ref{lem:unique-fixed} to each type, showing that the fixed points of the best response must satisfy certain utility conditions.

\subsubsection{Existence of a Nash Equilibrium}
Before investigating the uniqueness properties of equilibria under different inventory scenarios, we first establish that at least one equilibrium always exists in our setting, ensuring the model's predictions are well-defined.
\begin{proposition}[Existence of Nash Equilibrium]
\label{prop:existenceNash}
A Nash equilibrium $(q_1^*, q_2^*)$ exists in this model.
\end{proposition}
\begin{proof}[Proof sketch]
The existence of a Nash equilibrium follows from Brouwer's Fixed Point Theorem. Although the best-response mappings $\mathrm{BR}_1$ and $\mathrm{BR}_2$ are set-valued correspondences, Lemma~\ref{lem:unique-fixed} allows us to construct a continuous function $f(q_1,q_2) = (f_1(q_2),f_2(q_1))$ mapping $[0,1]^2$ to itself, where $f_1(q_2)$ is the unique fixed point of $\mathrm{BR}_1(q_1,q_2)$ for fixed $q_2$, and similarly for $f_2(q_1)$.

The continuity of $f_1$ and $f_2$ follows from the continuity of utility functions and the uniqueness of fixed points. By Brouwer's theorem, $f$ has a fixed point $(q_1^*,q_2^*)$ where $q_1^* = f_1(q_2^*)$ and $q_2^* = f_2(q_1^*)$, which by construction satisfies the Nash equilibrium conditions. See Appendix \ref{app:ne-existence} for the complete proof.
\end{proof}

The challenge encountered here is closely related to that in \citet{HassinRoetGreen2017}, who study a queueing game in which customers choose among join, inspect, and balk. In their model, as in ours, individual utilities satisfy an avoid-the-crowd property in each dimension, but the symmetric strategy is two-dimensional and the standard one-dimensional avoid-the-crowd uniqueness arguments do not apply. They resolve this by analyzing the geometry of an expected-utility set in two dimensions. Both models illustrate how avoid-the-crowd analysis extends to multidimensional strategy spaces.
It is natural, after establishing the existence of a Nash equilibrium, to ask whether that equilibrium is unique.  
This question is easier to address by considering different scenarios for the target inventory levels $S_i$.

\subsubsection{Zero Inventory Scenario: Continuum or Unique Equilibrium}

When both product types maintain zero inventory targets, the customer equilibrium structure depends critically on the relationship between cost-to-margin ratios and service capacity. This serves as the baseline case for understanding more complex inventory policies.
\begin{proposition}
\label{prop:S0-S0}
Let $S_1 = S_2 = 0$. The Nash equilibrium for the joining probabilities $(q_1^*, q_2^*)$ can be described as follows:

\begin{enumerate}
\item If
\[
\frac{c_1}{R_1 - p_1} = \frac{c_2}{R_2 - p_2}
\quad\text{and}\quad
\frac{c_1}{R_1 - p_1} < \mu < \frac{c_1}{R_1 - p_1} + \Lambda_1+ \Lambda_2,
\]
there is a continuum of Nash equilibria satisfying
\[
q_1^*\Lambda_1  + q_2^*\Lambda_2  = \mu-\frac{c_1}{R_1 - p_1},
\quad
0 \le q_1^*,q_2^* \le 1.
\]

\item Otherwise, there is a unique Nash equilibrium that can be explicitly determined from Table~\ref{tab:S0_S0_solutions}.
\end{enumerate}
\end{proposition}
\begin{table}[ht]
\centering
\caption{Equilibrium Outcomes for $S_1 = 0, S_2 = 0$}
\label{tab:S0_S0_solutions}
\small
\begin{tabular}{|p{3.8cm}|p{3.9cm}|p{3.9cm}|}
\hline
\textbf{(0,0):} & 
\textbf{(0,$\boldsymbol{q_2^*}$):} & 
\textbf{(0,1):} \\
& $q_2^* = \frac{1}{\Lambda_2}(\mu - \frac{c_2}{R_2 - p_2})$: & \\
$\mu\le \frac{c_1}{R_1 - p_1}$ & 
$\mu > \frac{c_2}{R_2 - p_2}$ & 
$\mu \ge \frac{c_2}{R_2 - p_2} + \Lambda_2$ \\
$\mu\le \frac{c_2}{R_2 - p_2}$ & 
$\mu < \frac{c_2}{R_2 - p_2} + \Lambda_2$ & 
$\mu \le \frac{c_1}{R_1 - p_1} + \Lambda_2$ \\
& $\frac{c_1}{R_1 - p_1} \ge \frac{c_2}{R_2 - p_2}$ & \\
\hline
\textbf{($\boldsymbol{q_1^*}$,0),} & 
\textbf{($\boldsymbol{q_1^*}$,$\boldsymbol{q_2^*}$),} & 
\textbf{($\boldsymbol{q_1^*}$,1),} \\
$q_1^* = \frac{1}{\Lambda_1}(\mu - \frac{c_1}{R_1 - p_1})$: & 
$q_1^*\Lambda_1 + q_2^*\Lambda_2 = \mu - \frac{c_1}{R_1 - p_1}$, & 
$q_1^* = \frac{1}{\Lambda_1}(\mu - \Lambda_2 - \frac{c_1}{R_1 - p_1})$: \\
& $0 < q_1^*, q_2^* < 1$: & \\
$\mu > \frac{c_1}{R_1 - p_1}$ & 
$\mu > \frac{c_1}{R_1 - p_1}$ & 
$\mu > \frac{c_1}{R_1 - p_1} + \Lambda_2$ \\
$\mu < \frac{c_1}{R_1 - p_1} + \Lambda_1$ & 
$\mu < \frac{c_1}{R_1 - p_1} + \Lambda_1 + \Lambda_2$ & 
$\mu < \frac{c_1}{R_1 - p_1} + \Lambda_1 + \Lambda_2$ \\
$\frac{c_1}{R_1 - p_1} \le \frac{c_2}{R_2 - p_2}$ & 
$\frac{c_1}{R_1 - p_1} = \frac{c_2}{R_2 - p_2}$ & 
$\frac{c_1}{R_1 - p_1} \ge \frac{c_2}{R_2 - p_2}$ \\
\hline
\textbf{(1,0):} & 
\textbf{(1,$\boldsymbol{q_2^*}$),} & 
\textbf{(1,1):} \\
& $q_2^* = \frac{1}{\Lambda_2}(\mu - \Lambda_1 - \frac{c_2}{R_2 - p_2})$: & \\
$\mu \ge \frac{c_1}{R_1 - p_1} + \Lambda_1$ & 
$\mu > \frac{c_2}{R_2 - p_2} + \Lambda_1$ & 
$\mu \ge \frac{c_1}{R_1 - p_1} + \Lambda_1 + \Lambda_2$ \\
$\mu \le \frac{c_2}{R_2 - p_2} + \Lambda_1$ & 
$\mu < \frac{c_2}{R_2 - p_2} + \Lambda_1 + \Lambda_2$ & 
$\mu \ge \frac{c_2}{R_2 - p_2} + \Lambda_1 + \Lambda_2$ \\
& $\frac{c_1}{R_1 - p_1} \le \frac{c_2}{R_2 - p_2}$ & \\
\hline
\end{tabular}
\end{table}

\begin{proof}[Proof sketch]
The key to this analysis is examining how the utility functions simplify when $S_1 = S_2 = 0$:
\[
U_i(q_1, q_2)
=
R_i - p_i
-
c_i \frac{1}{\mu - q_1\Lambda_1 - q_2\Lambda_2}.
\]

This enables explicit computation of all possible equilibrium configurations. When $\frac{c_1}{R_1 - p_1} = \frac{c_2}{R_2 - p_2}$, the condition $U_1=U_2=0$ yields a line segment of equilibria. In all other cases, a unique equilibrium exists. The complete case analysis appears in Appendix \ref{app:s0-s0}.
\end{proof}
This result reveals that even for zero inventories scenario the system can support multiple equilibria only when both customer types have identical patience ratios (the ratio of waiting cost to potential surplus, $\frac{c_i}{R_i - p_i}$). In this special case, customers of both types respond identically to waiting time relative to their product's value, creating a line of equilibrium points that represents different ways to allocate the same total service utilization. Each point on this line maintains the critical balance where total effective arrivals equal $\mu$ minus the patience threshold.

When patience ratios differ, however, a unique equilibrium emerges that typically favors the more patient customer type (lower $\frac{c_i}{R_i - p_i}$ ratio). This uniqueness highlights how even small differences in customer waiting sensitivity create predictable service patterns without requiring external coordination.

\subsubsection{Positive Inventory Scenario: Unique Equilibrium}

When at least one product maintains a positive inventory level, the strategic landscape changes significantly compared to the zero inventory scenario.

\begin{proposition}
\label{prop:positive-inventory}
When at least one product maintains a positive inventory level:
\begin{enumerate}
\item If $S_i > 0$ for some product $i$, then $q_i^* \in (0,1]$ in any Nash equilibrium (i.e., complete balking is impossible).
\item If $S_1 > 0$ and $S_2 > 0$, there exists a unique Nash equilibrium for the customers' joining probabilities $(q_1^*, q_2^*)$.
\item If exactly one target inventory level is positive (either $S_1 > 0, S_2 = 0$ or $S_1 = 0, S_2 > 0$), there still exists a unique Nash equilibrium.
\end{enumerate}
\end{proposition}

Unlike in the zero inventory case, where we derived explicit expressions for equilibrium strategies, the equilibria in positive-inventory scenarios must generally be found numerically. The utility functions become more complex, involving higher powers of the joining probabilities, which precludes simple closed-form solutions.

\begin{proof}[Proof sketch]
The proof proceeds in three steps:

\medskip
\noindent
\textit{Step 1: Proving that $S_i > 0$ implies $q_i > 0$.}
When $S_i > 0$, a customer who arrives to find the system empty of other type-$i$ customers would receive immediate service with probability 1, yielding utility $R_i - p_i > 0$. This makes joining strictly preferable to balking, ruling out $q_i = 0$ as an equilibrium.

\medskip
\noindent
\textit{Step 2: Uniqueness when both $S_1, S_2 > 0$.}
For interior equilibria where both types partially join $(q_1,q_2) \in (0,1)^2$, we analyze the implicit functions $q_1(q_2)$ and $q_2(q_1)$ defined by $U_1(q_1,q_2) = 0$ and $U_2(q_1,q_2) = 0$. The composition of these functions creates a contractive mapping with derivative strictly between 0 and 1, ensuring at most one fixed point.

For boundary equilibria where at least one type fully joins (e.g., $q_1 = 1$), we show that no such equilibrium can coexist with any other equilibrium by analyzing the behavior best responses at the boundaries.

\medskip
\noindent
\textit{Step 3: Uniqueness with asymmetric inventory.}
When exactly one product has positive inventory, similar monotonicity arguments apply with appropriate modifications to handle the possibility of $q_i = 0$ for the zero inventory product.

The complete proof with detailed analysis appears in Appendix \ref{app:positive-inventory}.
\end{proof}

The intuition behind this uniqueness result lies in how positive inventory transforms the mathematical structure of customer responses. When $S_i > 0$, the utility function involves higher powers of joining probabilities through terms like $(q_i\Lambda_i)^{S_i}$, creating a more complex but well-behaved relationship between strategies. This structure ensures that customers' best responses interact in a way that admits only one fixed point.  Unlike the zero inventory case where multiple equilibria can emerge under symmetric patience ratios, the positive-inventory scenario always yields a unique, predictable outcome regardless of whether customer types have similar or different preferences.

\section{Producer's Decision Problem}\label{sec:producer}
Having characterized how customers strategically respond to inventory levels through Nash equilibria, we now examine the producer's perspective. The interaction between producer and customers naturally forms a Stackelberg game --- a sequential decision model where one player (the leader) commits to a strategy before others (the followers) respond. In our context, the producer acts as the leader by setting inventory targets, while customers follow by determining their joining strategies. Unlike simultaneous-move games, the leader anticipates the followers' responses and incorporates this knowledge into its decision, making this a particularly suitable framework for analyzing producer-customer interactions where commitment to inventory levels precedes customer decisions.

\subsection{Profit Function}

The producer's profit can be written as
\[
\Pi(S) = p_1\Lambda_1q_1^{(S_1,S_2)} + p_2\Lambda_2q_2^{(S_1,S_2)} - h_1\mathbb{E}[I_1(S)] - h_2\mathbb{E}[I_2(S)],
\]
where $S = (S_1, S_2)$ denotes the target inventory vector, $p_i$ is the price charged for product type $i$, and $\Lambda_i q_i^{(S_1,S_2)}$ is the effective arrival rate of type-$i$ customers given their equilibrium joining probability $q_i^{(S_1,S_2)}$ at inventory $S$. The term $h_i$ is the holding cost per unit of inventory for product $i$, while $\mathbb{E}[I_i(S)]$ is the expected inventory level derived in Section~\ref{sec:inv_lev}:

\[\mathbb{E}[I_i(S)] = S_i - \frac{\Lambda_i q_i^{(S_1,S_2)}}{\mu - \Lambda_1 q_1^{(S_1,S_2)} - \Lambda_2 q_2^{(S_1,S_2)}} \left( 1 - \left( \frac{\Lambda_i q_i^{(S_1,S_2)}}{\mu - \Lambda_j q_j^{(S_1,S_2)}} \right)^{S_i} \right), \quad j\neq i.\]

The Stackelberg game structure appears in how the profit function is evaluated. The producer first selects inventory levels $(S_1, S_2)$, then customers determine their joining probabilities by reaching a Nash equilibrium among themselves, resulting in $(q_1^{(S_1,S_2)}, q_2^{(S_1,S_2)})$. Only after this sequential process is the producer's profit actually realized. The notation $q_i^{(S_1,S_2)}$ emphasizes that these probabilities are not independent variables but functions determined by the inventory choice.

This sequential decision process creates a fundamental tradeoff: higher inventory levels improve customer experience and increase joining probabilities $q_i^{(S_1,S_2)}$ (thus boosting revenue), but also raise holding costs. The optimization challenge is intensified because the equilibrium joining probabilities $q_i^{(S_1,S_2)}$ generally lack closed-form expressions except in special cases like $S_1=S_2=0$. Instead, they must be computed numerically for each inventory pair through the methods described in Section~\ref{sec:customer}, making the producer's optimization problem analytically intricate.

Unless $(S_1, S_2) = (0, 0)$ and the ratio conditions specified in Proposition \ref{prop:S0-S0} hold (which allows multiple equilibria), each inventory pair $(S_1, S_2)$ determines a unique customer equilibrium $(q_1^{(S_1,S_2)}, q_2^{(S_1,S_2)})$. In the exceptional scenario with multiple equilibria, we adopt the convention of selecting the equilibrium that yields the lowest profit to ensure robust performance guarantees.

\subsection{Target Inventory Thresholds}

To solve this Stackelberg game efficiently, we establish natural bounds on the producer's inventory decisions.

\begin{proposition}\label{prop:thresholds}
There exist thresholds $\bar{S}_1$ and $\bar{S}_2$ such that for all $S_1 \ge \bar{S}_1$ and $S_2 \ge \bar{S}_2$, customers of both types always join the system 
\[
\left(q_1^{(S_1,S_2)} = q_2^{(S_1,S_2)} = 1\right),
\]
and the producer's profit does not exceed
\[
\Pi(S) \le \Pi\left(\bar{S}_1, \bar{S}_2\right)
\quad
\text{for all } S_1 \ge \bar{S}_1, \ S_2 \ge \bar{S}_2.
\]
Moreover, these thresholds can be chosen as
\begin{equation}\label{eq:S_bar_D}
\bar{S}_1 
= 
\min\{S_1 \mid U_1^{S_1}(1,1) \ge 0\},
\quad
\bar{S}_2 
= 
\min\{S_2 \mid U_2^{S_2}(1,1) \ge 0\}.
\end{equation}
\end{proposition}
The proof appears in Appendix~\ref{app:thresholds}.
Once inventory exceeds these thresholds, adding more units increases costs without attracting more customers. Once sufficient inventory is available to make all potential customers willing to join, further increases merely raise costs without attracting new demand.

In the context of the Stackelberg game, this result identifies bounds on the leader's rational strategy space. The first-mover advantage of the producer is effectively limited --- beyond the thresholds $\bar{S}_1$ and $\bar{S}_2$, the followers' responses remain unchanged at full participation, rendering additional inventory strategically irrelevant. This property significantly simplifies the optimization problem by restricting the search to a finite set of candidate inventory levels.

\subsection{Profit-Maximizing Policy}
The producer's strategy involves determining the optimal inventory pair $(S_1^*, S_2^*)$ that maximizes profit, considering how customers will respond to each possible inventory configuration. By Proposition~\ref{prop:thresholds} , searching within $S_1 \in [0, \bar{S}_1]$ and $S_2 \in [0, \bar{S}_2]$ is sufficient for finding the optimal target inventory levels.

The complete solution procedure involves:

\begin{enumerate}
	\item Calculate the upper bounds $\bar{S}_1$ and $\bar{S}_2$ using \eqref{eq:S_bar_D}.
	\item For each inventory pair $(s_1, s_2)$ within these bounds:
	\begin{enumerate}
    	\item Determine the resulting Nash equilibrium joining probabilities $(q_1^{(S_1,S_2)}, q_2^{(S_1,S_2)})$.
   		\item Evaluate the producer's profit under these equilibrium responses.
	\end{enumerate}
	\item Select the inventory pair that yields the highest profit.
\end{enumerate} 

\section{Social Planner's Problem}\label{sec:social}
The previous two sections examined decentralized decision-making where the producer optimizes profit while customers independently maximize their individual utilities. We now consider a centralized perspective --- how would a social planner, concerned with maximizing overall welfare, coordinate inventory levels and customer participation rates?

Because the analysis now considers overall (social) performance, 
individual customers' joining strategies are not modeled directly. 
Instead, the fraction of arrivals of type~$i$ that actually join 
is represented by the effective arrival rate $\lambar{i}$. 
Hence, $\lambar{i}$ replaces $\Lambda_i q_i$ and becomes the decision variable indicating how many potential arrivals 
the planner decides to serve. The planner also chooses the target inventory levels  $(S_1, S_2)$.

\subsection{Social Welfare Function}
The social welfare function is:

\begin{equation*}
\SWe(\lambar{1}, \lambar{2}, S_1, S_2) = \lambar{1} R_1 + \lambar{2} R_2 - C(\lambar{1}, \lambar{2}, S_1, S_2),
\end{equation*}
where $\lambar{i} R_i$ is the total reward obtained from type-$i$ customers and the corresponding total cost function $C(\lambar{1}, \lambar{2}, S_1, S_2)$ represents the overall costs incurred by both products:
\[
    C(\lambar{1}, \lambar{2}, S_1, S_2) = \sum_{i=1}^2 \left( h_i \mathbb{E}[I_i(\lambar{1}, \lambar{2})] + c_i \lambar{i} \mathbb{E}[W_i(\lambar{1}, \lambar{2})] \right).
\]

In the expression above, the first term $h_i \mathbb{E}[I_i(\lambar{1}, \lambar{2})]$ represents the total holding cost for product~$i$, calculated by multiplying the per-unit holding cost $h_i$ by the expected inventory level derived in Section \ref{sec:inv_lev}. The second term $c_i \lambar{i} \mathbb{E}[W_i(\lambar{1}, \lambar{2})]$ captures the total waiting cost incurred by type-$i$ customers, combining the per-customer waiting cost rate $c_i$ with the effective arrival rate $\lambar{i}$ and the expected waiting time from Proposition~\ref{prop:waiting_time}.

When these components are combined for both products, we obtain the complete cost function:
\begin{align}\label{eq:total_cost}
    C(\lambar{1}, \lambar{2}, S_1, S_2) = \sum_{i=1}^2 \Bigg[ h_i S_i - h_i  \frac{\lambar{i}}{\mu - \lambar{1} - \lambar{2}} &\left( 1 - \left( \frac{\lambar{i}}{\mu - \lambar{j}} \right)^{S_i} \right) \nonumber\\
    &+ c_i  \frac{\lambar{i}}{\mu - \lambar{1} - \lambar{2}} \left( \frac{\lambar{i}}{\mu - \lambar{j}} \right)^{S_i} \Bigg],
\end{align}
and the social welfare function:
\begin{align}\label{eq:social_welfare}
\SWe(\lambar{1}, \lambar{2}, S_1, S_2) = 
\sum_{i=1}^2 \Bigg[ \lambar{i} R_i &- h_i S_i + h_i  \frac{\lambar{i}}{\mu - \lambar{1} - \lambar{2}} \left( 1 - \left( \frac{\lambar{i}}{\mu - \lambar{j}} \right)^{S_i} \right) \nonumber \\
&- c_i  \frac{\lambar{i}}{\mu - \lambar{1} - \lambar{2}} \left( \frac{\lambar{i}}{\mu - \lambar{j}} \right)^{S_i} \Bigg].
\end{align}

The social welfare function captures economic efficiency by measuring the total value generated in the system rather than just monetary transfers. While prices appear in the profit function, they cancel out in the welfare calculation since they represent transfers between customers and the producer. Instead, social welfare accounts for the full customer reward ($R_i$), holding costs, and waiting costs --- regardless of which party experiences them. Consequently, the welfare-maximizing solution may allocate capacity differently than what emerges under pure profit maximization.

The optimization of $\SWe(\lambar{1}, \lambar{2}, S_1, S_2)$ follows a hierarchical approach. First, we derive optimal inventory levels $S_1^*(\lambar{1}, \lambar{2})$ and $S_2^*(\lambar{1}, \lambar{2})$ for any given arrival rates, establishing upper bounds $\bar{S}_1$ and $\bar{S}_2$. Since these optimal inventory levels change discretely with arrival rates, they partition the feasible region into subdomains with constant inventory values. We then maximize welfare within each subdomain with respect to $\lambar{1}$ and $\lambar{2}$, and identify the global optimum by comparing these local maxima.

\subsection{Optimal Target Inventory Levels}
From a social planner's perspective, the ideal inventory level for each product balances three factors: holding costs, customer waiting costs, and the probability of immediate service. This tradeoff has an elegant mathematical characterization.

\begin{proposition}\label{prop:optimal_Si}
 For any positive fixed $\lambar{1}, \lambar{2}$ the optimal target inventory levels $(S_1^*, S_2^*)$  minimizing the total cost \eqref{eq:total_cost} are given by
\[
S_i^*(\lambar{1}, \lambar{2}) =  \left\lceil \frac{\log\left( \frac{h_i}{h_i + c_i} \right)}{\log\left(\frac{\lambar{i}}{\mu-\lambar{j}}\right)} \right\rceil-1 \quad \text{for } i = 1, 2.
\]\end{proposition}
\begin{proof}[Proof sketch]
The total cost can be decomposed as $C(S_1, S_2) = C_1(S_1) + C_2(S_2)$, allowing separate optimization for each product. Examining the cost difference $\Delta C_i(S_i) = C_i(S_i+1) - C_i(S_i)$, we identify the threshold where increasing $S_i$ no longer reduces cost. The complete derivation appears in Appendix \ref{app:optimal-si}.
\end{proof}

The formula's structure shows several properties of optimal inventory management. The ratio $\frac{h_i}{h_i + c_i}$ in the logarithm term creates a natural balance point between holding and waiting costs. As this ratio approaches zero (when $c_i \gg h_i$), optimal inventory increases to reduce waiting time. The dependence on $\frac{\lambar{i}}{\mu-\lambar{j}}$ reflects how congestion affects inventory decisions --- as the effective service capacity $\mu-\lambar{j}$ decreases due to the other product's demand, higher inventory becomes necessary to maintain service levels. The ceiling function ensures we obtain a practical integer solution.

Finally, recall that $\Lambda_1+\Lambda_2<\mu$. This implies 
\begin{equation}\label{eq:S_C_bound}
S_i <  \left\lceil \frac{\log\left( \frac{h_i}{h_i + c_i} \right)}{\log\left(\frac{\Lambda_i}{\mu}\right)} \right\rceil \quad \text{for } i = 1, 2
\end{equation}

for any feasible $\lambar{1}, \lambar{2}$.

\subsection{Welfare Optimization Approach}

Having established the optimal inventory levels for each product, we now need to determine the effective arrival rates $\lambar{1}$ and $\lambar{2}$ that maximize social welfare. This optimization is challenging for two reasons: $S_i^*(\lambar{1}, \lambar{2})$ changes discretely across the feasible region, creating discontinuities in the objective function, and even within regions of constant inventory levels, the welfare function's complexity prevents analytical solutions.

We approach this problem by partitioning the feasible region 
\[
\mathcal{D} = \left\{ (\lambar{1}, \lambar{2}) : 0 \le \lambar{i} \le \Lambda_i,  \lambar{1} + \lambar{2} < \mu \right\}
\]
into subdomains where the optimal inventory levels remain constant. From Proposition~\ref{prop:optimal_Si}, we know that $(S_i^*=s_i)$ occurs when:
\[
  \left(\frac{\lambar{i}}{\mu-\lambar{j}}\right)^{s_i+1} \leq \frac{h_i}{h_i + c_i} < 
  \left(\frac{\lambar{i}}{\mu-\lambar{j}}\right)^{s_i}
\]
Rearranging these inequalities yields boundary curves in the $\lambar{1}$-$\lambar{2}$ plane. For $s_1,s_2>0$, these boundaries define subdomains:
\begin{align*}
\mathcal{D}_{s_1, s_2} = \Big\{ (\lambar{1}, \lambar{2}) \in \mathcal{D} : \quad & \mu - \frac{\lambar{1}}{\beta_1^{1/s_1}} < \lambar{2} \leq \mu - \frac{\lambar{1}}{\beta_1^{1/(s_1+1)}},\\
&\qquad
\mu - \frac{\lambar{2}}{\beta_2^{1/s_2}} < \lambar{1} \leq \mu - \frac{\lambar{2}}{\beta_2^{1/(s_2+1)}}\Big\},
\end{align*}
where $\beta_i = \frac{h_i}{h_i+c_i}$. For $s_i = 0$, the conditions simplify to $\lambar{j} \leq \mu - \frac{\lambar{i}}{\beta_i}$.

Within each subdomain, the social welfare function is smooth but still analytically complex. We identify critical points using numerical methods to solve the first-order conditions $\partial \SWe/\partial \lambar{1} = 0$ and $\partial \SWe/\partial \lambar{2} = 0$, while also examining boundary solutions. The global optimum is determined by comparing these local maxima across all subdomains.

This approach transforms a discontinuous optimization problem into a collection of smooth subproblems with well-defined boundaries. The computational complexity remains manageable because the upper bounds on inventory levels are typically modest in practical settings. Appendix \ref{app:optimal-joining-rates} provides further mathematical details on optimization techniques.

\subsection{Welfare-Maximizing Policy}

The social planner's optimization yields a coordinated policy $(\lambar{1}^*, \lambar{2}^*, S_1^*, S_2^*)$ that maximizes overall welfare rather than individual objectives. Unlike the decentralized approach where the producer sets inventory levels and customers independently determine joining probabilities, the welfare-maximizing policy directly controls both aspects of the system.

The solution procedure involves:
\begin{enumerate}
	\item Establishing upper bounds $\bar{S}_1$ and $\bar{S}_2$ on inventory levels based on cost ratios $\frac{h_i}{h_i+c_i}$ and system parameters.
	\item For each candidate inventory pair $(s_1, s_2)$ within these bounds:
	\begin{enumerate}
   		\item Identifying the subdomain $\mathcal{D}_{s_1,s_2}$ where $(S_1^*, S_2^*) = (s_1, s_2)$;
   		\item Optimizing effective arrival rates $(\lambar{1}, \lambar{2})$ within this subdomain;
   		\item Computing the resulting social welfare.
   	\end{enumerate}

	\item Selecting the inventory-arrival rate combination that yields the highest welfare.
\end{enumerate}

This centralized approach can generate different outcomes than the decentralized Stackelberg equilibrium. While a profit-maximizing producer might restrict service to maintain high prices, the social planner typically admits more customers and maintains different inventory levels to balance holding costs against waiting costs.

\subsection{Bridging Decentralized and Centralized Solutions}\label{subsec:alignment}

The decentralized Stackelberg outcome from Section~\ref{sec:producer} generally differs from the welfare optimum: both the effective arrival rates $\lambar{i}$ and the inventory levels $S_i$ may be higher or lower than their socially optimal counterparts, depending on parameters (cf.\ Figure~\ref{fig:app-dec-shares} in our numerical experiments, illustrating that the deviation can go in either direction; see also \citet{Oz2015} for the single-product case). We now show that, for any fixed target inventory levels $S_1,S_2$, a simple correction scheme can implement the welfare-maximizing arrival rates.

Fix a target inventory pair $(S_1,S_2)$ and let $(\lambar{1}^{*},\lambar{2}^{*})$ maximize $\SWe(\lambar{1},\lambar{2},S_1,S_2)$ over feasible $(\lambar{1},\lambar{2})$, with $q_i^{*}:=\lambar{i}^{*}/\Lambda_i$. Suppose a regulator applies a type-$i$ admission correction $\tau_i$ (a toll if $\tau_i>0$, a subsidy if $\tau_i<0$), so the utility from joining becomes $U_i(q_1,q_2)-\tau_i$. Since these transfers are redistributed, the welfare function remains unchanged. Choosing
\[
\tau_i \;=\; R_i-p_i - c_i\,\mathbb{E}\!\left[W_i\bigl(\lambar{1}^{*},\lambar{2}^{*}\bigr)\right],\qquad i=1,2,
\]
supports $(q_1^{*},q_2^{*})$ as a Nash equilibrium of the customer game at inventory $(S_1,S_2)$, implementing the planner's target effective arrival rates.

Implementing the full welfare optimum $(\lambar{1}^*,\lambar{2}^*,S_1^*,S_2^*)$ would additionally require inducing the producer to choose $(S_1^*,S_2^*)$; the interdependence between inventory and customer joining behavior, combined with the discrete nature of inventory decisions, makes characterizing such a joint alignment scheme a nontrivial problem that we leave for future work.
\section{Numerical Experiment: Waiting-Cost Asymmetry}\label{sec:numerics}

This section examines system behavior when customer types differ in 
how much waiting they can tolerate. Examples of such asymmetry include business versus leisure travelers in transportation services, emergency versus elective patients in healthcare, and real-time versus batch processing workloads in computing resources. We compare the decentralized solution (profit-maximizing producer anticipating equilibrium behavior) with the centralized solution (welfare-maximizing planner), examining how profit, welfare, inventory decisions, and joining behavior vary across a range of asymmetry levels and potential utilization rates.

\subsection{Measuring Asymmetry}\label{sec:asymmetry-param}

A type-$i$ customer joins when the expected wait satisfies 
$\mathbb{E}[W_i(\cdot)] \le (R_i - p_i)/c_i$. Multiplying by the service rate $\mu$ 
gives a dimensionless quantity
\[
\nu_i := \frac{\mu(R_i - p_i)}{c_i},
\]
representing the maximum tolerable wait in units of mean service time $1/\mu$; 
larger $\nu_i$ means greater delay tolerance.

To isolate asymmetry between customer types, we track the ratio
\[
\kappa := \frac{\nu_1}{\nu_2}.
\]
Here $\kappa = 1$ means both types tolerate delay equally relative to their surplus; 
$\kappa > 1$ means type~2 is more delay-sensitive. 
In our symmetric setup ($R_1 = R_2$, $p_1 = p_2$), this simplifies to $\kappa = c_2/c_1$.

\subsection{Experimental Setup}\label{sec:setup}

We normalize $\mu = 1$ and impose symmetry in all parameters except type-2 waiting cost: $R_1 = R_2 = 10$, $p_1 = p_2 = 5$, $h_1 = h_2 = 0.4$, and $c_1 = 3$ (Appendix~\ref{app:reduced-h1-figs} reports results for an alternative configuration with $c_1 = 1$ and $h_1 = 0.05$).

We vary $\kappa \in [1, 20]$, from equal waiting tolerance ($\kappa = 1$) to type-2 tolerating only 1/20th as much waiting ($\kappa = 20$), implemented by setting $\nu_2 = \nu_1/\kappa$; under $(R_1, p_1) = (R_2, p_2)$ this simplifies to $c_2 = \kappa c_1$. The second axis is the potential utilization $\rho := (\Lambda_1 + \Lambda_2)/\mu \in [0.65, 0.90]$, with $\Lambda_1 = \Lambda_2 = \rho\mu/2$. Because customers may balk, realized utilization is typically below $\rho$ (see Figure~\ref{fig:app-rho-eff} in Appendix~\ref{app:baseline-figs}).

For each $(\kappa, \rho)$ we compute: (i) the \emph{decentralized} Stackelberg outcome $(S^{\textsc{dec}}_i, q^{\textsc{dec}}_i)$, where the producer maximizes profit anticipating Nash equilibrium joining rates; and (ii) the \emph{centralized} outcome $(S^{\textsc{cen}}_i, \lambda^{\textsc{cen}}_i)$, where a social planner maximizes welfare. For comparability we define $q^{\textsc{cen}}_i := \lambda^{\textsc{cen}}_i/\Lambda_i$ and refer to both $q^{\textsc{dec}}_i$ and $q^{\textsc{cen}}_i$ as joining rates.

Results are displayed as heatmaps over the $(\kappa, \rho)$ plane, computed on a grid with step $0.01$ in $\kappa$ and $0.001$ in $\rho$. Unless stated otherwise, holding costs are symmetric ($h_2/h_1 = 1$). Section~\ref{sec:sensitivity} plots cross-sections at $\rho = 0.9$ for $h_2/h_1 \in \{0.9, 1.0, 1.1\}$, both to examine sensitivity to holding-cost asymmetry and to provide a clearer view of how outcomes vary with $\kappa$ at fixed potential utilization.

\subsection{Profit and Social Welfare}\label{sec:profit-welfare}
\begin{figure}[!htbp]
\centering
\includegraphics[width=0.8\textwidth]{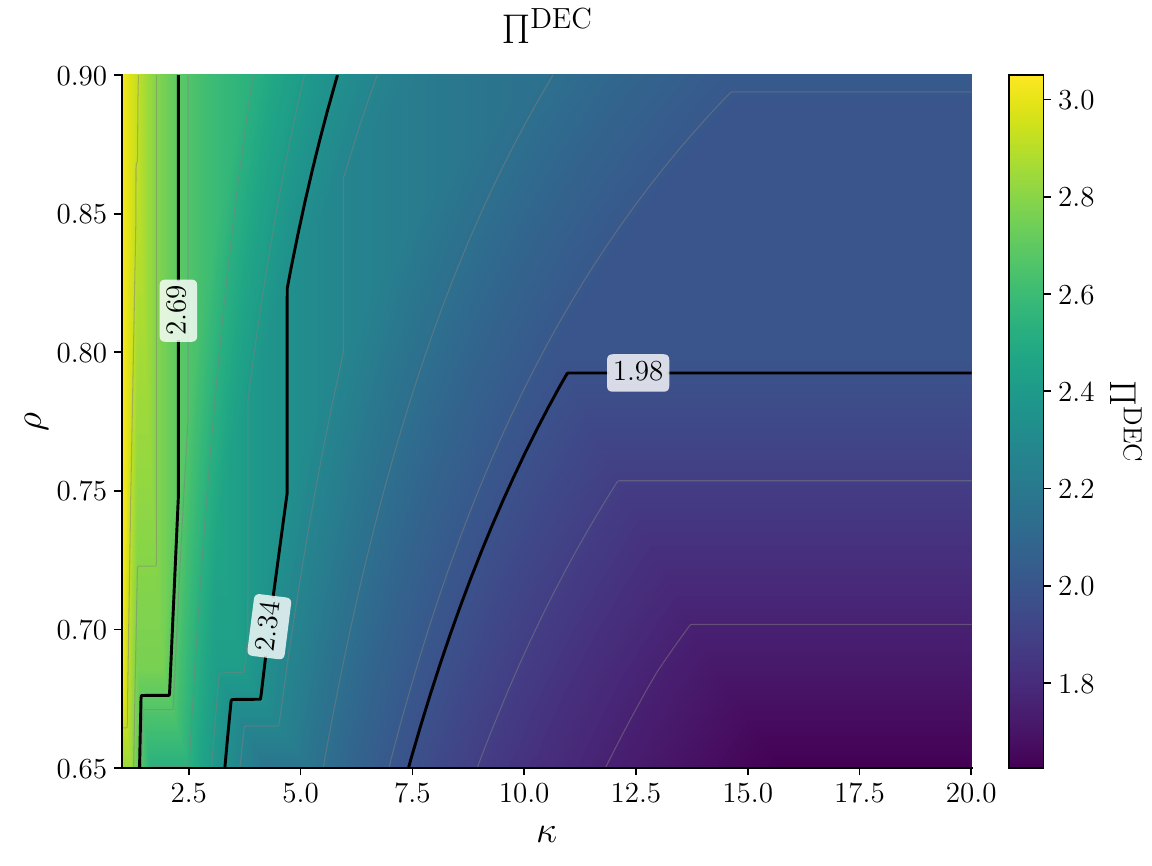}
\caption{Producer's profit $\Pi^{\textsc{dec}}$ under decentralization over the $(\kappa,\rho)$ plane. Thick black contours mark 1.98, 2.34, 2.69 (25\%, 50\%, 75\% of range); thin gray contours provide detail.}
\label{fig:profit-dec}
\end{figure}

\begin{figure}[!htbp]
\centering
\includegraphics[width=0.8\textwidth]{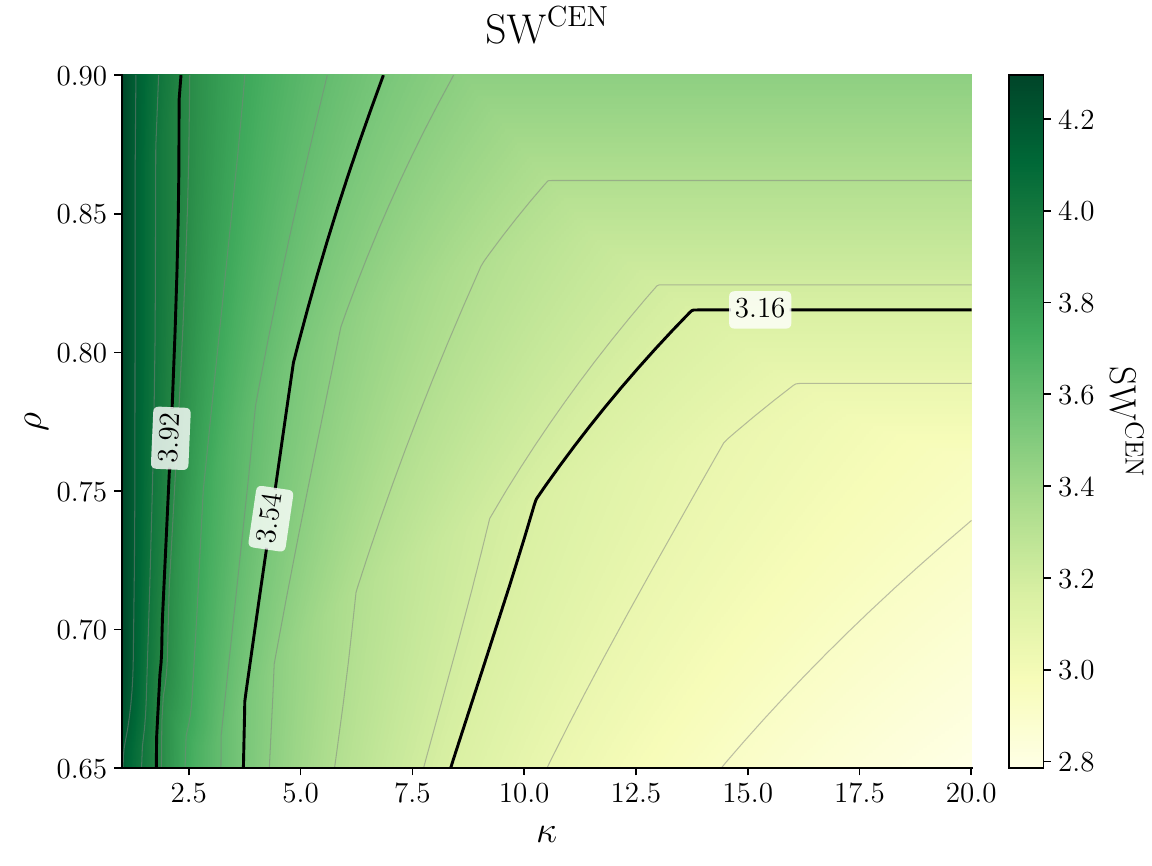}
\caption{Social welfare $\mathrm{SW}^{\textsc{cen}}$ under centralization over the $(\kappa,\rho)$ plane. Thick black contours mark 3.16, 3.54, 3.92 (25\%, 50\%, 75\% of range); thin gray contours provide detail.}
\label{fig:welfare-cen}
\end{figure}

Figures~\ref{fig:profit-dec} and~\ref{fig:welfare-cen} display the producer's profit under the decentralized solution ($\Pi^{\textsc{dec}}$) and the social welfare under the centralized solution ($\mathrm{SW}^{\textsc{cen}}$) across the $(\kappa, \rho)$ plane. Both metrics increase in~$\rho$: although higher~$\rho$ implies greater potential congestion, it also means a larger pool of potential customers ($\Lambda_1 + \Lambda_2 = \rho\mu$), and this demand effect dominates.

In contrast, both metrics decrease in~$\kappa$. As~$\kappa$ rises, type-2 customers tolerate less waiting, forcing the system to either hold more inventory (increasing costs) or accept lower type-2 participation (reducing revenue and rewards). Although both producer and planner adjust~$S_2$ to partially compensate (see Section~\ref{sec:decisions}), this is insufficient to fully offset the decline.

\begin{figure}[!htbp]
\centering
\includegraphics[width=0.8\textwidth]{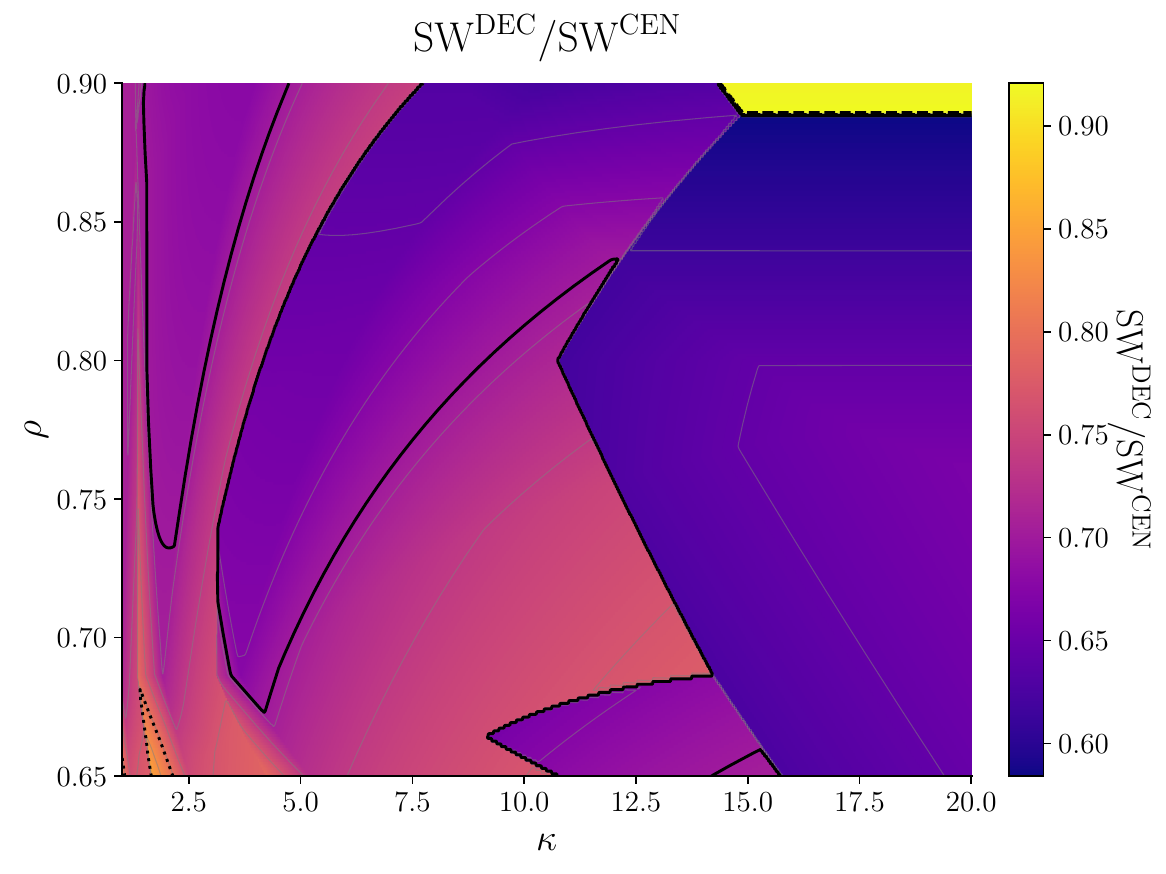}
\caption{Ratio of decentralized to centralized welfare, $\mathrm{SW}^{\textsc{dec}}/\mathrm{SW}^{\textsc{cen}}$ over the $(\kappa,\rho)$ plane. Contours mark $0.70$ (solid), $0.80$ (dotted), and $0.90$ (dashed); thin gray contours provide detail.}
\label{fig:efficiency}
\end{figure}

The social welfare achieved under the decentralized solution, $\mathrm{SW}^{\textsc{dec}}$ (Figure~\ref{fig:app-welfare-dec} in Appendix~\ref{app:baseline-figs}), lies below $\mathrm{SW}^{\textsc{cen}}$ throughout the parameter space. To quantify this welfare loss when all parties optimize individually rather than a planner coordinating for social welfare, we compute the ratio
\[
\frac{\mathrm{SW}^{\textsc{dec}}(\kappa, \rho)}{\mathrm{SW}^{\textsc{cen}}(\kappa, \rho)} \;\in\; (0,1].
\]
Figure~\ref{fig:efficiency} reveals that this ratio does not vary monotonically in either parameter. Plateau near~$0.92$ appears in the upper-right region, where both solutions exclude type-2 entirely ($S_2 = 0$, $q_2 = 0$). Once type-2 is excluded, varying~$\kappa$ no longer affects the system; both solutions reduce to a single-type problem with $S_1^{\textsc{dec}} = 1$, $S_1^{\textsc{cen}} = 2$, and $q_1 = 1$.

\FloatBarrier
\subsection{Target Inventories and Joining Behavior}\label{sec:decisions}

\begin{figure}[!htbp]
\centering
\includegraphics[width=\textwidth]{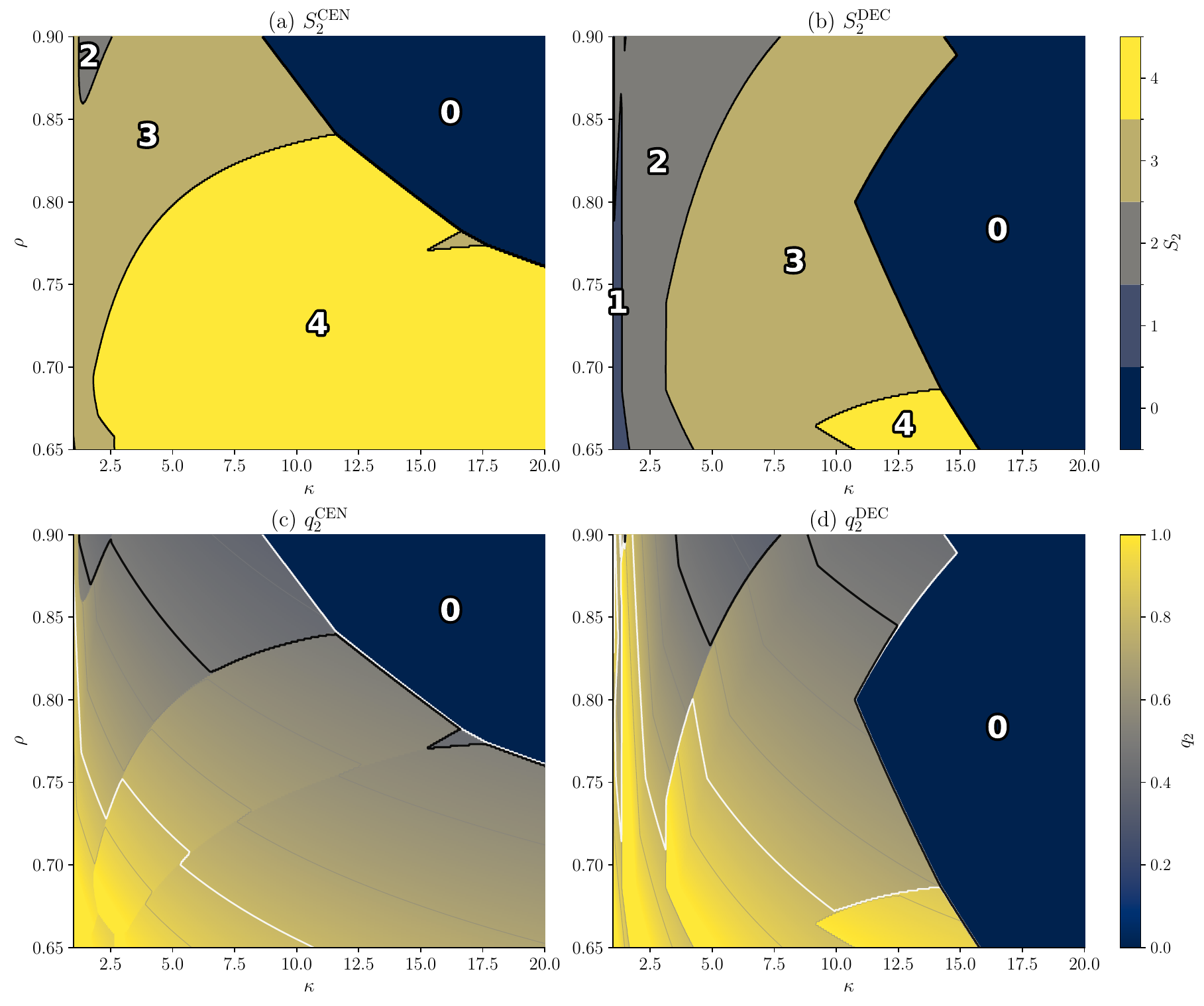}
\caption{Type-2 decisions over the $(\kappa,\rho)$ plane. Top row: inventory $S_2$; contours mark integer transitions. Bottom row: joining probability $q_2$; contours at $0.2$, $0.8$ (thick white) and $0.5$ (thick black); thin gray contours provide detail.}
\label{fig:type2-response}
\end{figure}

Figure~\ref{fig:type2-response} shows type-2 target inventory levels and joining rates. Both the producer and the planner respond to rising~$\kappa$ by increasing~$S_2$, holding more inventory to reduce expected waiting. Eventually, however, holding cost outweighs the benefit and the system gives up on type-2: inventory drops to zero. This compensate-then-drop pattern is the dominant feature of the heatmaps, though the boundary where exclusion occurs differs between producer and planner. At low $\rho$, the decentralized solution exhibits exclusion within our $\kappa$ range, while the centralized solution's exclusion threshold lies beyond $\kappa = 20$.

\begin{figure}[!htbp]
\centering
\includegraphics[width=\textwidth]{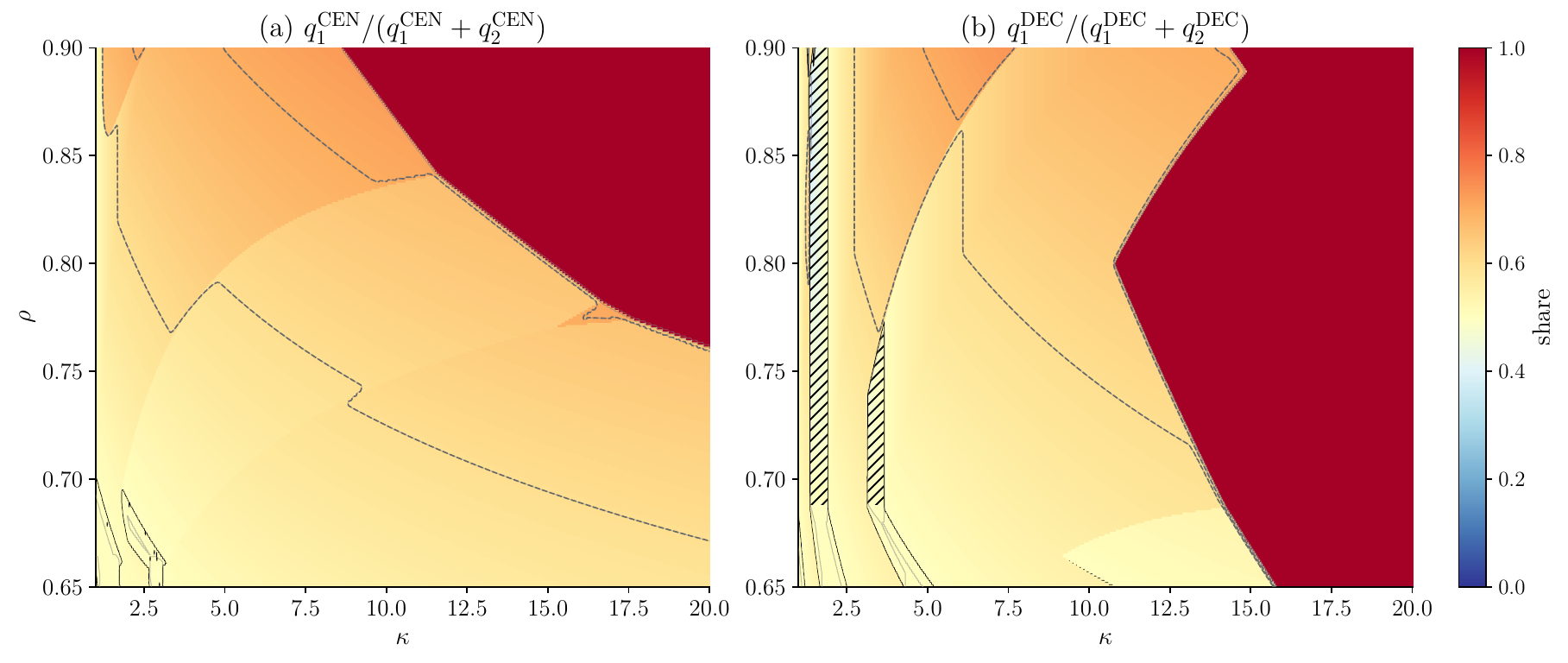}
\caption{Type-1 share $q_1/(q_1+q_2)$ over the $(\kappa,\rho)$ plane. Contours at $0.5$ (thin solid), $0.6$ and $0.7$ (dashed). Hatching indicates type-2 majority (share $< 0.5$).}
\label{fig:q1share}
\end{figure}

Type-1 joining rates stay close to~$1$ throughout most of the parameter space, with $S_1$ adjusting modestly (Figure~\ref{fig:app-type1-response} in Appendix~\ref{app:baseline-figs}). However, high $q_1$ does not guarantee type-1 dominance. Figure~\ref{fig:q1share} plots the type-1 share of total joining. The centralized solution always has type-1 at or above half, but the decentralized equilibrium exhibits a region at low~$\kappa$ where $q_1^{\textsc{dec}} < q_2^{\textsc{dec}}$ --- the type with lower waiting cost participates less. This inversion disappears as~$\kappa$ grows.

\FloatBarrier
\subsection{Holding-Cost Sensitivity}\label{sec:sensitivity}

\begin{figure}[!htbp]
\centering
\includegraphics[width=\textwidth]{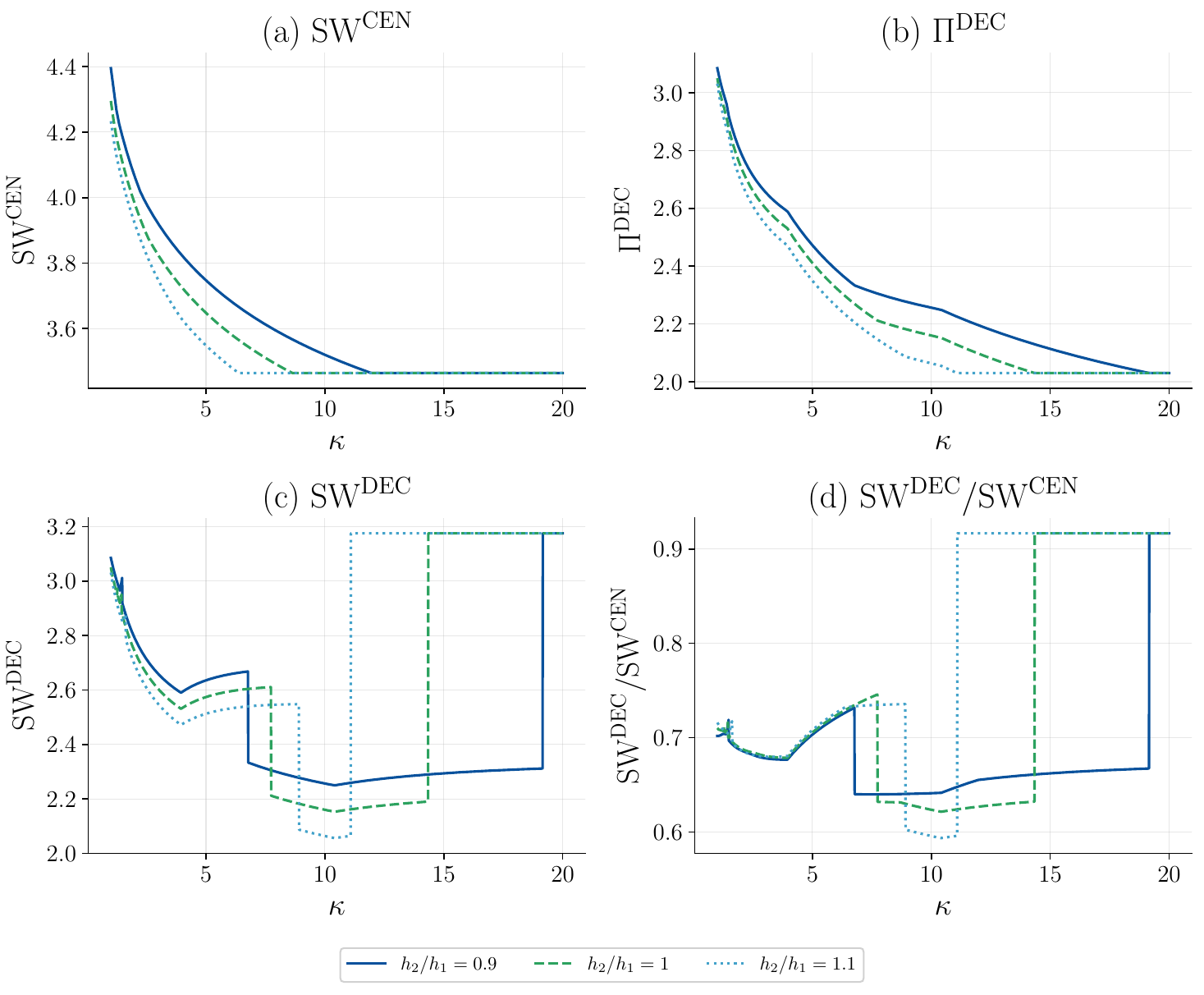}
\caption{Cross-section at $\rho = 0.9$ as functions of $\kappa$ for different $h_2/h_1$ ratios. Top row: $\mathrm{SW}^{\textsc{cen}}$ (left) and $\Pi^{\textsc{dec}}$ (right). Bottom row: $\mathrm{SW}^{\textsc{dec}}$ (left) and $\mathrm{SW}^{\textsc{dec}}/\mathrm{SW}^{\textsc{cen}}$ (right).}
\label{fig:cross-welf-profit}
\end{figure}

\begin{figure}[!htbp]
\centering
\includegraphics[width=\textwidth]{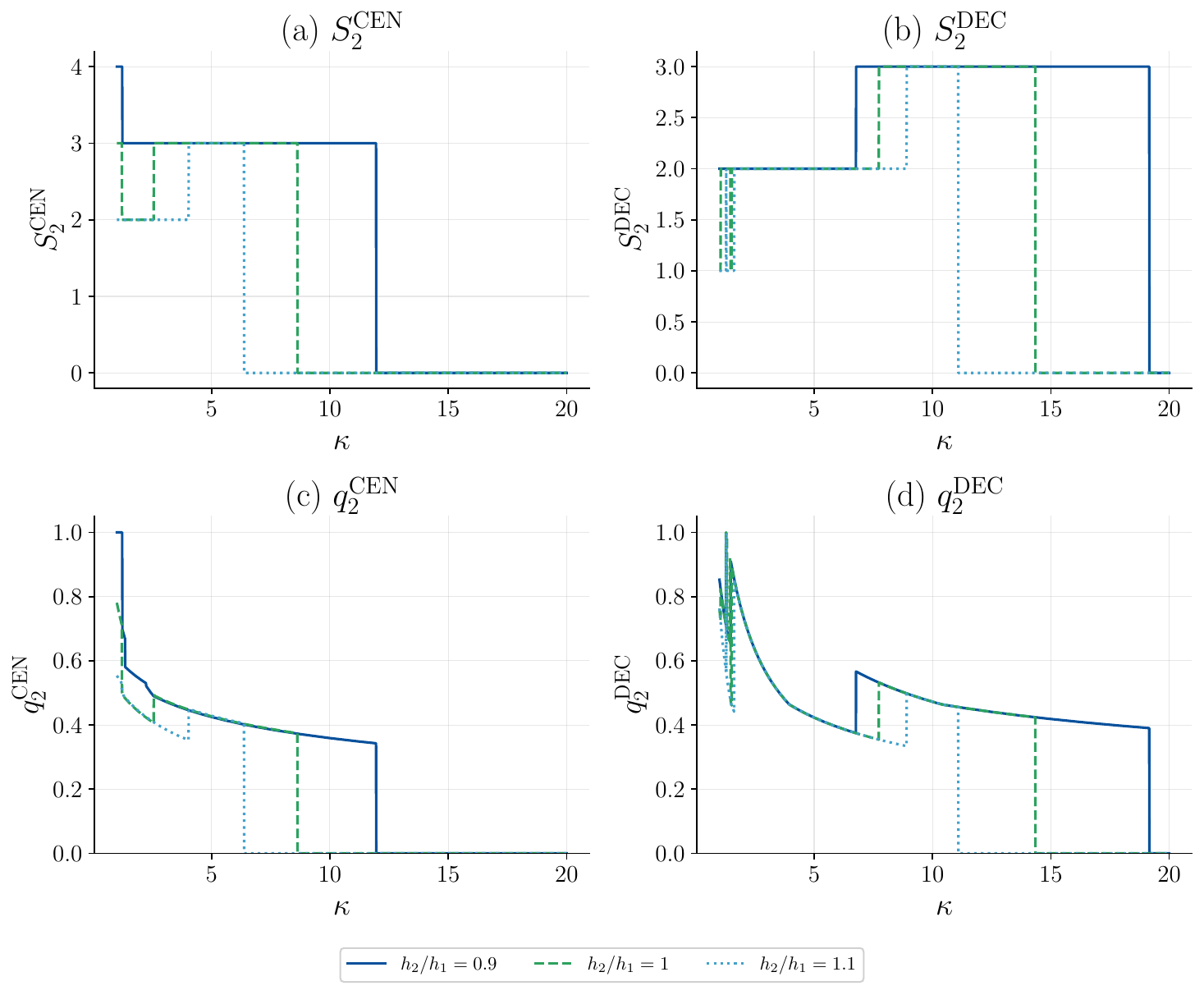}
\caption{Cross-section at $\rho = 0.9$: type-2 decisions as functions of $\kappa$ for different $h_2/h_1$ ratios. Top row: inventory $S_2$. Bottom row: joining probability $q_2$.}
\label{fig:cross-type2}
\end{figure}

The baseline experiment fixes symmetric holding costs ($h_2/h_1 = 1$). We now vary this ratio with $h_1$ fixed and $h_2$ varying. Higher $h_2/h_1$ makes type-2 inventory more expensive to hold, limiting the system's ability to buffer type-2 waiting times as~$\kappa$ rises. Conversely, lower $h_2/h_1$ makes inventory a cheaper tool for accommodating high-$c_2$ customers, extending the range of~$\kappa$ over which type-2 remains viable. Figure~\ref{fig:cross-welf-profit} reports outcome metrics as functions of~$\kappa$ at $\rho = 0.9$; Figure~\ref{fig:cross-type2} reports the associated type-2 decisions.

Changing $h_2/h_1$ primarily shifts the regime boundaries at which type-2 becomes unattractive to serve. Higher $h_2/h_1$ pushes the exclusion cutoff to lower~$\kappa$; lower $h_2/h_1$ expands the region where type-2 is served. The stepwise changes in~$S_2$ produce corresponding kinks and jumps in joining rates and welfare, reflecting transitions between integer inventory regimes. The welfare ratio $\mathrm{SW}^{\textsc{dec}}/\mathrm{SW}^{\textsc{cen}}$ is most sensitive to $h_2/h_1$ near these transitions, where producer and planner cutoffs diverge; away from the boundaries, both solutions react similarly and the ratio is comparatively stable.

Additional results --- including effective utilization, expected waiting times, and component-wise comparisons between centralized and decentralized solutions --- are reported in Appendix~\ref{app:baseline-figs}.

\subsection{Implications}\label{sec:numerics-summary}

Our numerical analysis yields several insights:

\begin{itemize}
\item When one customer type becomes more delay-sensitive, holding additional inventory for that type reduces expected waiting and sustains participation. However, holding costs eventually outweigh the benefit, making it optimal to stop serving the delay-sensitive type entirely. The transition is abrupt: small parameter changes near the boundary can shift the system from serving both types to serving only one.

\item At low asymmetry, decentralized equilibria exhibit a counterintuitive pattern: the customer type with lower waiting cost participates less. This inversion arises because individual joining decisions do not account for system-wide effects. The phenomenon disappears as asymmetry grows.

\item The parameter regions where a profit-maximizing firm stops serving the delay-sensitive type differ from those where a welfare-maximizing planner would do so. This disagreement about exclusion boundaries --- rather than differences in inventory levels or joining rates --- is the primary driver of welfare loss from decentralized decision-making.

\item Because welfare loss concentrates where profit-maximizing and welfare-maximizing policies make different serve-or-exclude decisions, tolls or subsidies would be most effective when the system is near the margin of abandoning a customer type.
\end{itemize}

\FloatBarrier

\section{Conclusion}

This paper analyzed strategic customer behavior in a two-product make-to-stock queueing system where products share production capacity. We established several key results:

First, we derived closed-form expressions for expected waiting times, inventory levels, and backlogs, showing how these measures respond to inventory levels and customer joining strategies. We then proved that Nash equilibria exist and are unique when at least one product has positive inventory, while a continuum of equilibria can arise in special zero-inventory cases.

For the producer, we demonstrated that profit-maximizing inventory levels have finite upper bounds, eliminating the need to search an unbounded parameter space. Similarly, the social planner's optimal inventory policy follows a clear structure based on the ratio of holding to waiting costs, even as effective arrival rates vary.

Our numerical experiments illustrated how inventory buffers waiting-cost asymmetry --- up to a point --- and how profit-maximizing and welfare-maximizing policies can disagree on when to exclude delay-sensitive customers.

The optimization methods described in this paper are implemented in Python and available at \href{https://github.com/ksrvskj/make-to-stock-2product}{github.com/ksrvskj/make-to-stock-2product}.

\newpage

\bibliographystyle{plainnat}
\bibliography{kan-kos-2p-mts1}

@book{Zipkin2000,
  author    = {Zipkin, Paul},
  title     = {Foundations of Inventory Management},
  publisher = {McGraw-Hill},
  address   = {New York},
  year      = {2000}
}

@book{Silver1998,
  author    = {Silver, Edward A. and Pyke, David F. and Peterson, Rein},
  title     = {Inventory Management and Production Planning and Scheduling},
  edition   = {3rd},
  publisher = {Wiley},
  address   = {New York},
  year      = {1998}
}

@book{Graves2000,
  editor    = {Graves, Stephen C. and Rinnooy Kan, A. H. G. and Zipkin, Paul H.},
  title     = {Logistics of Production and Inventory},
  series    = {Handbooks in Operations Research and Management Science},
  volume    = {4},
  publisher = {North-Holland},
  address   = {Amsterdam},
  year      = {1993}
}

@article{Naor1969,
  author    = {Naor, P.},
  title     = {The Regulation of Queue Size by Levying Tolls},
  journal   = {Econometrica},
  volume    = {37},
  number    = {1},
  pages     = {15--24},
  year      = {1969}
}

@article{Edelson1975,
  author    = {Edelson, Noel M. and Hilderbrand, David K.},
  title     = {Congestion Tolls for {P}oisson Queuing Processes},
  journal   = {Econometrica},
  volume    = {43},
  number    = {1},
  pages     = {81--92},
  year      = {1975}
}

@book{HassinHaviv2003,
  author    = {Hassin, Refael and Haviv, Moshe},
  title     = {To Queue or Not to Queue: Equilibrium Behavior in Queueing Systems},
  publisher = {Springer},
  address   = {Boston},
  year      = {2003}
}

@book{Buzacott1993,
  title     = {Stochastic Models of Manufacturing Systems},
  author    = {Buzacott, John A. and Shanthikumar, J. George},
  publisher = {Prentice Hall},
  year      = {1993},
  address   = {Englewood Cliffs, NJ}
}

@article{Caldentey2006,
  author    = {Caldentey, Ren{\'e} and Wein, Lawrence M.},
  year      = {2006},
  month     = {October},
  number	= {5},
  title     = {Revenue Management of a Make-to-Stock Queue},
  volume    = {54},
  pages     = {859--875},
  journal   = {Operations Research}
}

@inproceedings{Oz2015,
  title        = {On a Production/Inventory System with Strategic Customers and Unobservable Inventory Levels},
  author       = {{\"O}z, Can and Karaesmen, Fikri},
  booktitle    = {SMMSO 2015: 10th Conference on Stochastic Models of Manufacturing and Service Operations},
  year         = {2015},
  address     = {Volos, Greece},
  pages        = {161--168},
  publisher    = {University of Thessaly Press}
}

@article{Zhang2019,
  author  = {Zhang, Xuelu and Wang, Jinting},
  year    = {2019},
  month   = {May},
  title   = {Optimal inventory threshold for a dynamic service make-to-stock system with strategic customers},
  volume  = {35},
  pages   = {1103--1123},
  journal = {Applied Stochastic Models in Business and Industry}
}

@article{Arrow1951,
  author  = {Arrow, Kenneth J. and Harris, Theodore and Marschak, Jacob},
  title   = {Optimal Inventory Policy},
  journal = {Econometrica},
  volume  = {19},
  number  = {3},
  pages   = {250--272},
  year    = {1951}
}

@article{Scarf1960,
  author  = {Scarf, Herbert},
  title   = {The Optimality of {$(S, s)$} Policies in the Dynamic Inventory Problem},
  journal = {Mathematical Methods in the Social Sciences},
  year    = {1960},
  pages   = {196--202},
  editor  = {Arrow, K. and Karlin, S. and Suppes, P.},
  publisher = {Stanford University Press}
}

@article{Karlin1958,
  author  = {Karlin, Samuel},
  title   = {Optimal Policy for Dynamic Inventory Process with Stochastic Demands Subject to Seasonal Variations},
  journal = {Journal of the Society for Industrial and Applied Mathematics},
  volume  = {8},
  number  = {4},
  pages   = {611--629},
  year    = {1960}
}

@article{Veinott1965,
  author  = {Veinott, Arthur F.},
  title   = {Optimal Policy for a Multi-Product, Dynamic, Nonstationary Inventory Problem},
  journal = {Management Science},
  volume  = {12},
  number  = {3},
  pages   = {206--222},
  year    = {1965}
}

@article{Federgruen1996,
  author  = {Federgruen, Awi and Katalan, Zvi},
  title   = {The Stochastic Economic Lot Scheduling Problem: Cyclical Base-Stock Policies with Idle Times},
  journal = {Management Science},
  volume  = {42},
  number  = {6},
  pages   = {783--796},
  year    = {1996}
}

@article{Veatch1996,
  author  = {Veatch, Michael H. and Wein, Lawrence M.},
  title   = {Optimal Control of a Two-Station Tandem Production/Inventory System},
  journal = {Operations Research},
  volume  = {42},
  number  = {2},
  pages   = {337--350},
  year    = {1996}
}

@article{Hassin1985,
  author  = {Hassin, Refael},
  title   = {On the Optimality of First-Come Last-Served Queues},
  journal = {Econometrica},
  volume  = {53},
  number  = {1},
  pages   = {201--202},
  year    = {1985}
}

@article{Mendelson1985,
  author  = {Mendelson, Haim},
  title   = {Pricing Computer Services: Queueing Effects},
  journal = {Communications of the {ACM}},
  volume  = {28},
  number  = {3},
  pages   = {312--321},
  year    = {1985}
}

@article{Burnetas2007,
  author  = {Burnetas, Apostolos N. and Economou, Alexandros},
  title   = {Equilibrium Customer Strategies in a Single-Server {M}arkovian Queue with Setup Times},
  journal = {Queueing Systems},
  volume  = {56},
  number  = {3--4},
  pages   = {213--228},
  year    = {2007}
}

@book{Hassin2016,
  author    = {Hassin, Refael},
  title     = {Rational Queueing},
  publisher = {CRC Press},
  year      = {2016},
  address   = {Boca Raton, FL}
}

@article{Ata2006,
  author  = {Ata, Baris and Shneorson, Sergey},
  title   = {Dynamic Control of an M/M/1 Service System with Adjustable Arrival and Service Rates},
  journal = {Management Science},
  volume  = {52},
  number  = {11},
  pages   = {1778--1791},
  year    = {2006}
}

@article{Cachon2009,
  author  = {Cachon, Gerard P. and Swinney, Robert},
  title   = {Purchasing, Pricing, and Quick Response in the Presence of Strategic Consumers},
  journal = {Management Science},
  volume  = {55},
  number  = {3},
  pages   = {497--511},
  year    = {2009}
}

@article{Li2016,
  author  = {Li, Qiang and Guo, Peng and Li, C. L. and Song, J.},
  title   = {Equilibrium Joining Strategies and Optimal Control of a Make-to-Stock Queue},
  journal = {Production and Operations Management},
  volume  = {25},
  number  = {9},
  pages   = {1513--1527},
  year    = {2016}
}

@article{WangCai2023,
  author  = {Wang, Yimin and Cai, Chen},
  title   = {Strategic Behavior and Optimal Inventory Level in a Make-to-Stock Queueing System with Retrial Customers},
  journal = {Axioms},
  volume  = {12},
  number  = {5},
  pages   = {319},
  year    = {2023},
}

@article{Benjaafar2012,
  author  = {Benjaafar, Saif and ElHafsi, Mounir},
  title   = {Production and Inventory Control of a Single Product Assemble-to-Order System with Multiple Customer Classes},
  journal = {Management Science},
  volume  = {52},
  number  = {12},
  pages   = {1896--1912},
  year    = {2006}
}

@article{Cai2020,
  author  = {Cai, Xiaolan and Li, Jian and Chen, Shenglin and Tang, Xiangtong},
  title   = {Joint Pricing and Inventory Control in a Make-to-Stock Queue with Delay-Sensitive Customers},
  journal = {Journal of the Operational Research Society},
  volume  = {72},
  number  = {8},
  pages   = {1841--1856},
  year    = {2020}
}

@article{Abouee2012,
  author  = {Abouee-Mehrizi, Hossein and Balc{\i}o{\u{g}}lu, Ahmet Bora and Baron, Opher},
  title   = {Strategies for a Centralized Single Product Multiclass {M/G/1} Make-to-Stock Queue},
  journal = {Operations Research},
  volume  = {60},
  number  = {4},
  pages   = {803--812},
  year    = {2012}
}

@article{Zhao2011,
  author  = {Zhao, Ning and Lian, Zhe},
  title   = {A Queueing-Inventory System with Two Classes of Customers},
  journal = {International Journal of Production Economics},
  volume  = {129},
  number  = {1},
  pages   = {225--231},
  year    = {2011}
}

@article{Zare2021,
  author  = {Zare, Amir G. and Abouee-Mehrizi, Hossein and Berman, Oded},
  title   = {Two-Echelon Production Inventory Systems with Strategic Customers},
  journal = {Probability in the Engineering and Informational Sciences},
  volume  = {35},
  number  = {2},
  pages   = {258--275},
  year    = {2021}
}

@article{Gong2018,
  author    = {Abhishek, Vineet and Kash, Ian and Key, Peter},
  title     = {Fixed and market pricing for cloud services},
  journal.  = {Proceedings - IEEE INFOCOM},
  year      = {2012},
  pages     = {157--162},
  number	= {1}
}

@article{Hong2017,
  author  = {Lin, Fuhong and Zhou, Xianwei and Huang, Daochao and Song, Wei and Han, Dongsheng},
  title   = {Service Scheduling in Cloud Computing based on Queuing Game Model
},
  journal = {KSII Transactions on Internet and Information Systems},
  volume  = {8},
  number  = {5},
  pages   = {1554--1566},
  year    = {2014}
}

@article{Harris1913,
  author  = {Harris, Ford W.},
  title   = {How Many Parts to Make at Once},
  journal = {Factory, The Magazine of Management},
  volume  = {10},
  number  = {2},
  pages   = {135--136},
  year    = {1913}
}

@article{HassinRoetGreen2017,
  author  = {Hassin, Refael and Roet-Green, Ricky},
  title   = {The Impact of Inspection Cost on Equilibrium, Revenue, and Social Welfare in a Single-Server Queue},
  journal = {Operations Research},
  volume  = {65},
  number  = {3},
  pages   = {804--820},
  year    = {2017}
}

@article{Ha1997,
  author  = {Ha, Albert Y.},
  title   = {Inventory Rationing in a Make-to-Stock Production System with Several Demand Classes and Lost Sales},
  journal = {Management Science},
  volume  = {43},
  number  = {8},
  pages   = {1093--1103},
  year    = {1997}
}

@article{DeVericourt2002,
  author  = {de V{\'e}ricourt, Francis and Karaesmen, Fikri and Dallery, Yves},
  title   = {Optimal Stock Allocation for a Capacitated Supply System},
  journal = {Management Science},
  volume  = {48},
  number  = {11},
  pages   = {1486--1501},
  year    = {2002}
}

\newpage

\begin{appendices}
\section{System Performance Proofs}\label{app:system-analysis}
\subsection{Queue Position Distribution}\label{app:position-distribution}
\begin{proof}[Proof of \eqref{eq:joint-position}]
Without loss of generality, we consider the case where $i=1$; the argument for $i=2$ follows analogously by symmetry.

When a newly arriving type-1 customer finds product~1 out of stock, the system ensures there are exactly $S_1$ type-1 jobs in the production queue designated for replenishment. If a customer joins, they become the newest type-1 backorder. Their unit is supplied by the completion of the oldest replenishment type-1 job currently in the queue (after any earlier backorders). The join decision also triggers a new type-1 job appended at the tail, which becomes the newest replenishment job to maintain the inventory position at $S_1$.

The existing jobs in queue came from independent Poisson arrival streams with rates $\lambar{1}$ and $\lambar{2}$. Each of the current $N(t)=n$ jobs in queue is type-1 with probability $\lambar{1}/(\lambar{1}+\lambar{2})$ or type-2 with probability $\lambar{2}/(\lambar{1}+\lambar{2})$, subject to the condition that at least $S_1$ are type-1. This is the binomial thinning property. Without additional constraints, the probability of exactly $k$ type-1 jobs among $m$ is
\[
\binom{m}{k}
\left(\frac{\lambar{1}}{\lambar{1} + \lambar{2}}\right)^k
\left(\frac{\lambar{2}}{\lambar{1} + \lambar{2}}\right)^{m-k}.
\]

Once $N(t)=n$ is given, suppose $N_1^a(t) = n^a$ and $N_1^b(t) = n^b$, so $n^a + n^b = n$. Because the system is out of stock for product~1 and a potential customer is considering this position in queue, there must be exactly $S_1$ type-1 replenishment jobs in total. The last job in the $n^a$ group must be the first of these replenishment jobs, and exactly $S_1-1$ type-1 replenishment jobs must be among the $n^b$ jobs that would be behind the potential customer's position.

Hence, the conditional probability of $\{N_1^a = n^a, N_1^b = n^b\}$, given that $N(t)=n$, involves choosing $(S_1 - 1)$ among $(n^b)$ positions to be type~1.  This selection yields a binomial coefficient $\binom{n^b}{S_1 - 1}$, multiplied by the usual powers of $\tfrac{\lambar{1}}{\lambar{1}+\lambar{2}}$ and $\tfrac{\lambar{2}}{\lambar{1}+\lambar{2}}$ for type-1 vs.\ type-2 classification:
\[
\mathbb{P}\left\{N_1^a = n^a, N_1^b = n^b \big|  N=n\right\}
=
\binom{n^b}{S_1 - 1}
\left(\frac{\lambar{1}}{\lambar{1} + \lambar{2}}\right)^{S_1}
\left(\frac{\lambar{2}}{\lambar{1} + \lambar{2}}\right)^{n^b - S_1+1}
\]
if $n^a+n^b=n$ and zero otherwise.

Combining \textbf{(i)} the M/M/1 stationary distribution of $N(t)$ (see \eqref{eq:distr_jobs_total}), \textbf{(ii)} the binomial-thinning argument, and \textbf{(iii)} the backlog-replacement constraint yields the unconditional joint distribution of $\left(N_1^a(t), N_1^b(t)\right)$:
\[
\mathbb{P}\{N_1^a = n^a,N_1^b = n^b\}
=
\binom{n^b}{S_1-1}
\frac{(\lambar{1}+\lambar{2})^{n^a-1}\lambar{1}^{S_1}\lambar{2}^{n^b - S_1+1}\left(\mu - \lambar{1} - \lambar{2}\right)}
     {\mu^{n^a + n^b + 1}},
\]
valid for $n^a \ge 1$ and $n^b \ge S_1-1$, and zero otherwise.  Intuitively, the total queue length is $n^a + n^b$, with $n^a$ tasks up to (and including) the target job, and $n^b$ jobs strictly behind.  Among those $n^b$ trailing jobs, $S_1-1$ are forced to be type~1 tasks, and the rest type~2.
\end{proof}

\subsection{Waiting Time Derivation}\label{app:waiting-time}
\begin{proof}[Proof of Proposition~\ref{prop:waiting_time}]
\noindent
\textit{Step 1: Expressing $\mathbb{E}[N_1^a]$.}  
From Section~\ref{subsec:arriving-cust-position}, the joint distribution of $(N_1^a, N_1^b)$ for a type-1 backlog arrival is
\[
\mathbb{P}\{N_1^a = n^a,N_1^b = n^b\}
=
\binom{n^b}{S_1-1}
\frac{(\lambar{1}+\lambar{2})^{n^a-1}\lambar{1}^{S_1}\lambar{2}^{n^b - S_1 + 1}\left(\mu - \lambar{1} - \lambar{2}\right)}
     {\mu^{n^a + n^b + 1}},
\]
valid for $n^a \ge 1$ and $n^b \ge S_1-1$, and zero otherwise.
Hence,
\[
\mathbb{E}[N_1^a]
=
\sum_{n^b = S_1-1}^\infty
\sum_{n^a = 1}^\infty
n^a
\mathbb{P}\{N_1^a = n^a,N_1^b = n^b\}.
\]
Substitute the expression for $\mathbb{P}\{N_1^a = n^a,N_1^b = n^b\}$ and group factors accordingly.

\medskip
\noindent
\textit{Step 2: Separating the Double Sum.}  
Note that
\[
\mathbb{E}[N_1^a]
=
(\mu - \lambar{1} - \lambar{2}) 
\lambar{1}^{S_1}
\left[\sum_{n^a = 1}^\infty n^a \frac{(\lambar{1}+\lambar{2})^{n^a-1}}{\mu^{n^a}}\right]
\left[\sum_{n^b = S_1-1}^\infty \binom{n^b}{S_1-1} \frac{\lambar{2}^{n^b - (S_1 - 1)}}{\mu^{n^b + 1}}\right].
\]
This follows by rewriting $\mu^{-(n^a + n^b + 1)}$ as $\mu^{-n^a}\mu^{-n^b - 1}$ and factoring out terms that do not depend on $n^a$ or $n^b$.

\medskip
\noindent
\textit{Step 3: Summation Over ${n^a}$.}  
The sum over $n^a$ is
\[
\sum_{n^a=1}^\infty n^a \frac{(\lambar{1}+\lambar{2})^{n^a-1}}{\mu^{n^a}}
=
\frac{1}{\mu}\sum_{n^a=1}^\infty n^a \left(\frac{\lambar{1}+\lambar{2}}{\mu}\right)^{n^a-1}
=
\frac{1}{\mu}\cdot\frac{1}{\left(1 - \tfrac{\lambar{1}+\lambar{2}}{\mu}\right)^2}
=
\frac{\mu}{(\mu - \lambar{1} - \lambar{2})^2}.
\]
(The identity $\sum_{k=1}^\infty kx^k = \frac{x}{(1-x)^2}$ for $\lvert x\rvert<1$ applies.)

\medskip
\noindent
\textit{Step 4: Summation Over ${n^b}$.}  
Let $k = n^b - (S_1 - 1)$. Then $k$ runs from $0$ to $\infty$, and
\[
\sum_{n^b = S_1 - 1}^\infty
\binom{n^b}{S_1 - 1}
\frac{\lambar{2}^{n^b - (S_1 - 1)}}{\mu^{n^b + 1}}
=
\frac{1}{\mu^{S_1}}
\sum_{k=0}^\infty
\binom{k + S_1 - 1}{S_1 - 1}
\left(\frac{\lambar{2}}{\mu}\right)^k.
\]
By the binomial identity $\sum_{k=0}^\infty \binom{k+n}{n} x^k = \tfrac{1}{(1-x)^{n+1}}$ (for $\lvert x\rvert<1$), one obtains
\[
\frac{1}{\mu^{S_1}}
\left(\frac{\mu}{\mu - \lambar{2}}\right)^{S_1}
=
\left(\frac{1}{\mu - \lambar{2}}\right)^{S_1}.
\]

\medskip
\noindent
\textit{Step 5: Final Computation.}  
Combining the two sums gives
\[
\mathbb{E}[N_1^a]
=
(\mu - \lambar{1} - \lambar{2})\lambar{1}^{S_1}
\frac{\mu}{(\mu - \lambar{1} - \lambar{2})^2}
\left(\frac{1}{\mu - \lambar{2}}\right)^{S_1}.
\]
After simplifying,
\[
\mathbb{E}[N_1^a]
=
\frac{\mu \lambar{1}^{S_1}}
     {\left(\mu - \lambar{2}\right)^{S_1}\left(\mu - \lambar{1} - \lambar{2}\right)}.
\]
Since $\mathbb{E}[W_1(\lambar{1},\lambar{2})] = \tfrac{1}{\mu}\mathbb{E}[N_1^a]$, the result in \eqref{eq:waiting_time1} follows immediately.

\medskip
\noindent
\textit{Symmetry for Product~2.}
An entirely similar argument applies when a newly arriving type-2 customer finds product~2 out of stock, enforcing exactly $S_2$ replenishment jobs for product~2.  Interchanging the roles of $\lambar{1}$, $\lambar{2}$, and $S_1$, $S_2$ yields
\begin{equation}\label{eq:waiting_time2}
\mathbb{E}[W_2(\lambar{1}, \lambar{2})]
=
\frac{\lambar{2}^{S_2}}
     {\left(\mu - \lambar{1}\right)^{S_2}\left(\mu - \lambar{1} - \lambar{2}\right)}.
\end{equation}
\end{proof}

\section{Strategic Behavior Proofs}\label{app:nash-analysis}
\subsection{Utility and Monotonicity Properties}\label{app:lem-monotonicity}
\begin{proof}[Proof of Lemma~\ref{lem:monotonicity}]
Continuity follows because $U_1$ is a composition of continuous functions on a compact domain. Monotonicity is shown by examining partial derivatives under the two scenarios $S_1=0$ and $S_1>0$:

\medskip
\noindent
\textbf{Scenario $S_1 = 0$.}
In this scenario, the utility simplifies to:
\[
U_1(q_1, q_2) = R_1 - p_1 - c_1 \cdot \frac{1}{\mu - q_1 \Lambda_1 - q_2 \Lambda_2}.
\]

Its derivative with respect to $q_1$ is:
\[
\frac{\partial U_1}{\partial q_1} = -\frac{c_1 \Lambda_1}{(\mu - q_1 \Lambda_1 - q_2 \Lambda_2)^2} < 0,
\]
and similarly,
\[
\frac{\partial U_1}{\partial q_2} = -\frac{c_1 \Lambda_2}{(\mu - q_1 \Lambda_1 - q_2 \Lambda_2)^2} < 0.
\]

\medskip
\noindent
\textbf{Scenario $S_1 > 0$.}
Here, we have:
\[
U_1(q_1, q_2) = R_1 - p_1 - c_1 \cdot \frac{(q_1 \Lambda_1)^{S_1}}{(\mu - q_2 \Lambda_2)^{S_1} (\mu - q_1 \Lambda_1 - q_2 \Lambda_2)}.
\]

A direct calculation shows:
\[
\frac{\partial U_1}{\partial q_1} = -\frac{c_1 (q_1 \Lambda_1)^{S_1} [S_1 (\mu - q_1 \Lambda_1 - q_2 \Lambda_2) + q_1 \Lambda_1]}{q_1 (\mu - q_1 \Lambda_1 - q_2 \Lambda_2)^2 (\mu - q_2 \Lambda_2)^{S_1}} < 0,
\]
and
\[
\frac{\partial U_1}{\partial q_2} = -\frac{c_1 \Lambda_2 (q_1 \Lambda_1)^{S_1} [S_1 (\mu - q_1 \Lambda_1 - q_2 \Lambda_2) + q_2 \Lambda_2]}{(\mu - q_1 \Lambda_1 - q_2 \Lambda_2)^2 (\mu - q_2 \Lambda_2)^{S_1}} < 0.
\]
Thus, $U_1(q_1, q_2)$ is strictly decreasing in both $q_1$ and $q_2$.
\end{proof}

\subsection{Unique Fixed Points}\label{app:lem-unique-fixed}
\begin{proof}[Proof of Lemma~\ref{lem:unique-fixed}]
Fix any $q_2\in[0,1]$. By Lemma~\ref{lem:monotonicity}, the utility function $U_1(q_1,q_2)$ is strictly decreasing in $q_1$ for every fixed $q_2$. This property ensures that the mapping $q_1\mapsto \mathrm{BR}_1(q_1,q_2)$ has a unique fixed point, as detailed below.
\begin{enumerate}
\item If $U_1(0,q_2)\le 0$, then for every $q_1\in[0,1]$ the inequality 
    \[
    U_1(q_1,q_2)\le U_1(0,q_2)\le 0
    \]
    holds. Consequently, the best response is always $\{0\}$; that is, $\mathrm{BR}_1(q_1,q_2)=\{0\}$ for all $q_1$. In this situation, the fixed point is $q_1=0$ if and only if $U_1(0,q_2)\le 0$.

\item If $U_1(1,q_2)\ge 0$, then for every $q_1\in[0,1]$ one has 
    \[
    U_1(q_1,q_2)\ge U_1(1,q_2)\ge 0.
    \]
    Therefore, the best response is always $\{1\}$ (i.e., $\mathrm{BR}_1(q_1,q_2)=\{1\}$ for all $q_1$), and the unique fixed point is $q_1=1$ if and only if $U_1(1,q_2)\ge 0$.

\item Otherwise, suppose that 
    \[
    U_1(0,q_2)>0 \quad \text{and} \quad U_1(1,q_2)<0.
    \]
    Then, by continuity and strict monotonicity of $U_1(q_1,q_2)$, the Intermediate Value Theorem guarantees the existence of a unique $q_1^*\in(0,1)$ such that
    \[
    U_1(q_1^*,q_2)=0.
    \]
    For $q_1<q_1^*$, the strict decrease implies $U_1(q_1,q_2)>0$ so that $\mathrm{BR}_1(q_1,q_2)=\{1\}$; for $q_1>q_1^*$, one obtains $U_1(q_1,q_2)<0$ and hence $\mathrm{BR}_1(q_1,q_2)=\{0\}$. At $q_1=q_1^*$, since $U_1(q_1^*,q_2)=0$, it holds that $\mathrm{BR}_1(q_1^*,q_2)=[0,1]$. Thus the only self-consistent fixed point in $[0,1]$ is $q_1^*$, which is the unique solution of $U_1(q_1,q_2)=0$. 
\end{enumerate}
A symmetric argument applies for type~2 by fixing $q_1$.
\end{proof}

\subsection{Detailed Classification of Nash Equilibria}\label{app:ne-classification}
\begin{proof}[Proof of Proposition~\ref{prop:NE-summary}]
\textit{Step 1: Possible Coordinates via Lemma~\ref{lem:unique-fixed}.}  
Applying Lemma~\ref{lem:unique-fixed} to each type shows that, for any $q_2\in[0,1]$, $\mathrm{BR}_1(q_1,q_2)$ has exactly one fixed point in $q_1$.  This unique fixed point is either $0$, $1$, or the fraction in $(0,1)$ where $U_1(q_1^*,q_2)=0$.  
Similarly for $q_2$ at any fixed $q_1$.  
Hence each $q_i$ must lie in $\{0\}\cup(0,1)\cup\{1\}$.  
Combining both coordinates yields the $9$ cases of Table~\ref{tab:NE_cases}.

\medskip
\noindent
\textit{Step 2: Verifying Each Cell's Inequalities.} 
\begin{description}[font=\normalfont\itshape, labelindent=1em]
    \item[Corner Cells:] 
\begin{itemize}[label=-]
\item $\boldsymbol{(0,0)}$.  
{Equilibrium if and only if}  
$
(0,0)
=
\left(\mathrm{BR}_1(0,0),\mathrm{BR}_2(0,0)\right),
$
which happens precisely when
\[
U_1(0,0)\le0
\quad\text{and}\quad
U_2(0,0)\le0.
\]
All type-1 and type-2 customers balk.

\item $\boldsymbol{(0,1)}$.  
{Equilibrium if and only if}  
$
(0,1)
=
\left(\mathrm{BR}_1(0,1),\mathrm{BR}_2(0,1)\right),
$
which happens exactly when
\[
U_1(0,1)\le0
\quad\text{and}\quad
U_2(0,1)\ge0.
\]
Type-1 balks, type-2 joins.

\item $\boldsymbol{(1,0)}$.  
{Equilibrium if and only if}  
$
(1,0)
=
\left(\mathrm{BR}_1(1,0),\mathrm{BR}_2(1,0)\right),
$
which happens if and only if
\[
U_1(1,0)\ge0
\quad\text{and}\quad
U_2(1,0)\le0.
\]
Type-1 joins, type-2 balks.

 $\boldsymbol{(1,1)}$.  
{Equilibrium if and only if}  
$
(1,1)
=
\left(\mathrm{BR}_1(1,1),\mathrm{BR}_2(1,1)\right),
$
which occurs precisely when
\[
U_1(1,1)\ge0
\quad\text{and}\quad
U_2(1,1)\ge0.
\]
Both types join.
\end{itemize}

\item[Partial-Corner Cases:]
\begin{itemize}[label=-]
\item $\boldsymbol{(0,q_2^*)}$ {with} $q_2^*\in(0,1)$.  
{Equilibrium if and only if}  
\[
(0,q_2^*)
=
\left(\mathrm{BR}_1(0,q_2^*),\mathrm{BR}_2(0,q_2^*)\right),
\]x
which requires
\[
U_1(0,q_2^*)\le0,
\quad
U_2(0,0)>0, 
U_2(0,1)<0, 
U_2(0,q_2^*)=0.
\]
Strict monotonicity ensures that if $U_2(0,q_2^*)=0$ has a unique solution $q_2^*\in(0,1)$, the inequalities $U_2(0,0)>0$ and $U_2(0,1)<0$ automatically hold.  Hence they do not appear explicitly in Table~\ref{tab:NE_cases}. Similar inequalities are likewise omitted in other cells of the table where a unique interior solution exists - the strict monotonicity of the utility functions ensures that the appropriate boundary conditions are satisfied whenever there is a solution in $(0,1)$.

Type-1 balks, type-2 mixes in $(0,1)$.

\item $\boldsymbol{(q_1^*,0)}$ {with} $q_1^*\in(0,1)$.  
{Equilibrium if and only if}  
\[
(q_1^*,0)
=
\left(\mathrm{BR}_1(q_1^*,0),\mathrm{BR}_2(q_1^*,0)\right),
\]
which is equivalent to
\[
U_1(0,0)>0,\quad
U_1(1,0)<0,\quad
U_1(q_1^*,0)=0,
\quad
U_2(q_1^*,0)\le0.
\]
Type-1 is partial in $(0,1)$, type-2 balks.

\item $\boldsymbol{(q_1^*,1)}$ {with} $q_1^*\in(0,1)$.  
{Equilibrium if and only if}  
\[
(q_1^*,1)
=
\left(\mathrm{BR}_1(q_1^*,1),\mathrm{BR}_2(q_1^*,1)\right),
\]
which necessitates
\[
U_1(0,1)>0,\quad
U_1(1,1)<0,\quad
U_1(q_1^*,1)=0,
\quad
U_2(q_1^*,1)\ge0.
\]
Type-1 is partial, type-2 joins.

\item $\boldsymbol{(1,q_2^*)}$ with $q_2^*\in(0,1)$. Equilibrium if and only if
\[
(1,q_2^*)
=
\left(\mathrm{BR}_1(1,q_2^*),\mathrm{BR}_2(1,q_2^*)\right),
\]
which implies
\[
U_1(1,q_2^*)\ge0,
\quad
U_2(1,0)>0, 
U_2(1,1)<0, 
U_2(1,q_2^*)=0.
\]
Type-1 joins, type-2 is partial in $(0,1)$.
\end{itemize}

\item[Interior Case:] 
\begin{itemize}[label=-]
\item $\boldsymbol{(q_1^*,q_2^*)}$ {with} $q_1^*,q_2^*\in(0,1)$.  
{Equilibrium if and only if}  
\[
(q_1^*,q_2^*)
=
\left(\mathrm{BR}_1(q_1^*,q_2^*),\mathrm{BR}_2(q_1^*,q_2^*)\right),
\]
which holds precisely when
\[
U_1(0,q_2^*)>0,\quad U_1(1,q_2^*)<0,\quad U_1(q_1^*,q_2^*)=0,
\]
\[
U_2(q_1^*,0)>0,\quad U_2(q_1^*,1)<0,\quad U_2(q_1^*,q_2^*)=0.
\]
Both types are partial (interior) joiners.
\end{itemize}
\end{description}

\medskip
\noindent
\textit{Conclusion:} These nine cases completely characterize the possible Nash equilibria, with each equilibrium lying in exactly one of the cells of Table~\ref{tab:NE_cases}.

\end{proof}

\subsection{Complete Nash Equilibrium Existence Proof}\label{app:ne-existence}
\begin{proof}[Proof of Proposition~\ref{prop:existenceNash}]
The proof uses Brouwer's Fixed Point Theorem, which states that any continuous function mapping a compact, convex set into itself has a fixed point. Although $\mathrm{BR}_1$ and $\mathrm{BR}_2$ are set-valued correspondences, a continuous single-valued function can be constructed from them via Lemma~\ref{lem:unique-fixed}.

\textit{Step 1: Constructing a Continuous Function from the Best-Response Correspondences.}\\
For each fixed $q_2\in [0,1]$, Lemma~\ref{lem:unique-fixed} guarantees a unique fixed point in $q_1$ for $\mathrm{BR}_1(q_1,q_2)$. Denote this point by $f_1(q_2)$. Likewise, for each fixed $q_1$, let $f_2(q_1)$ be the unique fixed point in $q_2$ for $\mathrm{BR}_2(q_1,q_2)$. Define
\[
f: [0,1]^2  \to  [0,1]^2
\quad\text{by}\quad
f(q_1,q_2)  =  \left(f_1(q_2),f_2(q_1)\right).
\]

\textit{Step 2: Continuity of $f_1$ and $f_2$.}\\
Focus on $f_1$. Fix a point $q_2^0$ in $[0,1]$, and let $\{q_2^n\}$ be a sequence converging to $q_2^0$. By Lemma~\ref{lem:unique-fixed}, each $f_1(q_2^n)$ is determined by one of the following:
\begin{enumerate}
\item $f_1(q_2^n)=0$ if $U_1(0,q_2^n)\le 0$,
\item $f_1(q_2^n)=1$ if $U_1(1,q_2^n)\ge 0$,
\item $f_1(q_2^n)=q_1^n\in(0,1)$ if $U_1(q_1^n,q_2^n)=0$ otherwise.
\end{enumerate}
Consider any convergent subsequence $\{f_1(q_2^{n_k})\}\to y$. By the continuity of $U_1$ (Lemma~\ref{lem:monotonicity}), the limit $y$ must satisfy the same case at $q_2^0$, i.e.\ $f_1(q_2^0)=y$. Since every convergent subsequence converges to $f_1(q_2^0)$, it follows that $f_1(q_2^n)\to f_1(q_2^0)$. Thus $f_1$ is continuous. By symmetry, $f_2$ is also continuous.

\textit{Step 3: Applying Brouwer's Fixed Point Theorem.}\\
Because $f$ is a continuous map from the compact, convex set $[0,1]^2$ to itself, Brouwer's theorem guarantees a fixed point $(q_1^*,q_2^*)$ such that
\[
f\left(q_1^*,q_2^*\right)  =  \left(q_1^*,q_2^*\right).
\]

\textit{Step 4: Showing the Fixed Point is a Nash Equilibrium.}\\
At this fixed point, $q_1^* = f_1(q_2^*)$ and $q_2^* = f_2(q_1^*)$. By construction of $f_1$ and $f_2$, it follows that $q_1^*\in \mathrm{BR}_1(q_1^*,q_2^*)$ and $q_2^*\in \mathrm{BR}_2(q_1^*,q_2^*)$. Hence $(q_1^*,q_2^*)$ is indeed a Nash equilibrium.
\end{proof}

\subsection{Analysis for Zero Inventory Scenario}\label{app:s0-s0}
\begin{proof}[Proof of Proposition~\ref{prop:S0-S0}]
\textit{Step 1: Detailed Classification of Equilibrium Cases.}

When $S_1 = 0$ and $S_2 = 0$, the utility functions reduce to:
\[
U_1(q_1, q_2)
=
R_1 - p_1
-
c_1 \frac{1}{\mu - q_1\Lambda_1 - q_2\Lambda_2},
\quad
U_2(q_1, q_2)
=
R_2 - p_2
-
c_2 \frac{1}{\mu - q_1\Lambda_1 - q_2\Lambda_2}.
\]
In this simplified setting, they allow explicit solutions for one type's joining probability if the other type's probability is fixed. Consequently, the general equilibrium conditions from Table~\ref{tab:NE_cases} break down into the specific cases shown in Table~\ref{tab:S0_S0_solutions-app}.

\begin{table}[ht]
\centering
\caption{Equilibrium Outcomes for $S_1 = 0, S_2 = 0$}
\label{tab:S0_S0_solutions-app}
\small
\begin{tabular}{|p{3.8cm}|p{3.9cm}|p{3.9cm}|}
\hline
\textbf{(0,0):} & 
\textbf{(0,$\boldsymbol{q_2^*}$):} & 
\textbf{(0,1):} \\
& $q_2^* = \frac{1}{\Lambda_2}(\mu - \frac{c_2}{R_2 - p_2})$: & \\
$\mu\le \frac{c_1}{R_1 - p_1}$ & 
$\mu > \frac{c_2}{R_2 - p_2}$ & 
$\mu \ge \frac{c_2}{R_2 - p_2} + \Lambda_2$ \\
$\mu\le \frac{c_2}{R_2 - p_2}$ & 
$\mu < \frac{c_2}{R_2 - p_2} + \Lambda_2$ & 
$\mu \le \frac{c_1}{R_1 - p_1} + \Lambda_2$ \\
& $\frac{c_1}{R_1 - p_1} \ge \frac{c_2}{R_2 - p_2}$ & \\
\hline
\textbf{($\boldsymbol{q_1^*}$,0),} & 
\textbf{($\boldsymbol{q_1^*}$,$\boldsymbol{q_2^*}$),} & 
\textbf{($\boldsymbol{q_1^*}$,1),} \\
$q_1^* = \frac{1}{\Lambda_1}(\mu - \frac{c_1}{R_1 - p_1})$: & 
$q_1^*\Lambda_1 + q_2^*\Lambda_2 = \mu - \frac{c_1}{R_1 - p_1}$, & 
$q_1^* = \frac{1}{\Lambda_1}(\mu - \Lambda_2 - \frac{c_1}{R_1 - p_1})$: \\
& $0 < q_1^*, q_2^* < 1$: & \\
$\mu > \frac{c_1}{R_1 - p_1}$ & 
$\mu > \frac{c_1}{R_1 - p_1}$ & 
$\mu > \frac{c_1}{R_1 - p_1} + \Lambda_2$ \\
$\mu < \frac{c_1}{R_1 - p_1} + \Lambda_1$ & 
$\mu < \frac{c_1}{R_1 - p_1} + \Lambda_1 + \Lambda_2$ & 
$\mu < \frac{c_1}{R_1 - p_1} + \Lambda_1 + \Lambda_2$ \\
$\frac{c_1}{R_1 - p_1} \le \frac{c_2}{R_2 - p_2}$ & 
$\frac{c_1}{R_1 - p_1} = \frac{c_2}{R_2 - p_2}$ & 
$\frac{c_1}{R_1 - p_1} \ge \frac{c_2}{R_2 - p_2}$ \\
\hline
\textbf{(1,0):} & 
\textbf{(1,$\boldsymbol{q_2^*}$),} & 
\textbf{(1,1):} \\
& $q_2^* = \frac{1}{\Lambda_2}(\mu - \Lambda_1 - \frac{c_2}{R_2 - p_2})$: & \\
$\mu \ge \frac{c_1}{R_1 - p_1} + \Lambda_1$ & 
$\mu > \frac{c_2}{R_2 - p_2} + \Lambda_1$ & 
$\mu \ge \frac{c_1}{R_1 - p_1} + \Lambda_1 + \Lambda_2$ \\
$\mu \le \frac{c_2}{R_2 - p_2} + \Lambda_1$ & 
$\mu < \frac{c_2}{R_2 - p_2} + \Lambda_1 + \Lambda_2$ & 
$\mu \ge \frac{c_2}{R_2 - p_2} + \Lambda_1 + \Lambda_2$ \\
& $\frac{c_1}{R_1 - p_1} \le \frac{c_2}{R_2 - p_2}$ & \\
\hline
\end{tabular}
\end{table}

\noindent
To clarify how the cells in Table~\ref{tab:S0_S0_solutions-app} arise, consider the following points.
\begin{itemize}
    \item An interior equilibrium $(q_1^*, q_2^*) \in (0,1)^2$ must satisfy 
\[
U_1(q_1^*,q_2^*)=0
\quad\text{and}\quad
U_2(q_1^*,q_2^*)=0.
\]

Subtracting these two equations implies
\[
\frac{c_1}{R_1 - p_1}
=
\frac{c_2}{R_2 - p_2}.
\]

Under this condition, substituting back into $U_1=0$ (or $U_2=0$) yields a line:
\[
q_1^*\Lambda_1 + q_2^*\Lambda_2 = \mu - \frac{c_1}{R_1 - p_1}.
\]

To ensure at least one point of this line lies in $(0,1)^2$, it must satisfy:
\[
q_2^*(0) > 0
\quad\text{and}\quad
q_2^*(1) < 1.
\]

Plugging in $q_1^*=0$ and $q_1^*=1$ gives:
\[
\frac{1}{\Lambda_2}\left(\mu - \tfrac{c_1}{R_1-p_1}\right) > 0,
\quad
\frac{1}{\Lambda_2}\left(\mu - \tfrac{c_1}{R_1-p_1} - \Lambda_1\right) < 1,
\]
which is equivalent to
\[
\frac{c_1}{R_1 - p_1} < \mu < \frac{c_1}{R_1 - p_1} + \Lambda_1 + \Lambda_2.
\]

\item To see how other cells arise, consider $(q_1^*, 1)$. For it to be an equilibrium, 
\[
U_1(q_1^*,1)=0
\quad\text{and}\quad
U_2(q_1^*,1)\ge0,
\]
which translates into
\[
R_1 - p_1
-
c_1 \frac{1}{\mu - q_1^*\Lambda_1 - \Lambda_2} = 0
\quad\text{and}\quad
R_2 - p_2
-
c_2 \frac{1}{\mu - q_1^*\Lambda_1 - \Lambda_2} \ge 0.
\]
Hence,
\[
q_1^*
= 
\frac{1}{\Lambda_1}\left(\mu - \Lambda_2 - \tfrac{c_1}{R_1 - p_1}\right)
\quad\text{and}\quad
\frac{c_1}{R_1 - p_1} 
\ge 
\frac{c_2}{R_2 - p_2}.
\]
Since $q_1^*$ must lie in $(0,1)$, 
\[
\frac{c_1}{R_1 - p_1} + \Lambda_2 
<
\mu 
<
\frac{c_1}{R_1 - p_1} + \Lambda_1 + \Lambda_2.
\]

\item Similar computations under other assumptions on $(q_1^*,q_2^*)$ yield the remaining cells in Table~\ref{tab:S0_S0_solutions-app}.
\end{itemize}

\medskip
\noindent
\textit{Step 2: Equilibrium Coexistence Analysis.}
Equilibria may or may not coexist, depending on how
$\tfrac{c_1}{R_1 - p_1}$ compares to $\tfrac{c_2}{R_2 - p_2}$. Three cases arise:
\[
(1)    \frac{c_1}{R_1 - p_1} 
>
\frac{c_2}{R_2 - p_2}, 
\quad
(2)    \frac{c_1}{R_1 - p_1} 
<
\frac{c_2}{R_2 - p_2},
\quad
\text{or}
\quad
(3)    \frac{c_1}{R_1 - p_1}
=
\frac{c_2}{R_2 - p_2}.
\]
Cases (1) and (2) are symmetric by exchanging the roles of products~1 and~2, so the details below focus on (1) and (3).
\begin{description}[font=\normalfont\itshape, labelindent=1em]
    \item[\textit{Case 1:}] $\tfrac{c_1}{R_1 - p_1} > \tfrac{c_2}{R_2 - p_2}$.
Because the conditions in other cells conflict with this inequality, 
only $(0,0)$, $(0,q_2^*)$, $(0,1)$, $(q_1^*,1)$, and $(1,1)$ remain feasible. Intuitively, type-1 has a higher cost ratio and is therefore less inclined to join, so $q_1^* \le q_2^*$. In this setting:
\begin{itemize}[label=-]
\item 
${(0,0)}$ requires 
$\mu \le \tfrac{c_1}{R_1-p_1}$ {and} $\mu \le \tfrac{c_2}{R_2-p_2}$. 
Since 
$\tfrac{c_1}{R_1-p_1} > \tfrac{c_2}{R_2-p_2}$,
the binding constraint is 
$\mu \le \tfrac{c_2}{R_2-p_2}$.
All other cells in this case require $\mu$ to be strictly above $\tfrac{c_2}{R_2-p_2}$ (either to exceed it or to exceed $\tfrac{c_1}{R_1-p_1}$ or to add $\Lambda_i$), so $(0,0)$ cannot coincide with any of them.

\item 
${(1,1)}$ requires 
$
\mu 
\ge 
\frac{c_1}{R_1 - p_1} + \Lambda_1 + \Lambda_2$ 
 {and} 
$\mu 
\ge 
\frac{c_2}{R_2 - p_2} + \Lambda_1 + \Lambda_2.
$
Because 
$\tfrac{c_1}{R_1-p_1} > \tfrac{c_2}{R_2-p_2}$, 
the first condition is stricter, so 
$\mu \ge \tfrac{c_1}{R_1 - p_1} + \Lambda_1 + \Lambda_2.$
All other cells in this case require $\mu$ to be below that threshold (or within a smaller range), so $(1,1)$ also cannot coexist with any other cell.
\item $(0,q_2^*)$ and $(0,1)$ cannot coexist because by Lemma~\ref{lem:unique-fixed}, for fixed $q_1=0$, there is exactly one fixed point of $\mathrm{BR}_2(0,q_2)$ in $q_2$.
\item $(0,1)$ and $(q_1^*,1)$ cannot coexist because by Lemma~\ref{lem:unique-fixed}, for fixed $q_2=1$, there is exactly one fixed point of $\mathrm{BR}_1(q_1,1)$ in $q_1$.
\item 
${(0,q_2^*)}$ and ${(q_1^*,1)}$ have mutually exclusive $\mu$-ranges, since
\[
(0,q_2^*): 
\quad 
\tfrac{c_2}{R_2 - p_2} 
<
\mu 
<
\tfrac{c_2}{R_2 - p_2} + \Lambda_2,
\]
while
\[
(q_1^*,1): 
\quad 
\tfrac{c_1}{R_1 - p_1} + \Lambda_2
< 
\mu 
<
\tfrac{c_1}{R_1 - p_1} + \Lambda_1 + \Lambda_2.
\]
Because 
$\tfrac{c_1}{R_1 - p_1} > \tfrac{c_2}{R_2 - p_2}$,
we have 
$\tfrac{c_1}{R_1 - p_1} + \Lambda_2 > \tfrac{c_2}{R_2 - p_2} + \Lambda_2$.
Thus, these intervals do not overlap, so $(0,q_2^*)$ and $(q_1^*,1)$ also cannot coexist.
\end{itemize}

Moreover, whenever one strategy is fixed at $0$ or $1$, the other type's best response must be unique (Lemma~\ref{lem:unique-fixed}). Hence, we cannot have two different equilibria of the form $(1,q_2)$ and $(1,q'_2)$, or $(0,q_2)$ and $(0,q'_2)$. Consequently, no two equilibria can coexist on the same boundary cell, ensuring uniqueness in this case.

\item[\textit{Case 3:}] $\tfrac{c_1}{R_1 - p_1} = \tfrac{c_2}{R_2 - p_2}$.
If 
\[
\mu  \le  \frac{c_1}{R_1 - p_1}
\quad\text{or}\quad
\mu  \ge  \frac{c_1}{R_1 - p_1} + \Lambda_1 + \Lambda_2,
\]
the equilibrium is either $(0,0)$ or $(1,1)$, and it is unique.

On the other hand, if 
\[
\tfrac{c_1}{R_1 - p_1} 
 <  
\mu 
 <  
\tfrac{c_1}{R_1 - p_1} + \Lambda_1 + \Lambda_2,
\]
the interior cell conditions hold. In this case, the entire line segment
\[
(q_1^*,q_2^*)
\quad\text{where}\quad
q_1^*\Lambda_1 + q_2^*\Lambda_2 = \mu - \frac{c_1}{R_1 - p_1},
\quad
0 \le q_1^*,q_2^* \le 1
\]
forms a continuum of equilibria. In particular, points on that line with $q_i^*$ at 0 or 1 lie on boundary cells in Table~\ref{tab:S0_S0_solutions-app} and thus connect the interior segment to two of the border equilibria.
\end{description}
\end{proof}

\subsection{Analysis for Positive Inventory Scenario}\label{app:positive-inventory}
\begin{proof}[Proof of Proposition~\ref{prop:positive-inventory}]
\noindent
\textit{Step 1: $S_i > 0$ implies $q_i > 0$ in any equilibrium.}

If $S_i>0$ and $q_i=0$, then
\[
U_i\left(0,q_j\right)
 = 
R_i  -  p_i 
 - 
c_i
\frac{(0)^{S_i}}{(\mu - q_j\Lambda_j)^{S_i}\left(\mu - q_j\Lambda_j\right)}
 = 
R_i - p_i
 > 0
\]
(since $R_i>p_i$).  
Hence type~$i$ would strictly prefer to join, making $q_i=0$ impossible in equilibrium.

\medskip
\noindent
\textit{Step 2: Uniqueness when both $S_1, S_2 > 0$.}

\medskip
\noindent
\textit{Step 2.1: Detailed Classification of Equilibrium Cases.} To illustrate the possible equilibrium scenarios when $S_1 > 0$ and $S_2 > 0$, first note that Table~\ref{tab:NE_cases} includes cells where $q_i=0$. However, under $S_i>0$, we have shown that $q_i=0$ cannot arise in equilibrium. Consequently, only cells with $q_1,q_2\in(0,1]$ remain relevant. These lead to the equilibrium configurations summarized in Table~\ref{tab:S-positive-cases}:

\begin{table}[h!]
\caption{General Forms of Equilibrium Conditions for $S_1 > 0, S_2 > 0$}
\label{tab:S-positive-cases}
\centering
\setcellgapes{4pt}  
\makegapedcells 
\begin{tabular}{|l|l|l|}
\hline
\makecell[l]{
  \textbf{($\boldsymbol{q_1^*}$,$\boldsymbol{q_2^*}$)}, $q_1^*,q_2^*\in(0,1)$ --- NE:\\[4pt]
  $q_1^* = \BRf{1}(q_1^*,q_2^*), q_2^* = \BRf{2}(q_1^*,q_2^*)$\\[8pt]
  $U_1(q_1^*,q_2^*) = 0$, $U_2(q_1^*,q_2^*) = 0$ 
}
&
\makecell[l]{
  \textbf{($\boldsymbol{q_1^*}$,1)}, $q_1^*\in(0,1)$ --- NE:\\[4pt]
  $q_1^* = \BRf{1}(q_1^*,1), 1 = \BRf{2}(q_1^*,1)$\\[8pt]
  $U_1(q_1^*,1) = 0$, $U_2(q_1^*,1) \ge 0$ 
}
\\ \hline
\makecell[l]{
  \textbf{(1,$\boldsymbol{q_2^*}$)}, $q_2^*\in(0,1)$ --- NE:\\[4pt]
  $1 = \BRf{1}(1,q_2^*), q_2^* = \BRf{2}(1,q_2^*)$\\[8pt]
  $U_1(1,q_2^*) \ge 0$, $U_2(1,q_2^*) = 0$ 
}
&
\makecell[l]{
  \textbf{(1,1)} --- NE:\\[4pt]
  $1 = \BRf{1}(1,1), 1 = \BRf{2}(1,1)$\\[8pt]
  $U_1(1,1) \ge 0$, $U_2(1,1) \ge 0$ 
}
\\ \hline
\end{tabular}
\end{table}

We next show that only one of these cases can actually occur, and when it does, the equilibrium is unique. The argument proceeds in two steps.

\medskip
\noindent
\textit{Step 2.2: Interior Analysis.} Consider the case where $\left(q_1^*,q_2^*\right)$ has both $q_1^*,q_2^*\in(0,1)$ and satisfies
\[
U_1\left(q_1^*,q_2^*\right) = 0
\quad\text{and}\quad
U_2\left(q_1^*,q_2^*\right) = 0.
\]

Here, the goal is to show that any such interior equilibrium must be unique {within} $(0,1)^2$.

To handle interior solutions $(q_1,q_2)\in(0,1)^2$, define two functions:
\begin{align}
& g_1(q_2): \text{ the unique } q_1 \text{ satisfying } U_1(q_1, q_2) = 0, \nonumber \\
& g_2(q_1): \text{ the unique } q_2 \text{ satisfying } U_2(q_1, q_2) = 0.
\end{align}

These functions are well-defined by monotonicity of $U_1$ and $U_2$ (Lemma~\ref{lem:monotonicity}). Then consider the composition
\[
\varphi(q_1) = g_1(g_2(q_1)).
\]

If multiple equilibria existed in $(0,1)^2$, there would be multiple fixed points of $\varphi$, i.e.\ multiple solutions to $\varphi(q_1)=q_1$.

Using the Implicit Function Theorem, one obtains:
\[
\varphi'(q_1)
=
\frac{\partial g_1}{\partial q_2}
\cdot
\frac{\partial g_2}{\partial q_1}
=
\left(-\frac{\frac{\partial U_1}{\partial q_2}}{\frac{\partial U_1}{\partial q_1}}\right)
\cdot
\left(-\frac{\frac{\partial U_2}{\partial q_1}}{\frac{\partial U_2}{\partial q_2}}\right).
\]

For $S_i>0$, the partial derivatives are:
\[
\frac{\partial U_1}{\partial q_1}
=
-\frac{
c_1 \left(q_1 \Lambda_1\right)^{S_1}
\left[S_1\left(\mu - q_1\Lambda_1 - q_2\Lambda_2\right) + q_1 \Lambda_1\right]
\left(\mu - q_2 \Lambda_2\right)^{-S_1}
}{
q_1 \left(q_1 \Lambda_1 + q_2 \Lambda_2 - \mu\right)^2
},
\]
\[
\frac{\partial U_1}{\partial q_2}
=
-
\frac{
c_1 \Lambda_2\left(q_1 \Lambda_1\right)^{S_1}
\left[S_1\left(\mu - q_1\Lambda_1 - q_2\Lambda_2\right) + q_2 \Lambda_2\right]
\left(\mu - q_2 \Lambda_2\right)^{-S_1}
}{
\left(q_1 \Lambda_1 + q_2 \Lambda_2 - \mu\right)^2
}.
\]

Similar expressions hold for $\frac{\partial U_2}{\partial q_1}$ and $\frac{\partial U_2}{\partial q_2}$. Substituting these into $\varphi'(q_1)$ and simplifying yields
\[
\varphi'(q_1)
=
\frac{
q_1\Lambda_1q_2\Lambda_2
\left(S_1(\mu-q_1\Lambda_1-q_2\Lambda_2) + \mu - q_2\Lambda_2\right)
\left(S_2(\mu-q_1\Lambda_1-q_2\Lambda_2) + \mu - q_1\Lambda_1\right)
}{
(q_1\Lambda_1-\mu)(\mu-q_2\Lambda_2)\left(S_1 (\mu-q_1\Lambda_1-q_2\Lambda_2) + q_1\Lambda_1\right)\left(S_2 (\mu-q_1\Lambda_1-q_2\Lambda_2) + q_2\Lambda_2\right)
}.
\]

Recall from Section~\ref{sec:model} that $\mu > \Lambda_1 + \Lambda_2$. Hence, for $q_1,q_2 \le 1$, $\mu > q_1\Lambda_1 + q_2\Lambda_2$ ensures all factors in the fraction remain positive, implying $\varphi'(q_1) > 0$.  
What remains is to show that $\varphi'(q_1) < 1$.

Let
\[
A = \mu - q_1\Lambda_1 - q_2\Lambda_2,\quad
B_1 = \mu - q_1\Lambda_1,\quad
B_2 = \mu - q_2\Lambda_2.
\]

Then
\[
\varphi'(q_1)
=
\frac{(B_1 - A)(B_2 - A)\left(S_1 A + B_2\right)\left(S_2 A + B_1\right)}
{B_1B_2\left(S_1 A + B_2 - A\right)\left(S_2 A + B_1 - A\right)}.
\]

Dividing the numerator and denominator by $\left(B_1 B_2 (S_1 A + B_2) (S_2 A + B_1)\right)$ simplifies the expression to:
\[
\varphi'(q_1)
=
\frac{
\left(1-\frac{A}{B_1}\right)
\left(1-\frac{A}{B_2}\right)
}{
\left(1-\frac{A}{S_2 A + B_1}\right)
\left(1-\frac{A}{S_1 A + B_2}\right)
},
\]
and since $\frac{A}{S_i A + B_j} < \frac{A}{B_j}$, each factor in the denominator is smaller than its counterpart in the numerator, implying $\varphi'(q_1) < 1$. 

Thus, $\varphi$ is a strictly increasing function with derivative less than $1$, indicating that its slope is always below that of the line $y = x$. Consequently, $\varphi$ can cross $y=x$ at most once, so $\varphi(q_1) = q_1$ has at most one solution. In other words, there is at most a single fixed point, which is the unique interior Nash equilibrium $\left(q_1^*,q_2^*\right)$.

\medskip
\noindent
\textit{Step 2.3: Boundary Analysis.} Suppose there is already an equilibrium $\left(q_1^*,1\right)$ with $q_1^*\in(0,1]$. 
Another equilibrium of the same form $\left(\widetilde{q_1^*},1\right)$ cannot arise, 
because by Lemma~\ref{lem:unique-fixed}, for fixed $q_2=1$, there is exactly one fixed point of $\mathrm{BR}_1(q_1,1)$ in $q_1$.
It remains to show that no second equilibrium of the form $\left(\widetilde{q_1^*},q_2^*\right)$ 
can exist with $\widetilde{q_1^*}\in(0,1]$ and $q_2^*\in(0,1)$.  

Assume by contradiction that such a second equilibrium $(\widetilde{q_1^*},q_2^*)$ does exist. By Lemma~\ref{lem:unique-fixed}, for any fixed $q_2$, let $f_1(q_2)$ denote the unique fixed point of $\mathrm{BR}_1(q_1,q_2)$ in $q_1$, and similarly for any fixed $q_1$, let $f_2(q_1)$ denote the unique fixed point of $\mathrm{BR}_2(q_1,q_2)$ in $q_2$. Define
\[
\psi(q_1) = f_1(f_2(q_1)).
\]

Since any Nash equilibrium $\left(q_1^*,q_2^*\right)$ must satisfy $q_2^*=f_2(q_1^*)$ and $q_1^*=f_1(q_2^*)$ (as these are fixed points of the respective best responses), having two equilibria would imply two distinct solutions to $\psi(q_1)=q_1$.

Next, observe that $\left(q_1^*,1\right)$ being an equilibrium means $f_2(q_1^*)=1$, while the other equilibrium $(\widetilde{q_1^*},q_2^*)$ has $f_2(\widetilde{q_1^*})<1$. By monotonicity of $U_2$ (Lemma~\ref{lem:monotonicity}), there must be a threshold $q_1'\in[q_1^*,1]$ such that
\[
f_2(q_1)
=
\begin{cases}
1, & 0 \le q_1 \le q_1',\\
g_2(q_1), & q_1 > q_1',
\end{cases}
\]
where $g_2(q_1)$ is the interior solution defined earlier (i.e., the value making $U_2=0$).

Hence,
\[
\psi(q_1)
=
f_1(f_2(q_1))
=
\begin{cases}
f_1(1), & 0 \le q_1 \le q_1',\\
f_1(g_2(q_1)), & q_1 > q_1'.
\end{cases}
\]

On $[0,q_1']$, $\psi$ is constant at the value $f_1(1)$, giving it slope $0$.
For $q_1 > q_1'$, $\psi(q_1)$ is either identically $1$ or given by $g_1\left(g_2(q_1)\right)$, which is strictly increasing with derivative less than $1$ (Step~2's argument). 
Moreover, $\psi$ is continuous at $q_1'$. 
Hence, across $[0,1]$, $\psi$ is piecewise constant and piecewise strictly increasing with slope $<1$. Therefore, it can intersect the line $y=x$ at most once.

Since $q_1=q_1^*$ is already one intersection (the known equilibrium $\left(q_1^*,1\right)$), a second intersection would contradict the above argument that $\psi$ crosses $y=x$ only once. Thus, no second solution is possible, completing the contradiction.

The same reasoning applies if the boundary equilibrium is $\left(1, q_2^*\right)$. Therefore, if $\left(q_1^*, 1\right)$ (or $\left(1, q_2^*\right)$) exists, no other equilibrium is possible, completing the uniqueness proof.

\medskip
\noindent
\textit{Step 3: Uniqueness with asymmetric inventory.}

Assume without loss of generality that $S_1 = 0$ and $S_2 > 0$.

\medskip
\noindent
\textit{Step 3.1: Detailed Classification of Equilibrium Cases.} All possible forms of equilibrium conditions appear in Table~\ref{tab:S1zero-S2positive-cases}, which is analogous to Table~\ref{tab:S-positive-cases} from the $S_1>0, S_2>0$ scenario. The only difference here is that $q_1$ may be zero, introducing an extra boundary row.

\begin{table}[ht]
\caption{Possible Equilibria for $S_1 = 0, S_2 > 0$}
\label{tab:S1zero-S2positive-cases}
\centering
\setcellgapes{4pt}  
\makegapedcells 
\begin{tabular}{|l|l|l|}
\hline
\makecell[l]{
  \textbf{(0, $\boldsymbol{q_2^*}$)}, $q_2^*\in(0,1)$ --- NE:\\[4pt]
  $0 = \BRf{1}(0,q_2^*), q_2^* = \BRf{2}(0,q_2^*)$.\\[8pt]
  $U_1(0,q_2^*) \le 0$, $U_2(0,q_2^*) = 0$ 
}
&
\makecell[l]{
  \textbf{(0,1)} --- NE:\\[4pt]
  $0 = \BRf{1}(0,1), 1 = \BRf{2}(0,1)$\\[8pt]
  $U_1(0,1) \le 0$, $U_2(0,1) \ge 0$ 
}
\\ \hline
\makecell[l]{
  \textbf{($\boldsymbol{q_1^*}$,$\boldsymbol{q_2^*}$)}, $q_1^*,q_2^*\in(0,1)$ --- NE:\\[4pt]
  $q_1^* = \BRf{1}(q_1^*,q_2^*), q_2^* = \BRf{2}(q_1^*,q_2^*)$\\[8pt]
  $U_1(q_1^*,q_2^*) = 0$, $U_2(q_1^*,q_2^*) = 0$ 
}
&
\makecell[l]{
  \textbf{($\boldsymbol{q_1^*}$,1)}, $q_1^*\in(0,1)$ --- NE:\\[4pt]
  $q_1^* = \BRf{1}(q_1^*,1), 1 = \BRf{2}(q_1^*,1)$\\[8pt]
  $U_1(q_1^*,1) = 0$, $U_2(q_1^*,1) \ge 0$ 
}
\\ \hline
\makecell[l]{
  \textbf{(1,$\boldsymbol{q_2^*}$)}, $q_2^*\in(0,1)$ --- NE:\\[4pt]
  $1 = \BRf{1}(1,q_2^*), q_2^* = \BRf{2}(1,q_2^*)$\\[8pt]
  $U_1(1,q_2^*) \ge 0$, $U_2(1,q_2^*) = 0$ 
}
&
\makecell[l]{
  \textbf{(1,1)} --- NE:\\[4pt]
  $1 = \BRf{1}(1,1), 1 = \BRf{2}(1,1)$\\[8pt]
  $U_1(1,1) \ge 0$, $U_2(1,1) \ge 0$ 
}
\\ \hline
\end{tabular}
\end{table}

Four main equilibrium categories remain from the strictly positive scenario: $(q_1^*, q_2^*)$, $(q_1^*, 1)$, $(1, q_2^*)$, and $(1, 1)$. In addition, two new boundary cells arise here: $(0,q_2^*)$ and $(0,1)$. The argument proceeds in two steps.

\medskip
\noindent
\textit{Step 3.2: Interior Analysis.} When $q_1>0$ and $q_2<1$, consider points where $U_1$ and $U_2$ are zero. To handle these interior solutions $(q_1,q_2)\in(0,1)^2$, define as before:
\begin{align}
& g_1(q_2): \text{ the unique } q_1 \text{ satisfying } U_1(q_1, q_2) = 0, \nonumber \\
& g_2(q_1): \text{ the unique } q_2 \text{ satisfying } U_2(q_1, q_2) = 0.
\end{align}

The partial derivative calculations are nearly identical to those in Step 2, except $\frac{\partial U_1}{\partial q_1}$ and $\frac{\partial U_1}{\partial q_2}$ simplify because $S_1=0$. Specifically,
\[
\frac{\partial U_1}{\partial q_1}
=
-
\frac{c_1\Lambda_1}{(\mu - q_1\Lambda_1 - q_2\Lambda_2)^2}, \quad \frac{\partial U_1}{\partial q_2} = -\frac{c_1 \Lambda_2}{(\mu - q_1\Lambda_1 - q_2\Lambda_2)^2}.
\]

This yields a similar expression for $\varphi'(q_1)$ in the composition $\varphi(q_1) = g_1(g_2(q_1))$. One obtains
\[
\varphi'(q_1)
=
\frac{(B_1 - A)\left(S_2A + B_1\right)}
{B_1\left(S_2A + B_1 - A\right)},
\]
where $A=\mu - q_1\Lambda_1 - q_2\Lambda_2$ and $B_1=\mu - q_1\Lambda_1$. By the same ratio argument as before, $0<\varphi'(q_1)<1$, so at most one interior equilibrium exists.

\medskip
\noindent
\textit{Step 3.3: Boundary Analysis.} By Lemma~\ref{lem:unique-fixed}, for any fixed $q_2$, let $f_1(q_2)$ denote the unique fixed point of $\mathrm{BR}_1(q_1,q_2)$ in $q_1$, and similarly for any fixed $q_1$, let $f_2(q_1)$ denote the unique fixed point of $\mathrm{BR}_2(q_1,q_2)$ in $q_2$.

Any equilibria from $\left(q_1^*,1\right)$, $\left(1,q_2^*\right)$ or $(1,1)$ cannot coexist by the argument as in Step 2.3.

Two equilibria of the form $(0,q_2^*)$ with $q_2^*\in(0,1]$ cannot coexist since by Lemma~\ref{lem:unique-fixed}, for fixed $q_1=0$, there is exactly one fixed point of $\mathrm{BR}_2(0,q_2)$ in $q_2$. It remains to prove that $(0,q_2^*)$ with $q_2^*\in(0,1]$ cannot coexist with any of $(q_1^*,1)$, $(1,q_2^*)$, $(1,1)$. 

Indeed, suppose $(0,q_2^*)$ where $q_2^*\in(0,1]$ and $(q_1^*, \widetilde{q_2^*})$ where $q_1^*, \widetilde{q_2^*} \in (0,1]$ coexist. Define
\[
\psi(q_1) = f_1(f_2(q_1)).
\]
By the same argument as before, $\psi$ is continuous and consists of (1) a constant piece $f_1(f_2(0))$ for small $q_1$; (2) a strictly increasing segment with slope less than 1 for larger $q_1$.  This contradicts the existence of two fixed points of $\psi$.

 Hence, there is a {unique} Nash equilibrium overall for $S_1=0$ and $S_2>0$. 
\end{proof}

\section{Profit Optimization Proofs}\label{app:producer-optimality}
\subsection{Inventory Threshold Derivation}\label{app:thresholds}
\begin{proof}[Proof of Proposition~\ref{prop:thresholds}]
\textit{Step 1: Positivity of the utility functions.}
Let $U_i^{S_i}(1,1)$ be the utility of a type-$i$ customer when the target inventory level is $S_i$ and all customers join. Concretely,
\begin{align*}
&U_1^{S_1}(1,1)
=
R_1 - p_1
-
c_1 \frac{\Lambda_1^{S_1}}{(\mu - \Lambda_2)^{S_1}\left(\mu - \Lambda_1 - \Lambda_2\right)},\\
&U_2^{S_2}(1,1)
=
R_2 - p_2
-
c_2 \frac{\Lambda_2^{S_2}}{(\mu - \Lambda_1)^{S_2}\left(\mu - \Lambda_1 - \Lambda_2\right)}.
\end{align*}

Each $U_i^{S_i}(1,1)$ increases strictly in $S_i$, and
\[
\lim_{S_i \to \infty} U_i^{S_i}(1,1) = R_i - p_i > 0,
\]
since $R_i > p_i$.

\medskip
\noindent
\textit{Step 2: Thresholds that guarantee full joining.}
Define
\[
\bar{S}_1 
= 
\min\{S_1 \mid U_1^{S_1}(1,1) > 0\},
\quad
\bar{S}_2 
= 
\min\{S_2 \mid U_2^{S_2}(1,1) > 0\}.
\]

Because $U_i^{S_i}(1,1)$ is strictly increasing and eventually positive, each $\bar{S}_i$ is finite. If $S_i \ge \bar{S}_i$, then $U_i^{S_i}(1,1)>0$, so type-$i$ customer strictly prefers to join. Hence, whenever $S_1 \ge \bar{S}_1$ and $S_2 \ge \bar{S}_2$, all customers join for both products, 
which is to say $q_1^{(S_1,S_2)} = q_2^{(S_1,S_2)} = 1$.

\medskip
\noindent
\textit{Step 3: Profit implications of large inventories.}
The total profit can be written as
$
\Pi(S_1,S_2)
=
\Pi_1(S_1)
+
\Pi_2(S_2),
$
where 
\[
\Pi_i(S_i)
=
p_i\Lambda_iq_i^{(S_1,S_2)} 
-
h_i\mathbb{E}[I_i(S_i)].
\]

 If $S_i \ge \bar{S}_i$, then all customers of type~$i$ join, making
$
p_i\Lambda_iq_i^{(S_1,S_2)}
=
p_i\Lambda_i
$
 constant. Meanwhile, the expected inventory 
\[
\mathbb{E}[I_i(S_i)]=S_i - \frac{\Lambda_i}{\mu - \Lambda_1 - \Lambda_2} \left( 1 - \left( \frac{\Lambda_i}{\mu - \Lambda_j} \right)^{S_i} \right)
\]
is an increasing function of $S_i$.  Hence, $\Pi_i(S_i)$ decreases once $S_i > \bar{S}_i$.  Consequently,
\[
\Pi_i(S_i) 
\le 
\Pi_i(\bar{S}_i)
\quad
\text{whenever } 
S_i \ge \bar{S}_i.
\]

Summing over both products,
\[
\Pi(S_1,S_2)
=
\Pi_1(S_1) + \Pi_2(S_2)
\le
\Pi_1(\bar{S}_1) + \Pi_2(\bar{S}_2)
=
\Pi\left(\bar{S}_1,\bar{S}_2\right),
\]
which completes the proof.
\end{proof}

\section{Welfare Optimization Proofs}\label{app:social-optimization}
\subsection{Optimal Target Inventory Levels Derivation}\label{app:optimal-si}
\begin{proof}[Proof of Proposition~\ref{prop:optimal_Si}]
\textit{Step 1: Separation of the Cost Function.} 
In this proof, $C(\lambar{1}, \lambar{2}, S_1, S_2)$ is treated as $C(S_1, S_2)$, reflecting that $\lambar{i}$ are constant.

Observe that the total cost can be written as
\[
C(S_1, S_2) = C_1(S_1) + C_2(S_2),
\]
where
\[
C_i(S_i) = h_i S_i - h_i  \frac{\lambar{i}}{\mu - \lambar{1} - \lambar{2}} \left( 1 - \left( \frac{\lambar{i}}{\mu - \lambar{j}} \right)^{S_i} \right) + c_i  \frac{\lambar{i}}{\mu - \lambar{1} - \lambar{2}} \left( \frac{\lambar{i}}{\mu - \lambar{j}} \right)^{S_i}.
\]
Since $C_1(S_1)$ is solely a function of $S_1$ and $C_2(S_2)$ solely of $S_2$, each can be minimized separately. The argument below focuses on $C_1(S_1)$; the same reasoning applies to $C_2(S_2)$.

\medskip
\noindent
\textit{Step 2: Cost Difference and Monotonicity.}
\[
C_1(S_1) = h_1 S_1 - h_1  \frac{\lambar{1}}{\mu - \lambar{1} - \lambar{2}} \left( 1 - \left( \frac{\lambar{1}}{\mu - \lambar{2}} \right)^{S_1} \right) + c_1  \frac{\lambar{1}}{\mu - \lambar{1} - \lambar{2}} \left( \frac{\lambar{1}}{\mu - \lambar{2}} \right)^{S_1}.
\]

Define the cost difference
\[
\Delta C_1(S_1)
=
C_1(S_1+1) - C_1(S_1).
\]

A direct computation shows
\begin{align*}
\Delta C_1(S_1) = h_1 &- \frac{\lambar{1}}{\mu - \lambar{1} - \lambar{2}} \Bigg[h_1\left(\left( 1 - \left( \frac{\lambar{1}}{\mu - \lambar{2}} \right)^{S_1+1} \right)-\left( 1 - \left( \frac{\lambar{1}}{\mu - \lambar{2}} \right)^{S_1} \right)\right)\\
&-c_1\left(\left( 1 - \left( \frac{\lambar{1}}{\mu - \lambar{2}} \right)^{S_1+1} \right) - \left( 1 - \left( \frac{\lambar{1}}{\mu - \lambar{2}} \right)^{S_1} \right)\right)\Bigg]\\
&= h_1 - \frac{\lambar{1}}{\mu - \lambar{1} - \lambar{2}}  \left( \frac{\lambar{1}}{\mu - \lambar{2}} \right)^{S_1}  (h_1+c_1)  \frac{\mu - \lambar{1} - \lambar{2}}{\mu-\lambar{2}}\\
& = h_1 - (h_1+c_1)  \left( \frac{\lambar{1}}{\mu - \lambar{2}} \right)^{S_1+1}.
\end{align*}

Moreover, $\Delta C_1(S_1)$ is strictly increasing in $S_1$. Indeed, 
\[
\Delta C_1(S_1+1) - \Delta C_1(S_1)
=
\left(h_1 + c_1\right)
\left( \frac{\lambar{1}}{\mu - \lambar{2}} \right)^{S_1+1}
\frac{\mu-\lambar{1}- \lambar{2}}{\mu - \lambar{2}} 
>0.
\]

\medskip
\noindent
\textit{Step 3: Identifying the Optimal $S_1^*$.}

To find the cost-minimizing $S_1$, examine the cost difference, $\Delta C_1(S_1)$ and choose the smallest integer $S_i$ such that $\Delta C_i(S_i) \ge 0$. For any $S_1 < S_1^*$, $\Delta C_1(S_1)<0$, so increasing $S_1$ lowers cost. Once $S_1\ge S_1^*$ makes $\Delta C_1(S_1)\ge 0$, further increments do not reduce cost. Because $\Delta C_1(S_1)$ is strictly increasing, there is exactly one integer $S_1^*$ satisfying
\[
\Delta C_1(S_1^*) \ge0,
\quad
\Delta C_1(S_1^*-1) <0.
\]

Rewriting $\Delta C_1(S_1^*)\ge 0$ shows
\begin{align*}
&h_1 - (h_1+c_1)  \left( \frac{\lambar{1}}{\mu - \lambar{2}} \right)^{S_1^*+1} \ge 0 
;\qquad
\left( \frac{\lambar{1}}{\mu - \lambar{2}} \right)^{S_1^*+1} \le \frac{h_1}{h_1 + c_1};\\
& S_1^* \ge \frac{\log\left( \frac{h_1}{h_1 + c_1} \right)}{\log\left(\frac{\lambar{1}}{\mu-\lambar{2}}\right)}-1.
\end{align*}

Hence the minimal integer solution is 
\[
S_1^* =  \left\lceil \frac{\log\left( \frac{h_1}{h_1 + c_1} \right)}{\log\left(\frac{\lambar{1}}{\mu-\lambar{2}}\right)} \right\rceil-1.
\]

An identical argument applies to $C_2(S_2)$.
\end{proof}

\subsection{Optimization Within Subdomains}\label{app:optimal-joining-rates}

In each subdomain $\mathcal{D}_{s_1,s_2}$, the target inventory levels remain fixed at $S_1^*=s_1$ and $S_2^*=s_2$. Substituting $S_1=s_1$ and $S_2=s_2$ into the social welfare function
\[
\SWe(\lambar{1},\lambar{2})
=
\sum_{i=1}^2
\left[
  \lambar{i} R_i
  - h_is_i
  + h_i\tfrac{\lambar{i}}{\mu - \lambar{1} - \lambar{2}}
    \left(1 - \left(\tfrac{\lambar{i}}{\mu - \lambar{j}}\right)^{s_i}\right)
  - c_i\tfrac{\lambar{i}}{\mu - \lambar{1} - \lambar{2}}
    \left(\tfrac{\lambar{i}}{\mu - \lambar{j}}\right)^{s_i}
\right]
\]
produces a smooth function of $\lambar{1},\lambar{2}$ in $\mathcal{D}_{s_1,s_2}$.  

Although setting 
$\partial \SWe/\partial \lambar{1} = 0$ 
and 
$\partial \SWe/\partial \lambar{2} = 0$
may not yield a closed-form solution, standard methods (e.g.\ gradient-based algorithms) can find a local maximum. Explicit partial derivatives are readily computed. For instance:
\begin{align*}
&\frac{\partial \SWe}{\partial \lambar{1}}
= 
R_1
 + 
\frac{1}{(\mu - \lambar{1} - \lambar{2})^2}
\bigl[
  h_1(\mu - \lambar{2}) 
  + h_2\lambar{2}
  - (c_1+h_1)\bigl(\mu - \lambar{2} 
  \\& + s_1(\mu - \lambar{1} - \lambar{2})\bigr)
    \left(\tfrac{\lambar{1}}{\mu - \lambar{2}}\right)^{s_1} 
   - (c_2+h_2)\left(\mu - \lambar{1} + s_2(\mu - \lambar{1} - \lambar{2})\right)
    \left(\tfrac{\lambar{2}}{\mu - \lambar{1}}\right)^{s_2+1}
\bigr],
\end{align*}
(and similarly for $\partial \SWe/\partial \lambar{2}$).

\section{Supplementary Figures}\label{app:additional-numerics}

This appendix collects additional figures from the numerical experiment in Section~\ref{sec:numerics}. We present two parameter configurations to explore how holding cost asymmetries affect system behavior.

\subsection{Baseline Configuration}\label{app:baseline-figs}

The figures in this subsection use the baseline parameters from Section~\ref{sec:numerics}. These results complement the main text by showing additional metrics --- effective utilization, expected waiting times, and component-wise comparison between centralized and decentralized solutions --- across the full $(\kappa, \rho)$ parameter space.
\begin{figure}[!htbp]
\centering
\includegraphics[width=\textwidth]{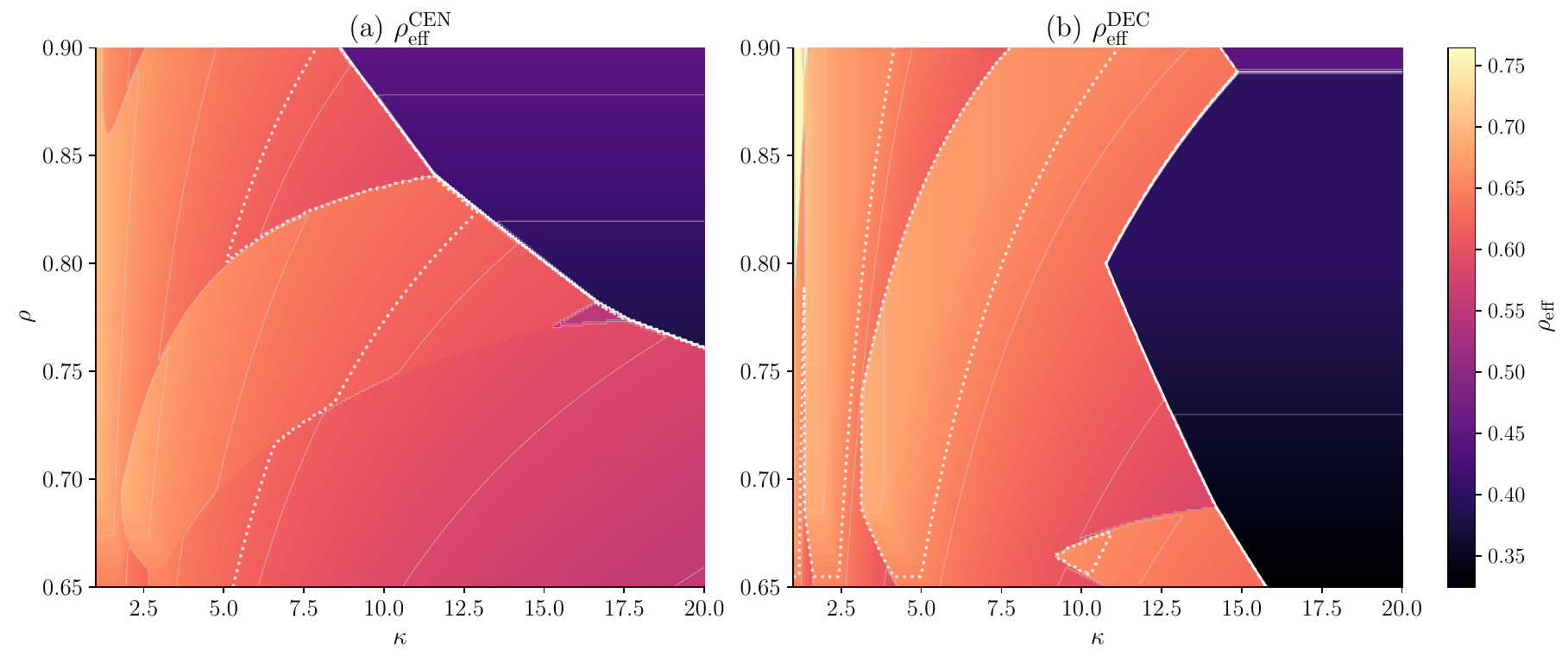}
\caption{Effective utilization $\rho_{\mathrm{eff}} = (q_1\Lambda_1 + q_2\Lambda_2)/\mu$ over the $(\kappa,\rho)$ plane. White contours at 25\% (solid), 50\% (dashed), 75\% (dotted) of range; thin gray contours provide detail.}
\label{fig:app-rho-eff}
\end{figure}

\begin{figure}[!htbp]
\centering
\includegraphics[width=0.8\textwidth]{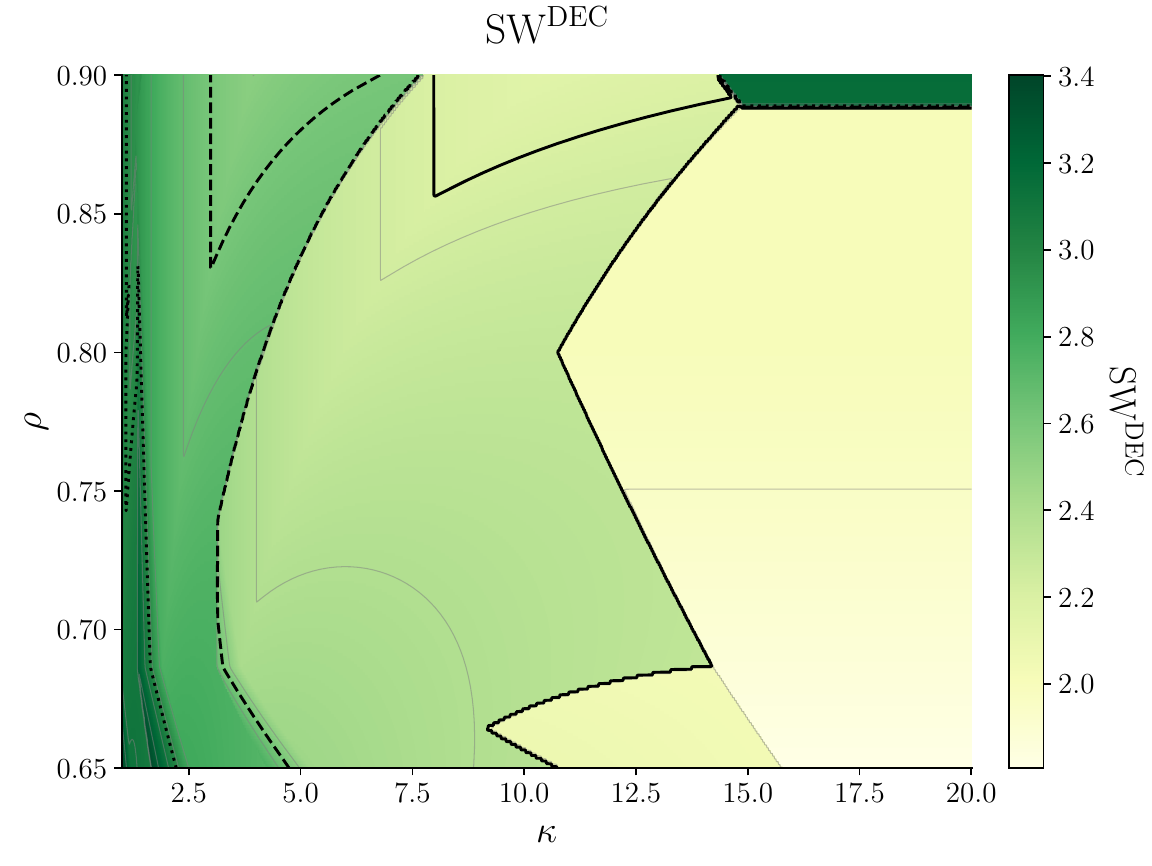}
\caption{Social welfare $\mathrm{SW}^{\textsc{dec}}$ under decentralization over the $(\kappa,\rho)$ plane. Thick black contours mark 2.20, 2.60, 3.00 (25\%, 50\%, 75\% of range); thin gray contours provide detail.}
\label{fig:app-welfare-dec}
\end{figure}

\begin{figure}[!htbp]
\centering
\includegraphics[width=\textwidth]{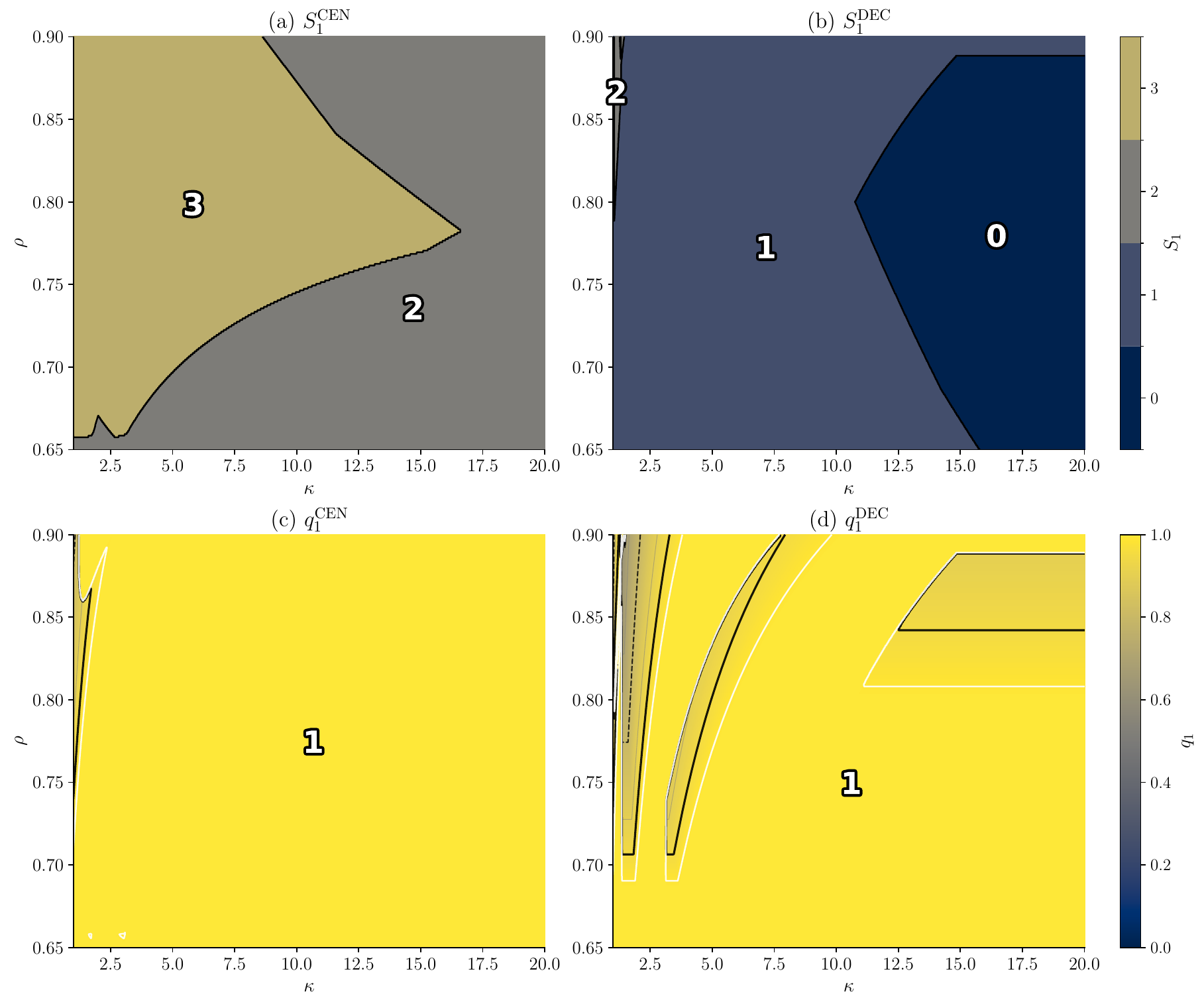}
\caption{Type-1 decisions over the $(\kappa,\rho)$ plane. Top row: inventory $S_1$; contours mark integer transitions. Bottom row: joining probability $q_1$; contours at $0.80$ (dashed black), $0.95$ (thick black), and $0.99$ (thick white); thin gray contours provide detail.}
\label{fig:app-type1-response}
\end{figure}

\begin{figure}[!htbp]
\centering
\includegraphics[width=\textwidth]{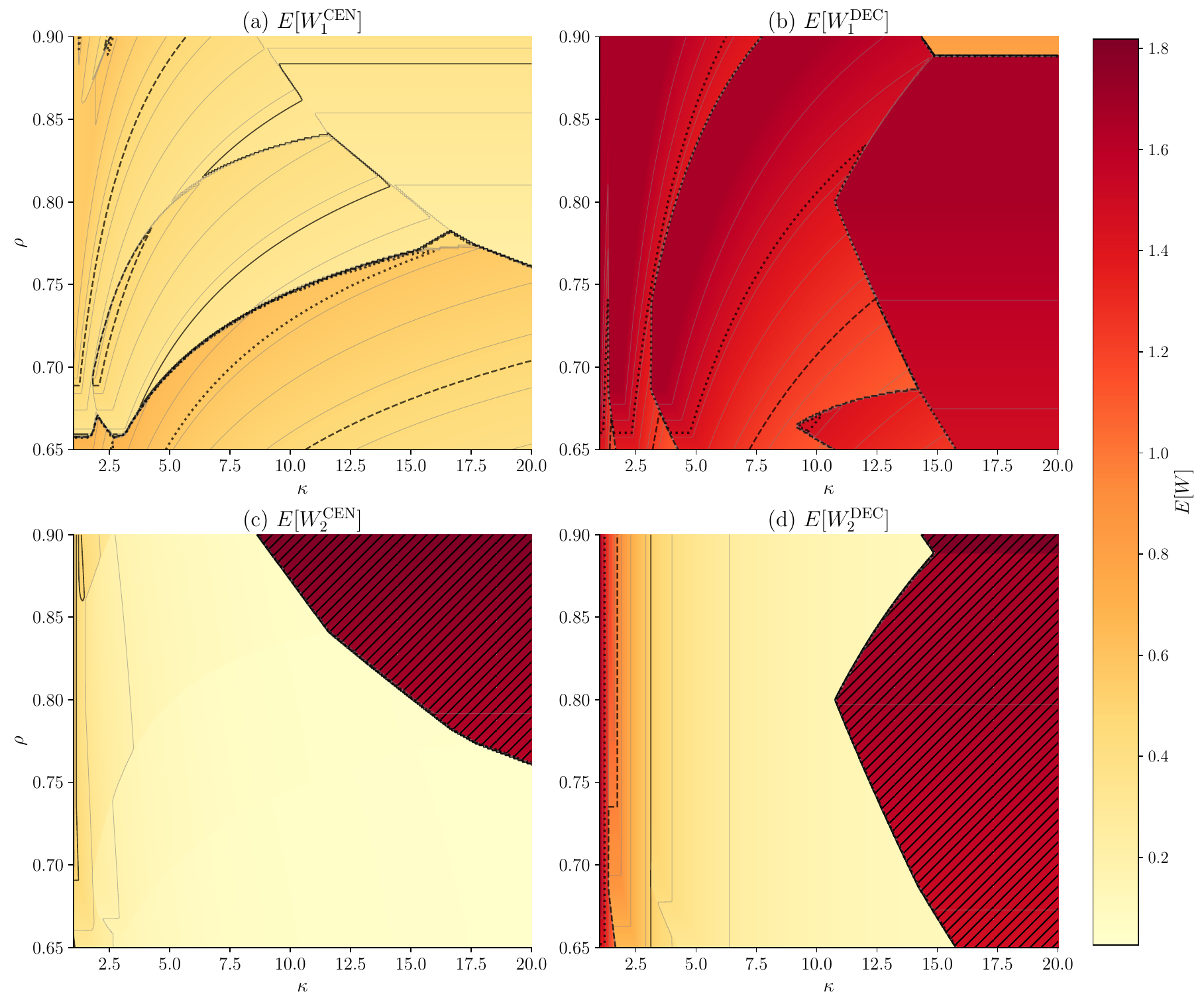}
\caption{Expected waiting times $E[W_i]$ over the $(\kappa,\rho)$ plane. Top row: type-1. Bottom row: type-2. Contours at 25\% (solid), 50\% (dashed), 75\% (dotted) of each panel's range; thin gray contours provide detail. Hatching indicates no joining.}
\label{fig:app-waiting-times}
\end{figure}

\begin{figure}[!htbp]
\centering
\includegraphics[width=\textwidth]{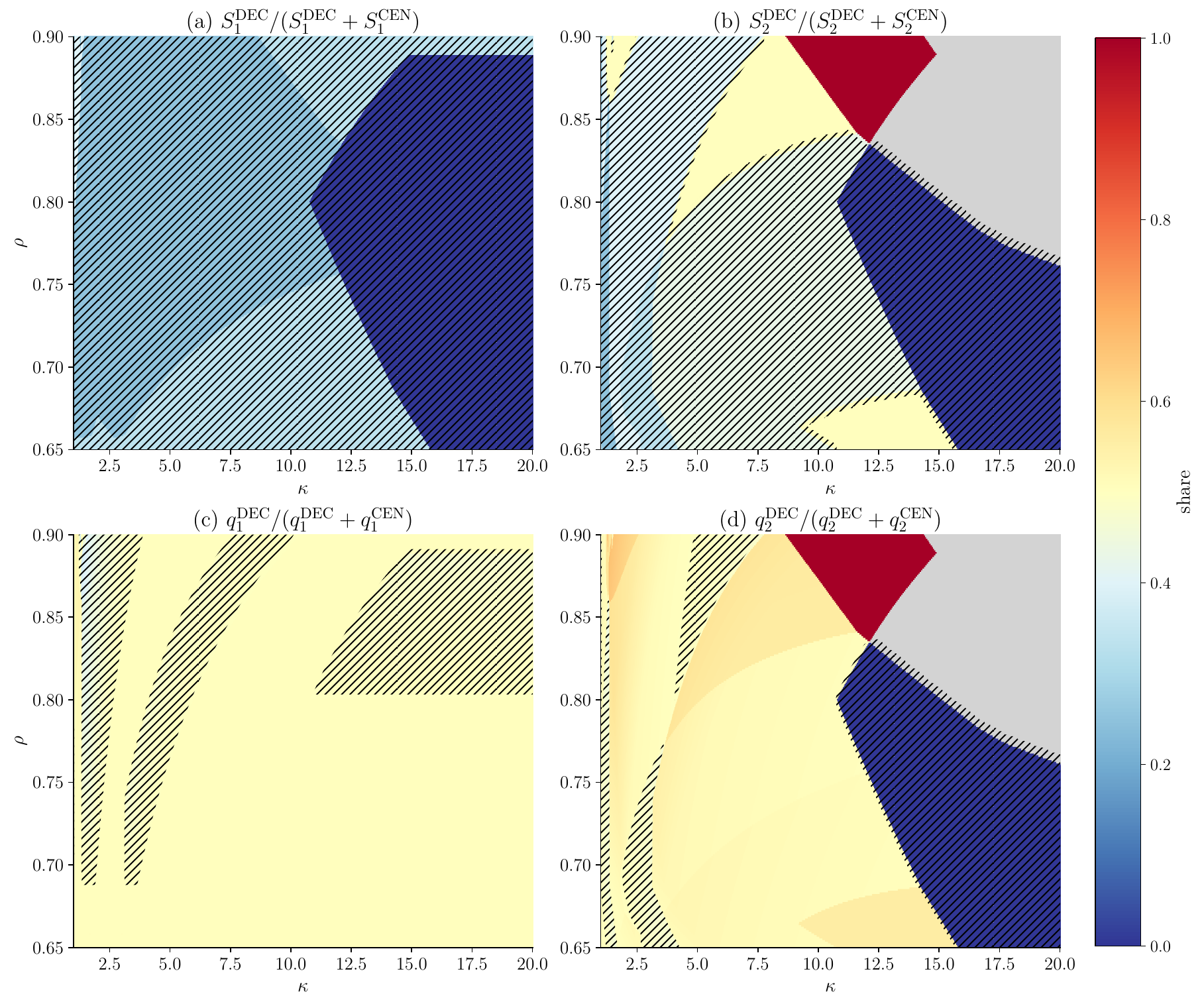}
\caption{Decentralized shares $X^{\textsc{dec}}/(X^{\textsc{dec}}+X^{\textsc{cen}})$ for $X \in \{S_1, S_2, q_1, q_2\}$. Top row: inventory shares. Bottom row: joining probability shares. Diverging colormap centered at $0.5$: blue indicates higher under centralization, red indicates higher under decentralization. Hatching where centralization dominates (share $< 0.5$). Light gray regions indicate undefined shares (zero values in both settings).}
\label{fig:app-dec-shares}
\end{figure}

\begin{figure}[!htbp]
\centering
\includegraphics[width=\textwidth]{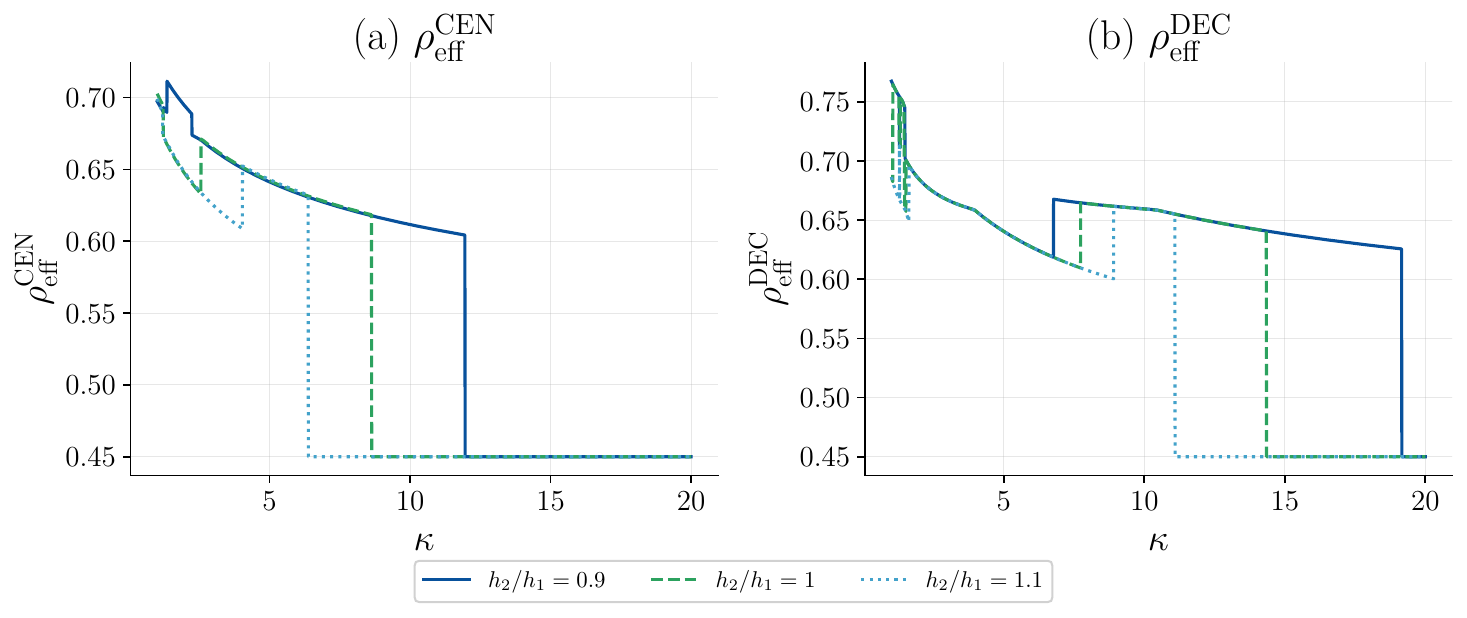}
\caption{Cross-section at $\rho = 0.9$: effective utilization $\rho_{\mathrm{eff}} = (q_1\Lambda_1 + q_2\Lambda_2)/\mu$ as function of $\kappa$ for different $h_2/h_1$ ratios. Left: centralized. Right: decentralized.}
\label{fig:cross-app-rho-eff}
\end{figure}

\begin{figure}[!htbp]
\centering
\includegraphics[width=\textwidth]{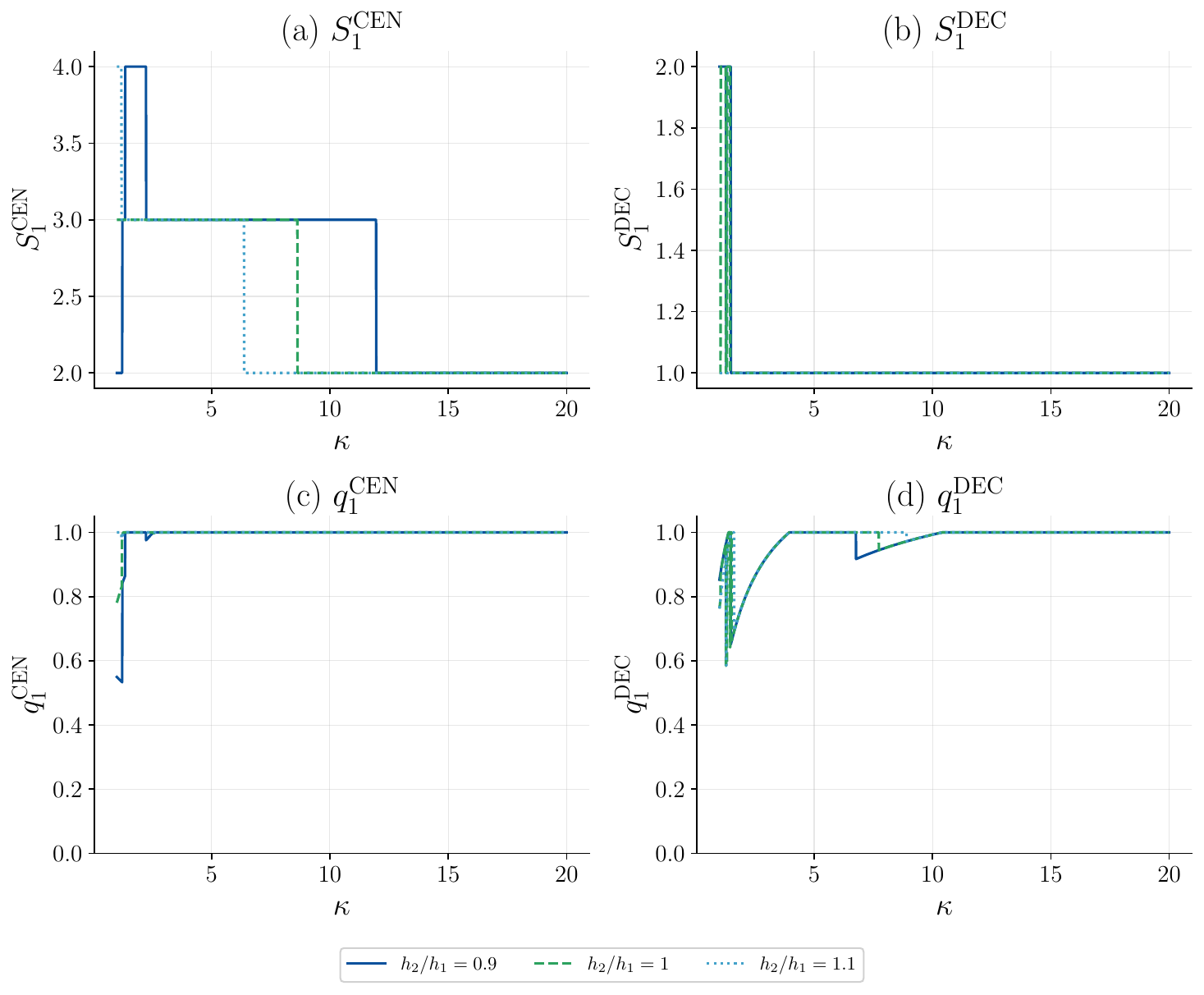}
\caption{Cross-section at $\rho = 0.9$: type-1 decisions as functions of $\kappa$ for different $h_2/h_1$ ratios. Top row: inventory $S_1$. Bottom row: joining probability $q_1$.}
\label{fig:cross-app-type1}
\end{figure}

\begin{figure}[!htbp]
\centering
\includegraphics[width=\textwidth]{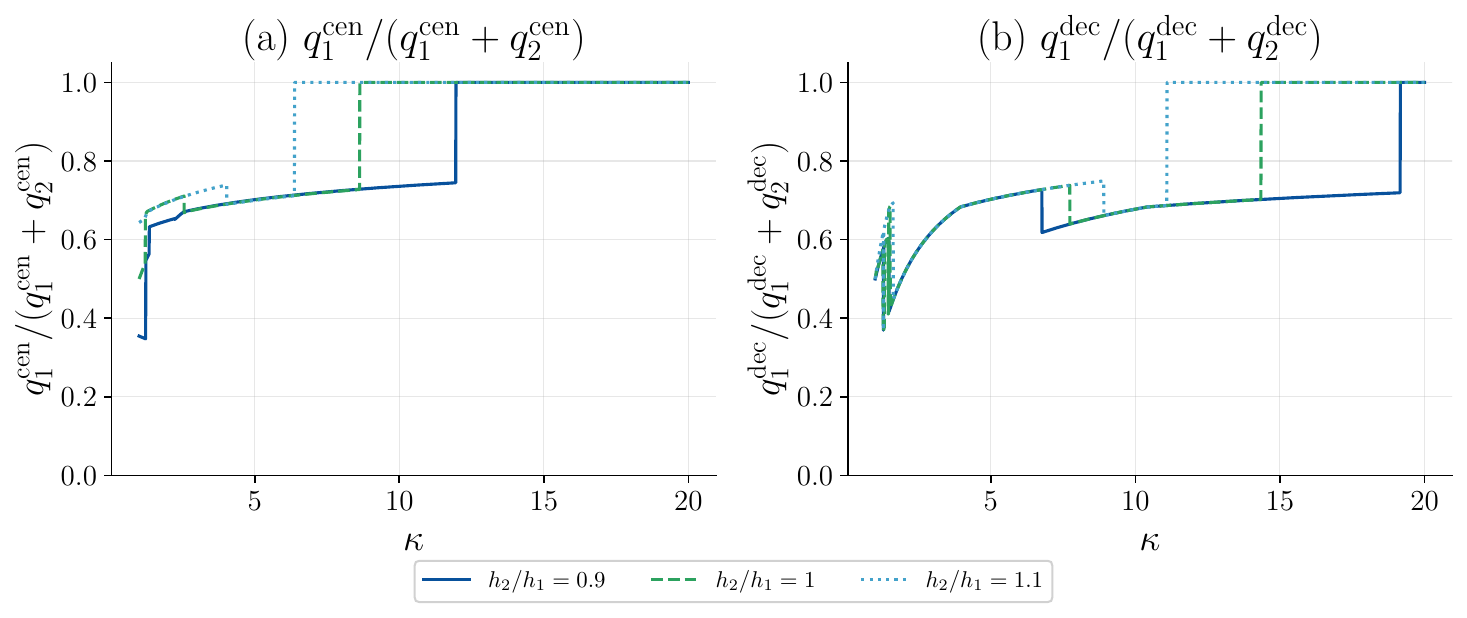}
\caption{Cross-section at $\rho = 0.9$: type-1 share $q_1/(q_1+q_2)$ as function of $\kappa$ for different $h_2/h_1$ ratios. Left: centralized. Right: decentralized.}
\label{fig:cross-app-q1share}
\end{figure}

\begin{figure}[!htbp]
\centering
\includegraphics[width=\textwidth]{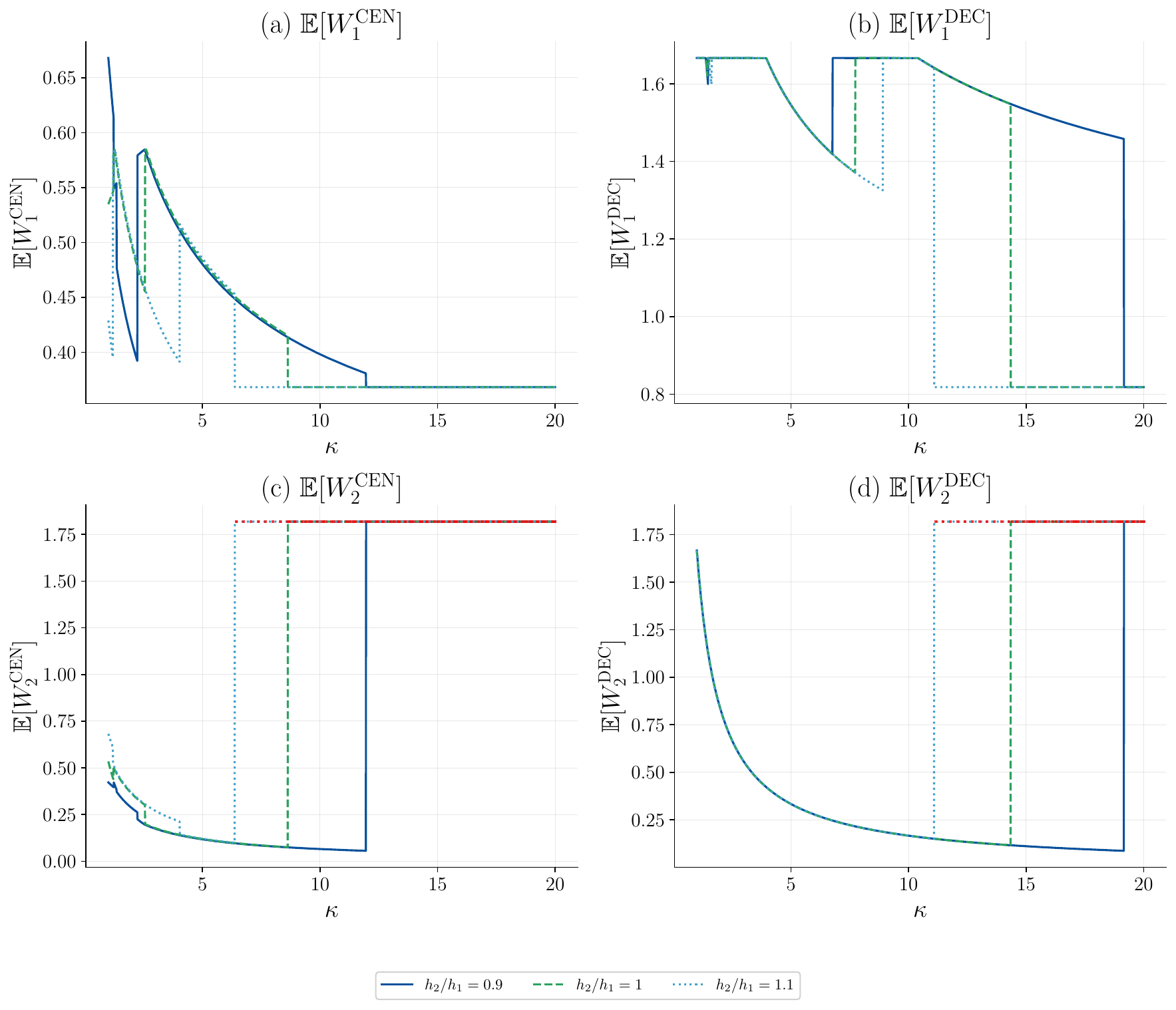}
\caption{Cross-section at $\rho = 0.9$: expected waiting times $E[W_i]$ as functions of $\kappa$ for different $h_2/h_1$ ratios. Red segments indicate no joining.}
\label{fig:cross-app-waiting}
\end{figure}

\begin{figure}[!htbp]
\centering
\includegraphics[width=\textwidth]{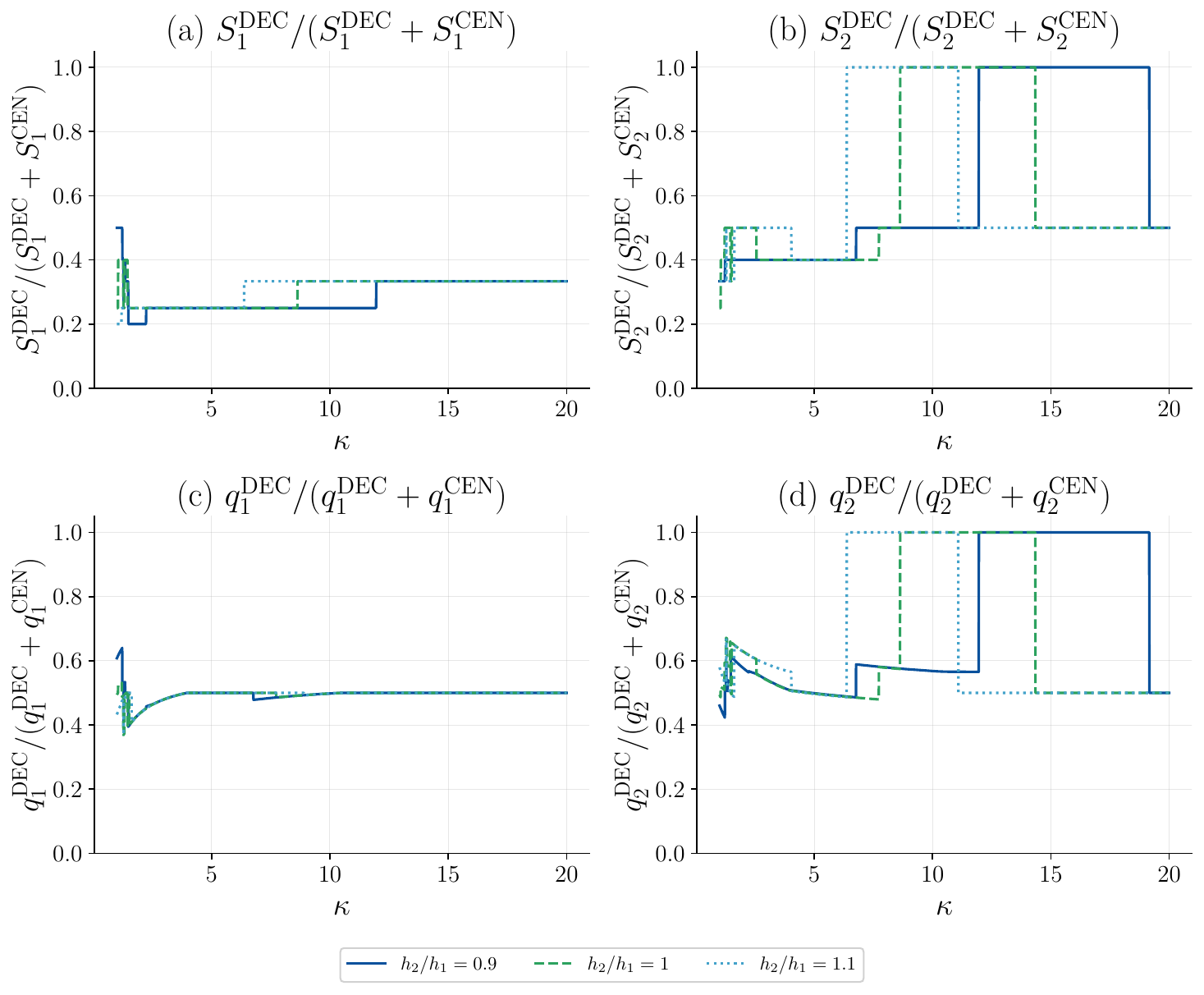}
\caption{Cross-section at $\rho = 0.9$: decentralized shares $X^{\textsc{dec}}/(X^{\textsc{dec}}+X^{\textsc{cen}})$ as functions of $\kappa$ for different $h_2/h_1$ ratios. Top row: inventory shares. Bottom row: joining probability shares.}
\label{fig:cross-app-dec-shares}
\end{figure}

\FloatBarrier

\subsection{Reduced Type-1 Holding Cost}\label{app:reduced-h1-figs}

The figures below repeat the analysis with $c_1 = 1$ and $h_1 = 0.05$ (compared to the baseline $c_1 = 3$ and $h_1 = 0.4$). The lower holding cost reduces the trade-off between retaining customers and saving on inventory storage for product~1: with smaller~$h_i$, the producer faces less pressure to limit inventory. Cross-sections are taken at $\rho = 0.8$.
\begin{figure}[!htbp]
\centering
\includegraphics[width=0.8\textwidth]{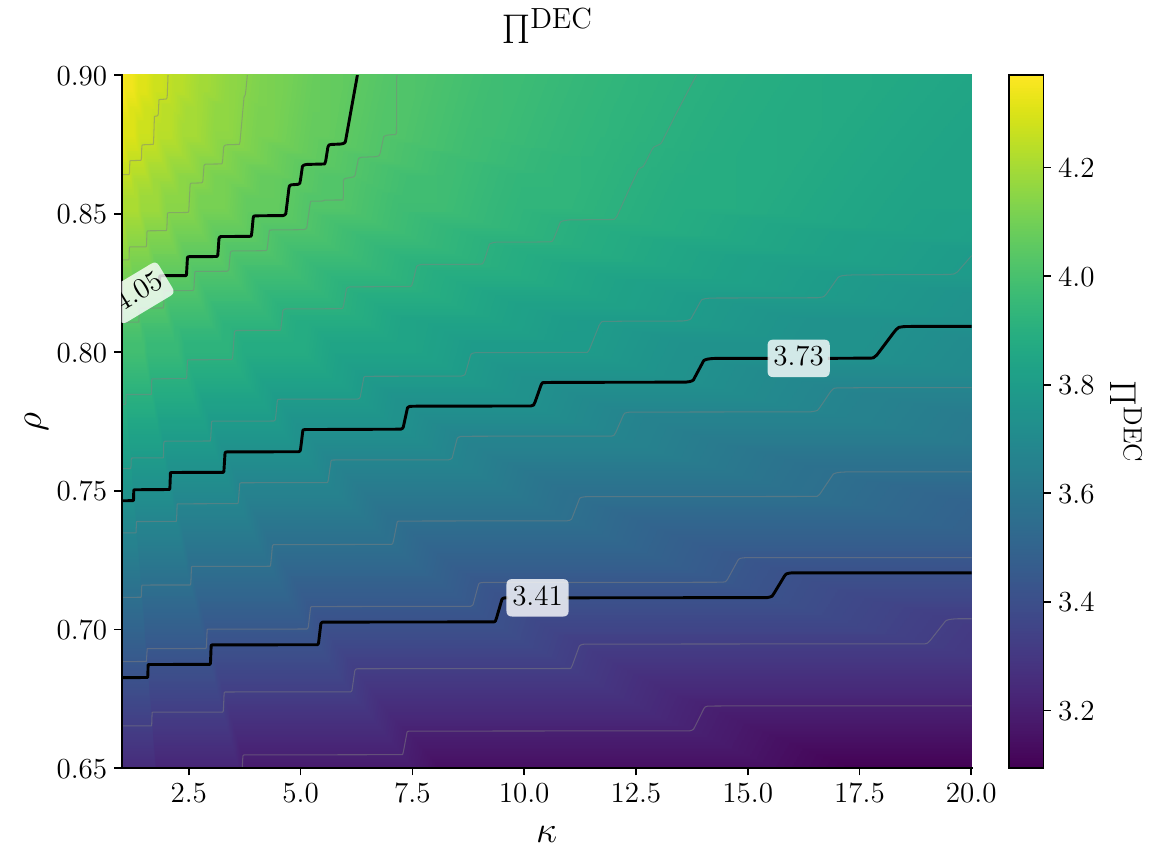}
\caption{Producer's profit $\Pi^{\textsc{dec}}$ under decentralization over the $(\kappa,\rho)$ plane (reduced $h_1$). Thick black contours mark 3.41, 3.73, 4.05 (25\%, 50\%, 75\% of range); thin gray contours provide detail.}
\label{fig:app2-profit-dec}
\end{figure}

\begin{figure}[!htbp]
\centering
\includegraphics[width=0.8\textwidth]{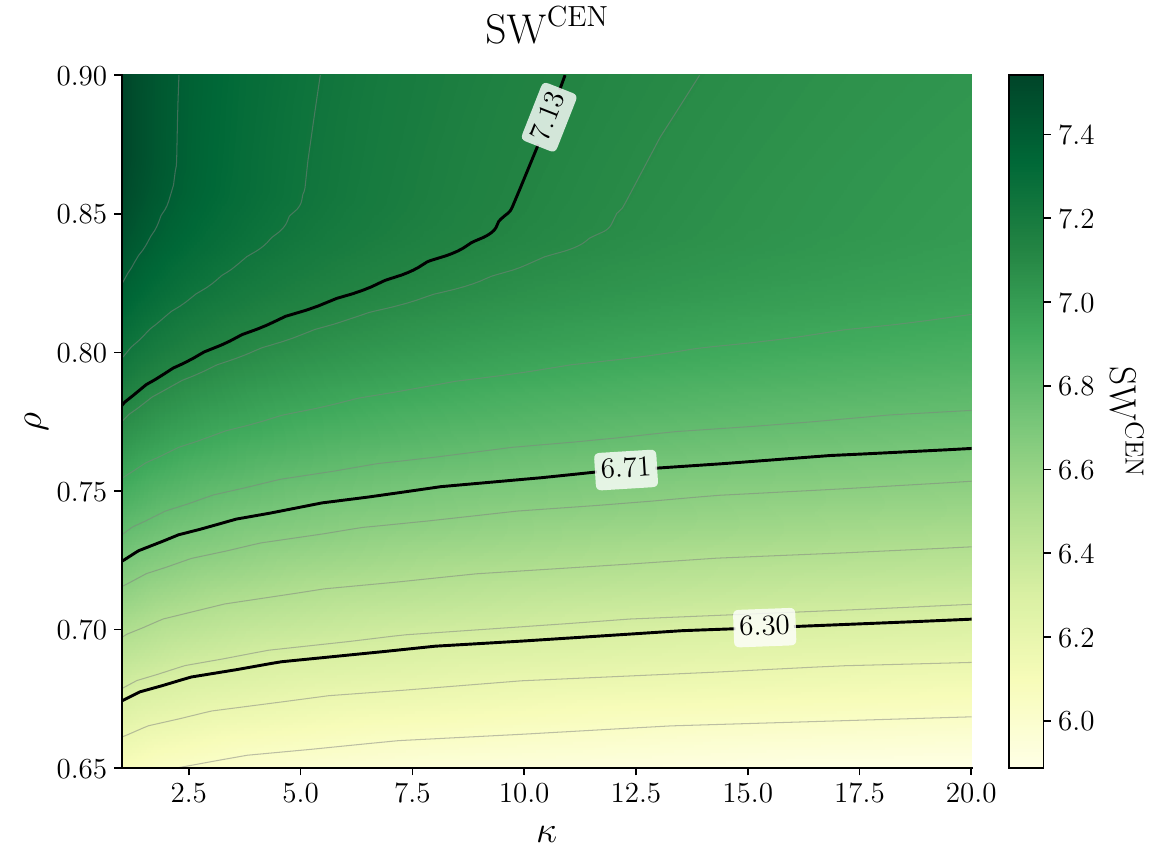}
\caption{Social welfare $\mathrm{SW}^{\textsc{cen}}$ under centralization over the $(\kappa,\rho)$ plane (reduced $h_1$). Thick black contours mark 6.30, 6.71, 7.13 (25\%, 50\%, 75\% of range); thin gray contours provide detail.}
\label{fig:app2-welfare-cen}
\end{figure}

\begin{figure}[!htbp]
\centering
\includegraphics[width=0.8\textwidth]{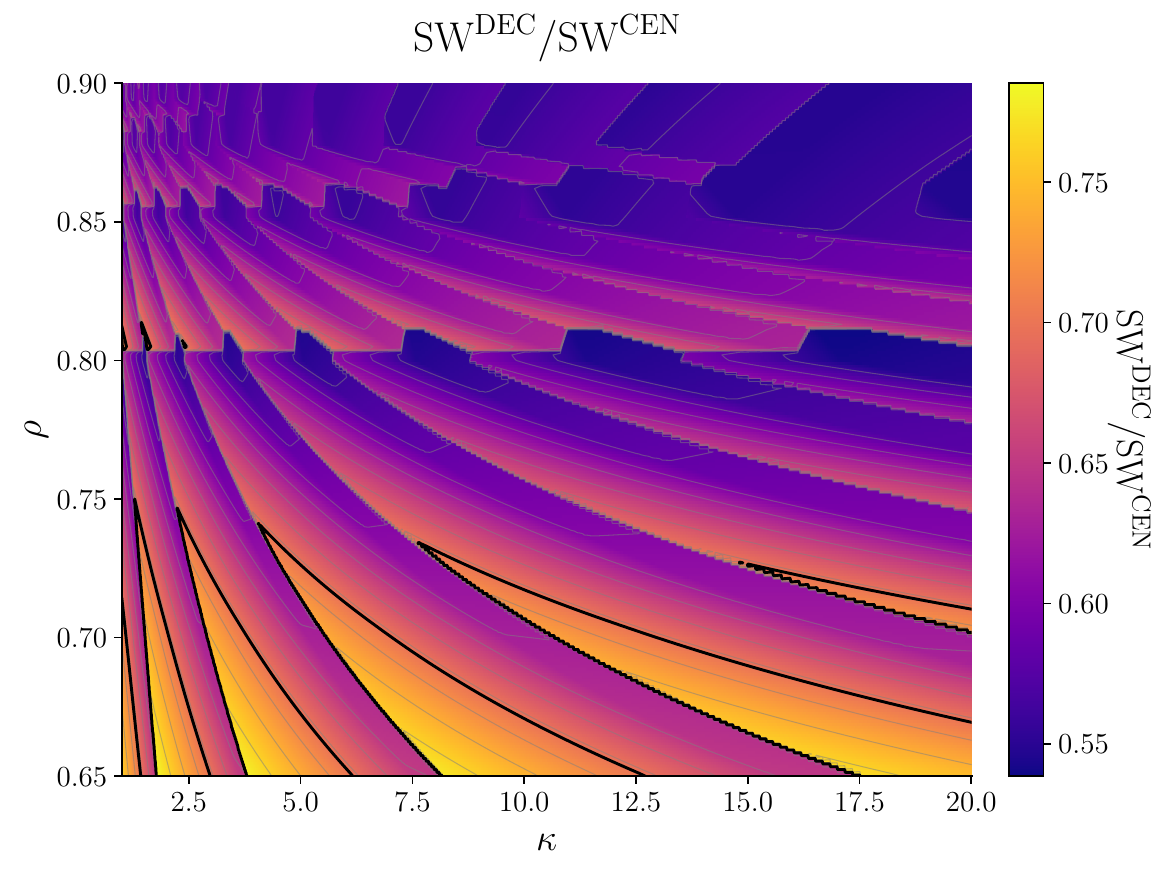}
\caption{Efficiency $\mathrm{SW}^{\textsc{dec}}/\mathrm{SW}^{\textsc{cen}}$ over the $(\kappa,\rho)$ plane (reduced $h_1$). Contours mark $0.70$ (solid), $0.80$ (dotted), and $0.90$ (dashed); thin gray contours provide detail.}
\label{fig:app2-efficiency}
\end{figure}

\begin{figure}[!htbp]
\centering
\includegraphics[width=\textwidth]{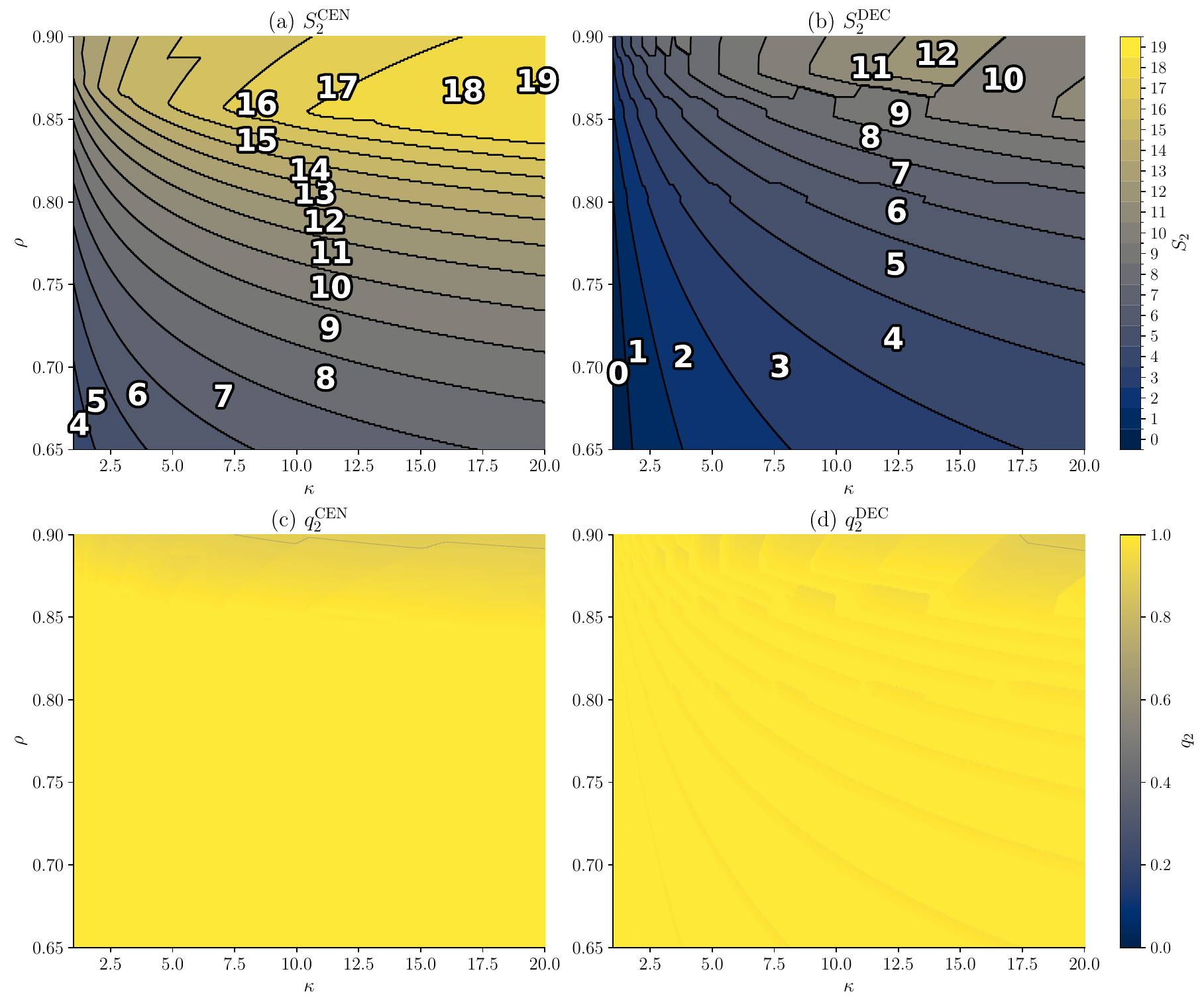}
\caption{Type-2 decisions over the $(\kappa,\rho)$ plane (reduced $h_1$). Top row: inventory $S_2$; contours mark integer transitions. Bottom row: joining probability $q_2$.}
\label{fig:app2-type2-response}
\end{figure}

\begin{figure}[!htbp]
\centering
\includegraphics[width=\textwidth]{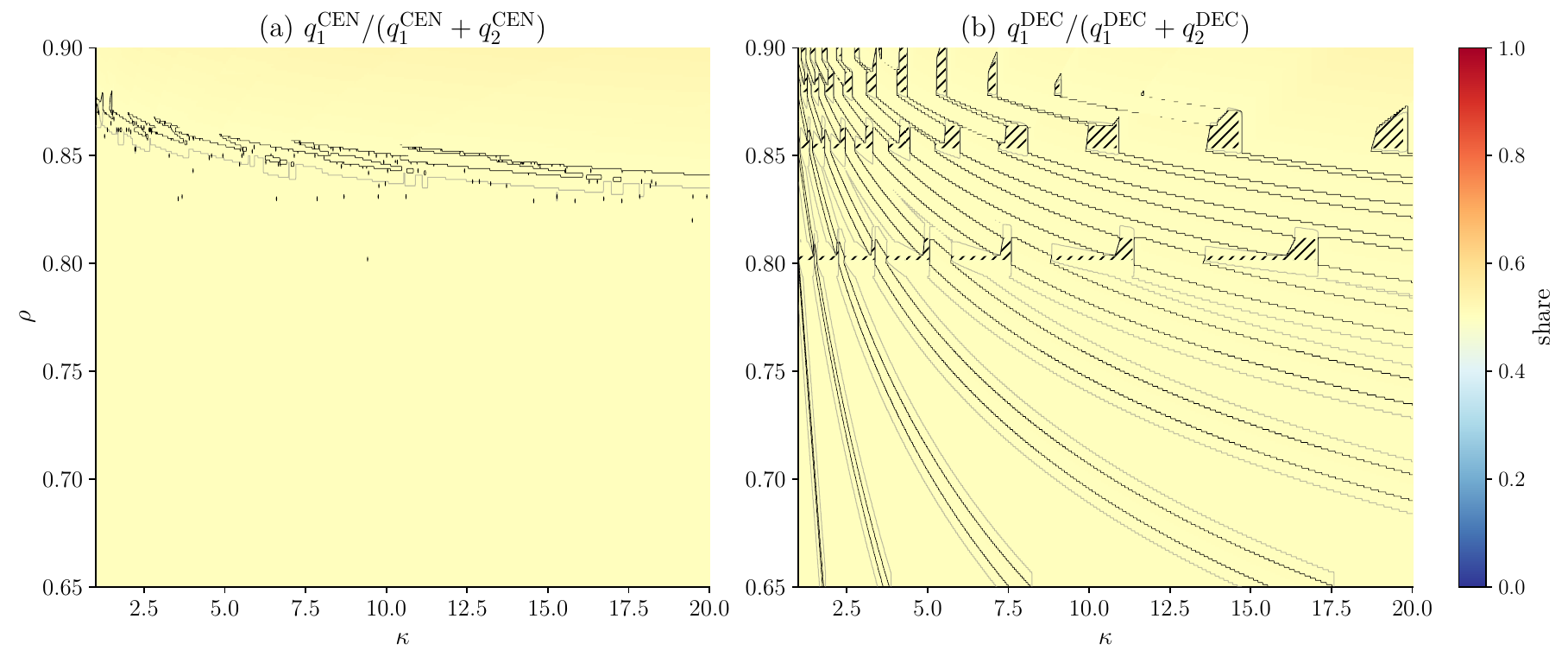}
\caption{Type-1 share $q_1/(q_1+q_2)$ over the $(\kappa,\rho)$ plane (reduced $h_1$). Contours at $0.5$ (thin solid), $0.6$ and $0.7$ (dashed). Hatching indicates type-2 majority (share $< 0.5$).}
\label{fig:app2-q1share}
\end{figure}

\begin{figure}[!htbp]
\centering
\includegraphics[width=\textwidth]{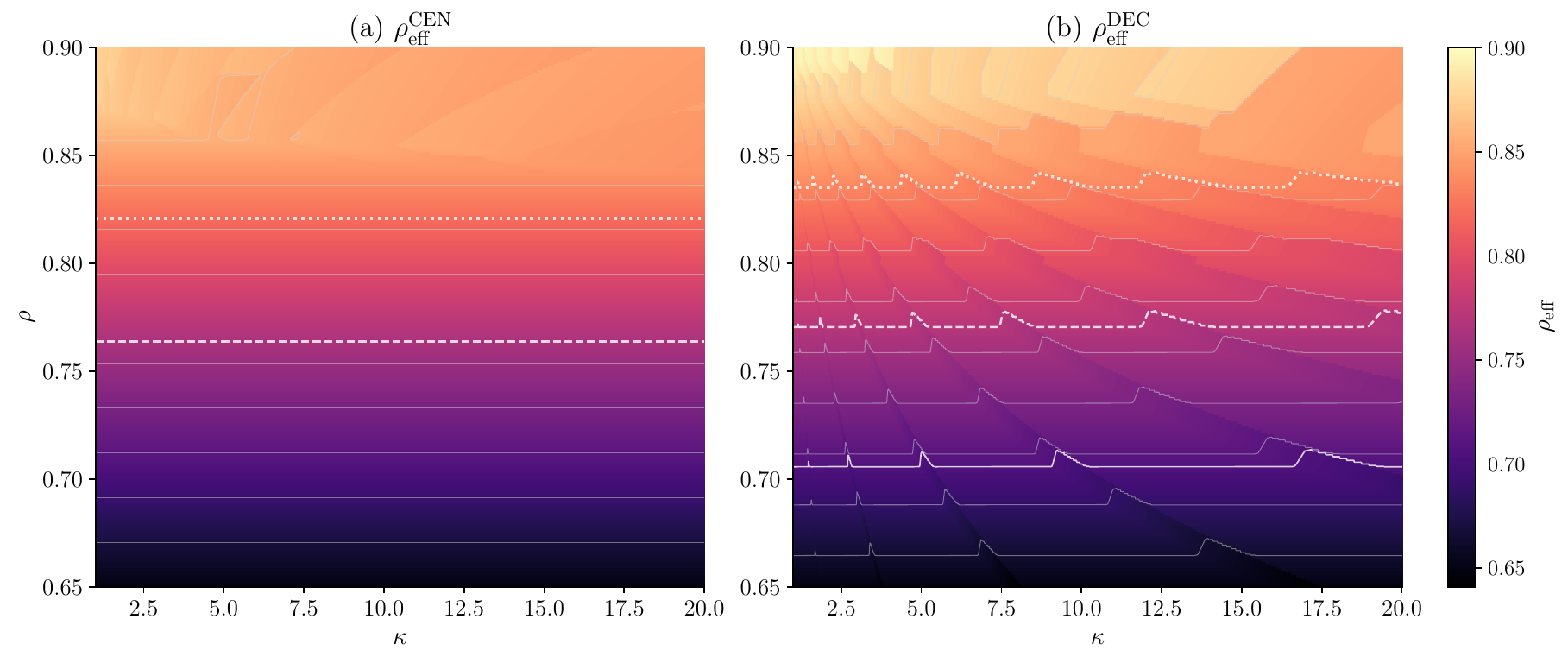}
\caption{Effective utilization $\rho_{\mathrm{eff}} = (q_1\Lambda_1 + q_2\Lambda_2)/\mu$ over the $(\kappa,\rho)$ plane (reduced $h_1$). White contours at 25\% (solid), 50\% (dashed), 75\% (dotted) of range; thin gray contours provide detail.}
\label{fig:app2-rho-eff}
\end{figure}

\begin{figure}[!htbp]
\centering
\includegraphics[width=0.8\textwidth]{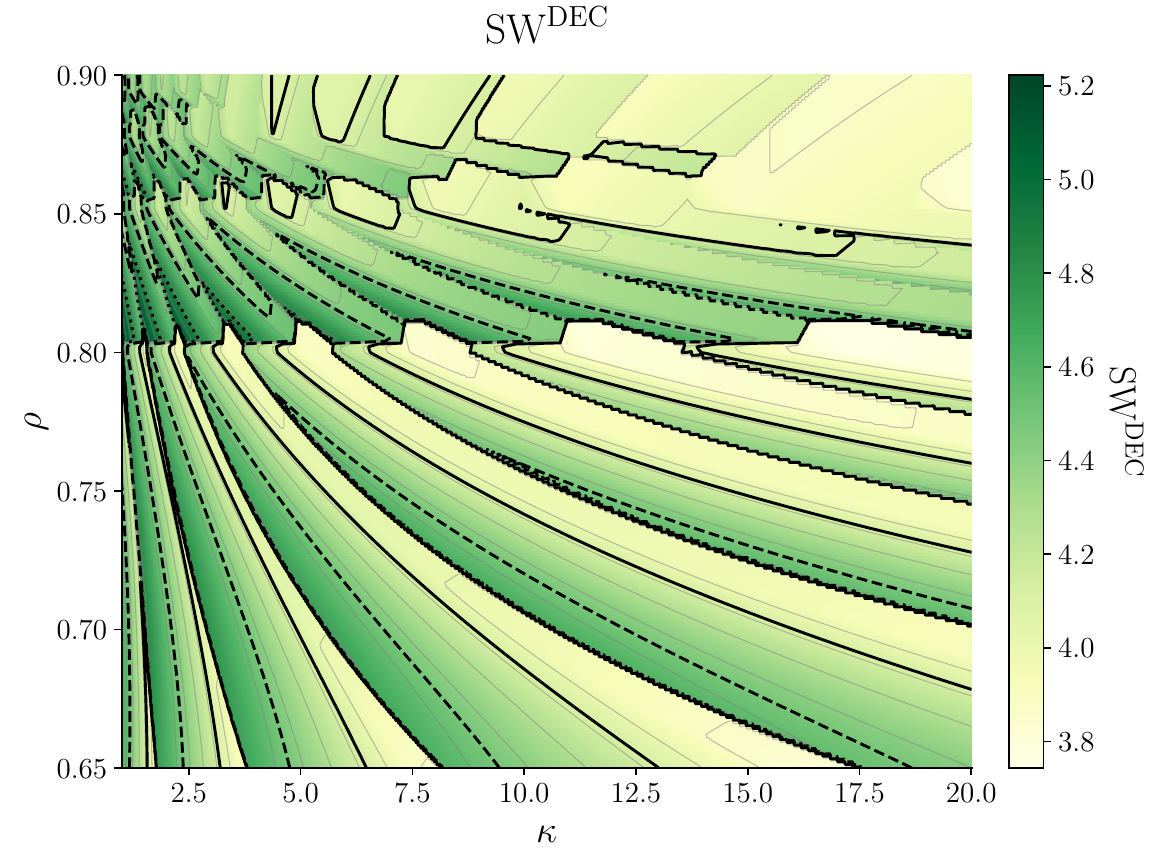}
\caption{Social welfare $\mathrm{SW}^{\textsc{dec}}$ under decentralization over the $(\kappa,\rho)$ plane (reduced $h_1$). Contours mark 4.11, 4.48, 4.85 (25\% solid, 50\% dashed, 75\% dotted); thin gray contours provide detail.}
\label{fig:app2-welfare-dec}
\end{figure}

\begin{figure}[!htbp]
\centering
\includegraphics[width=\textwidth]{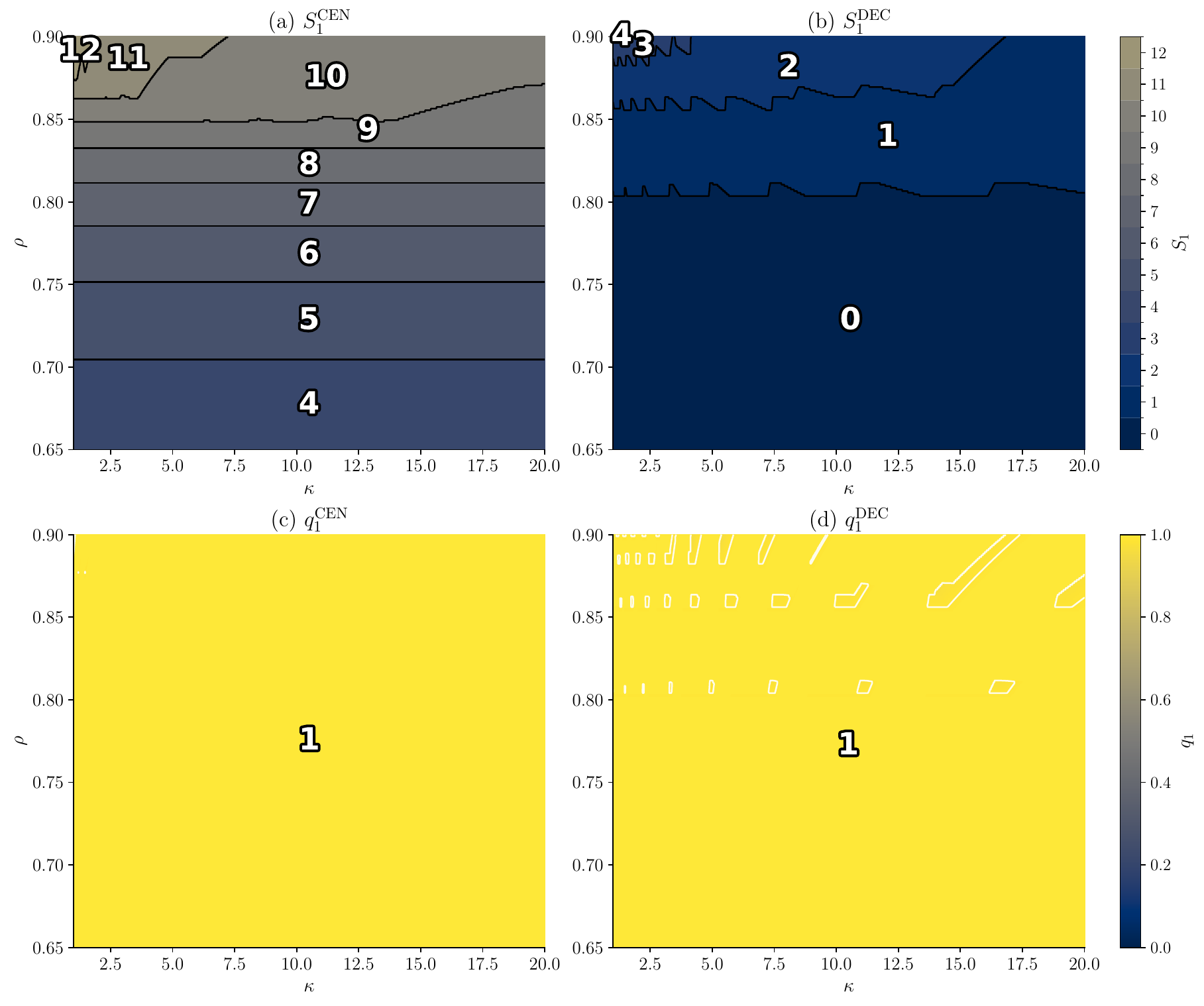}
\caption{Type-1 decisions over the $(\kappa,\rho)$ plane (reduced $h_1$). Top row: inventory $S_1$; contours mark integer transitions. Bottom row: joining probability $q_1$; contour at $0.99$ (thick white).}
\label{fig:app2-type1-response}
\end{figure}

\begin{figure}[!htbp]
\centering
\includegraphics[width=\textwidth]{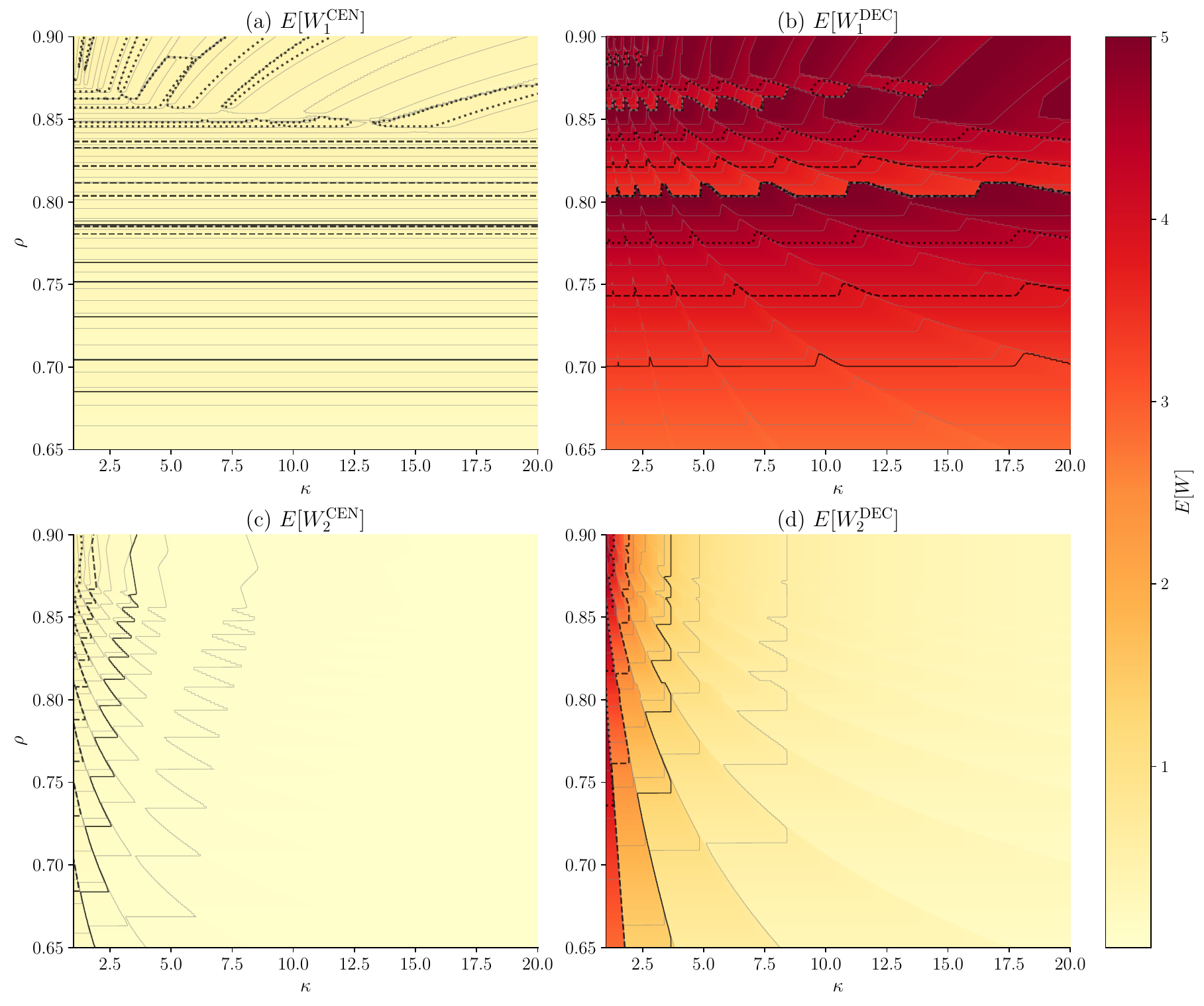}
\caption{Expected waiting times $E[W_i]$ over the $(\kappa,\rho)$ plane (reduced $h_1$). Top row: type-1. Bottom row: type-2. Contours at 25\% (solid), 50\% (dashed), 75\% (dotted) of each panel's range.}
\label{fig:app2-waiting-times}
\end{figure}

\begin{figure}[!htbp]
\centering
\includegraphics[width=\textwidth]{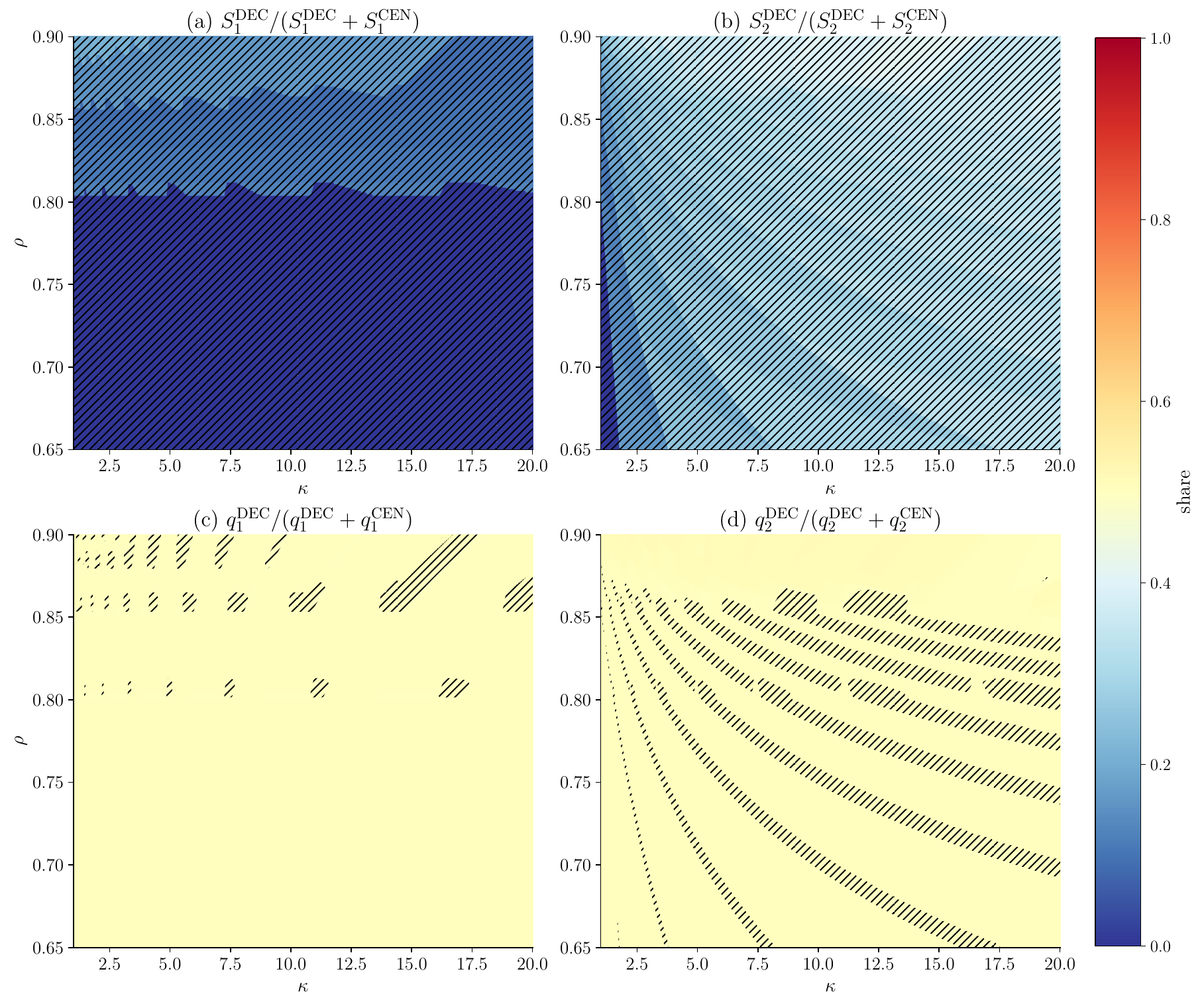}
\caption{Decentralized shares $X^{\textsc{dec}}/(X^{\textsc{dec}}+X^{\textsc{cen}})$ for $X \in \{S_1, S_2, q_1, q_2\}$ (reduced $h_1$). Top row: inventory shares. Bottom row: joining probability shares. Diverging colormap centered at $0.5$: blue indicates higher under centralization, red indicates higher under decentralization (not present). Hatching where centralization strictly dominates (share $< 0.5$).}
\label{fig:app2-dec-shares}
\end{figure}

\begin{figure}[!htbp]
\centering
\includegraphics[width=\textwidth]{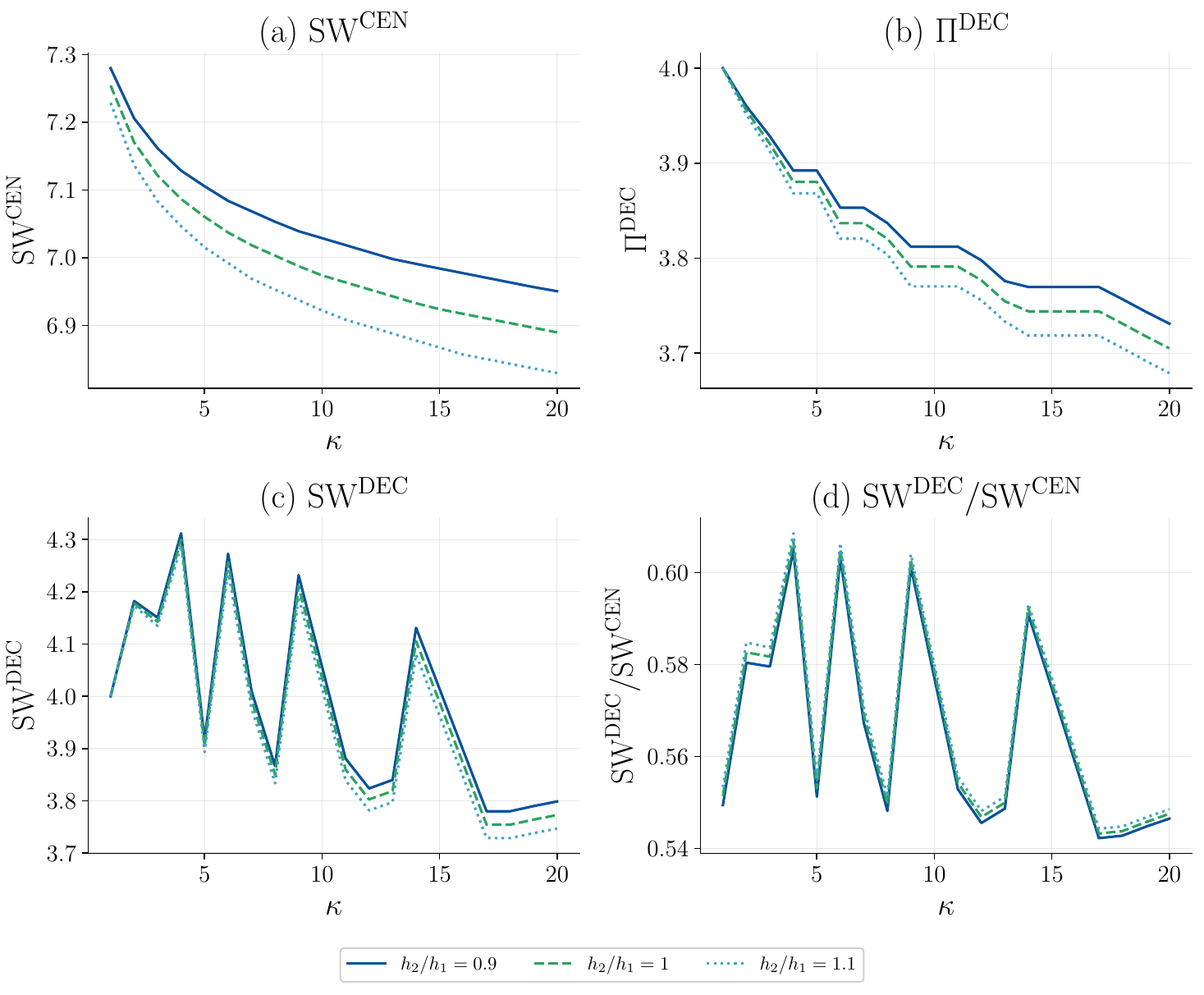}
\caption{Cross-section at $\rho = 0.8$ as functions of $\kappa$ for different $h_2/h_1$ ratios (reduced $h_1$). Top row: $\mathrm{SW}^{\textsc{cen}}$ (left) and $\Pi^{\textsc{dec}}$ (right). Bottom row: $\mathrm{SW}^{\textsc{dec}}$ (left) and efficiency ratio (right).}
\label{fig:cross-app2-welf-profit}
\end{figure}

\begin{figure}[!htbp]
\centering
\includegraphics[width=\textwidth]{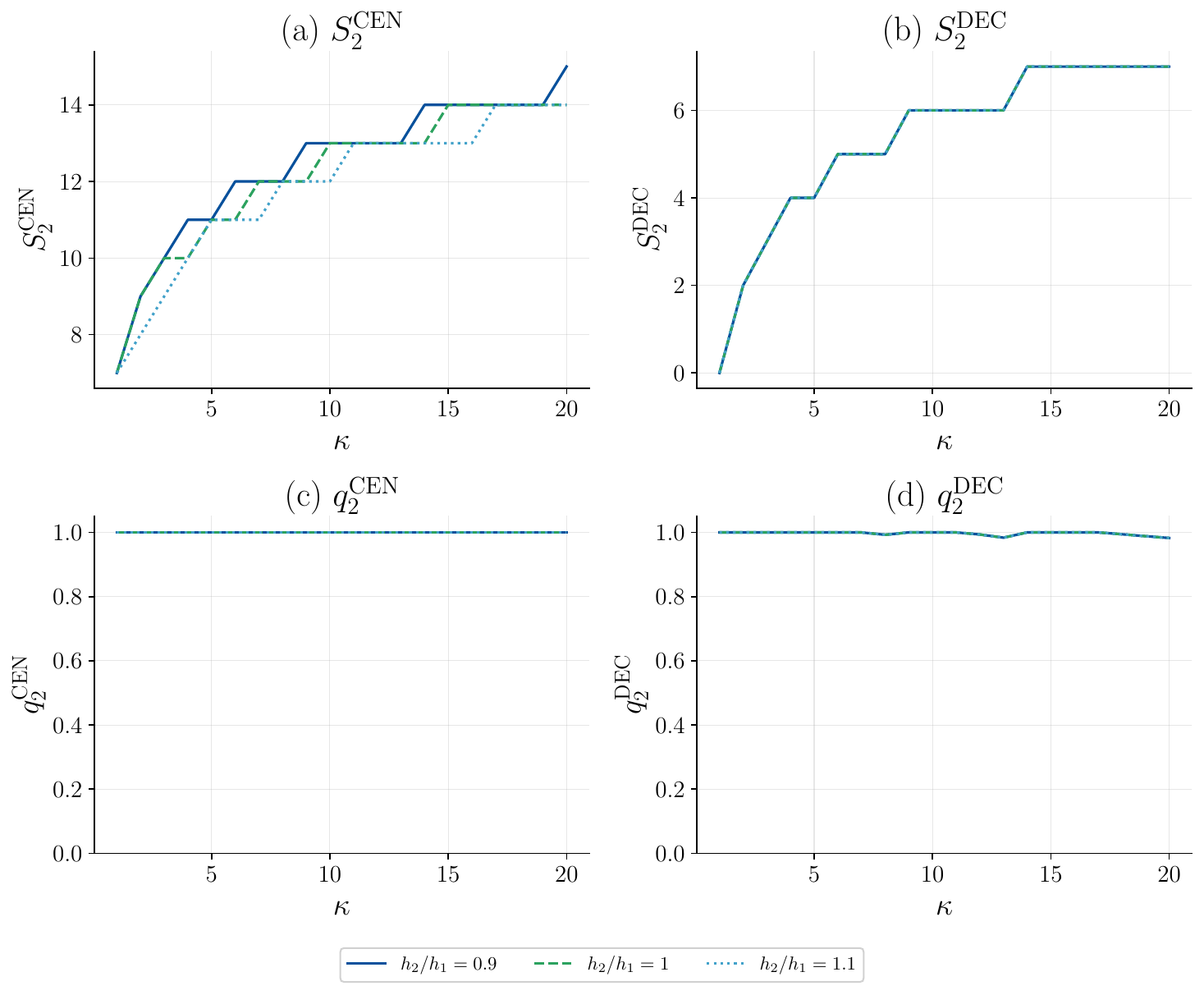}
\caption{Cross-section at $\rho = 0.8$: type-2 decisions as functions of $\kappa$ for different $h_2/h_1$ ratios (reduced $h_1$). Top row: inventory $S_2$. Bottom row: joining probability $q_2$.}
\label{fig:cross-app2-type2}
\end{figure}


\begin{figure}[!htbp]
\centering
\includegraphics[width=\textwidth]{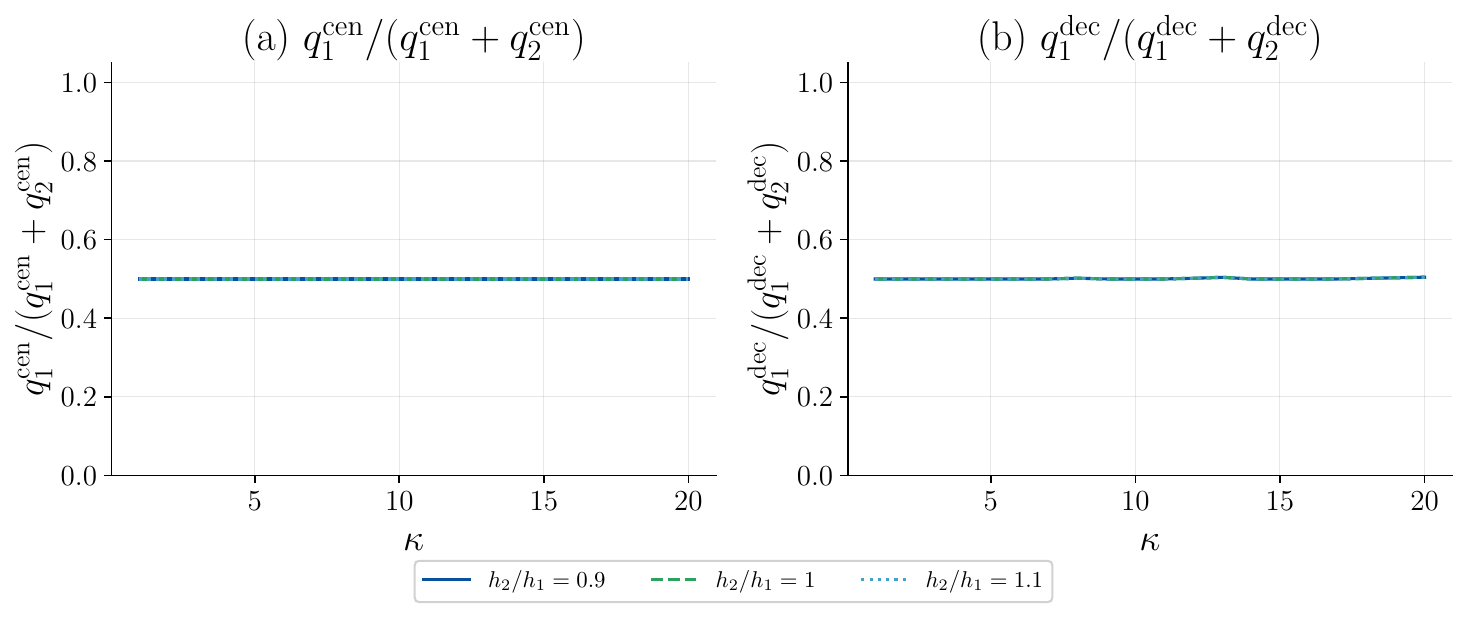}
\caption{Cross-section at $\rho = 0.8$: type-1 share $q_1/(q_1+q_2)$ as function of $\kappa$ for different $h_2/h_1$ ratios (reduced $h_1$). Left: centralized. Right: decentralized.}
\label{fig:cross-app2-q1share}
\end{figure}

\begin{figure}[!htbp]
\centering
\includegraphics[width=\textwidth]{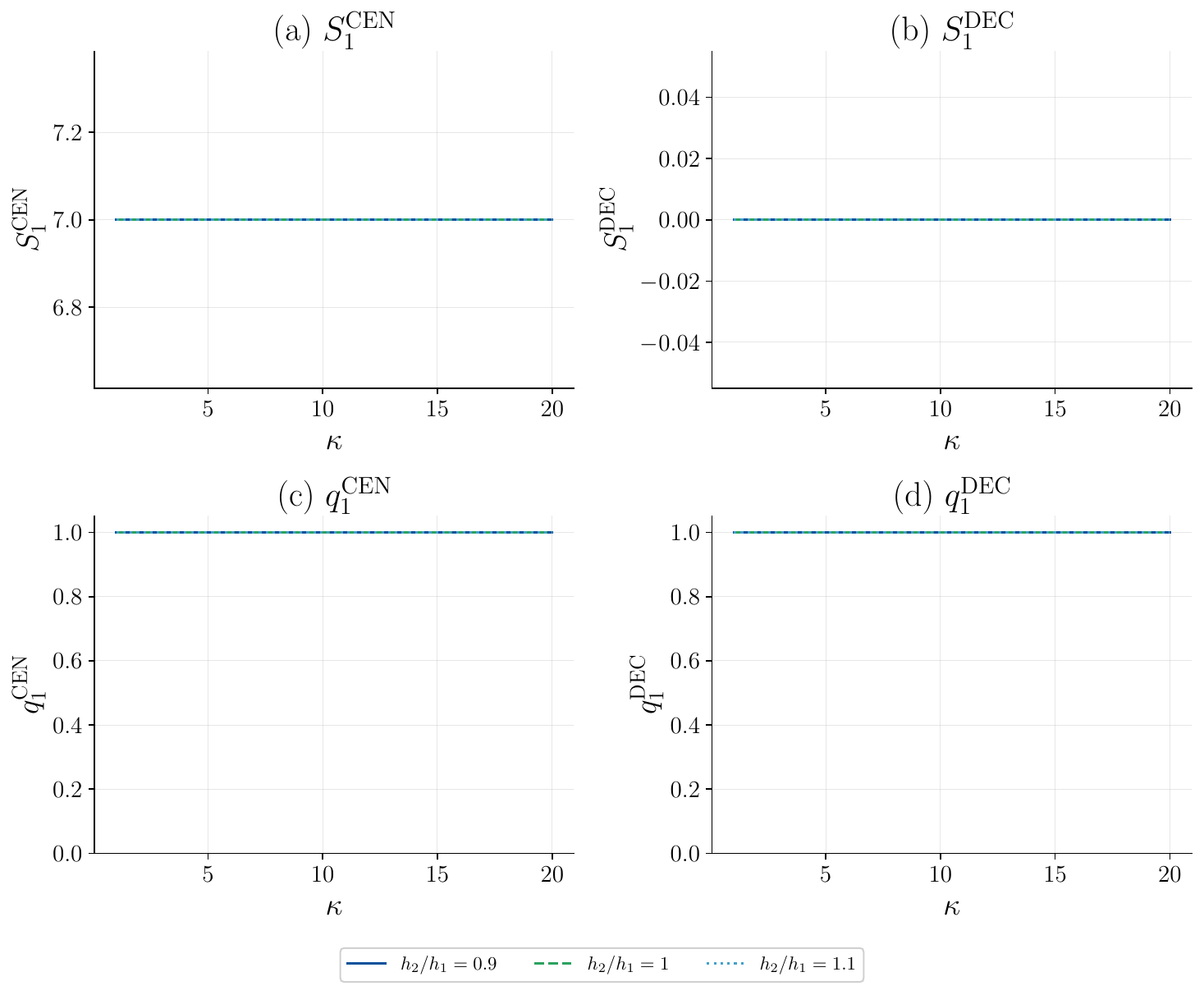}
\caption{Cross-section at $\rho = 0.8$: type-1 decisions as functions of $\kappa$ for different $h_2/h_1$ ratios (reduced $h_1$). Top row: inventory $S_1$. Bottom row: joining probability $q_1$.}
\label{fig:cross-app2-type1}
\end{figure}

\begin{figure}[!htbp]
\centering
\includegraphics[width=\textwidth]{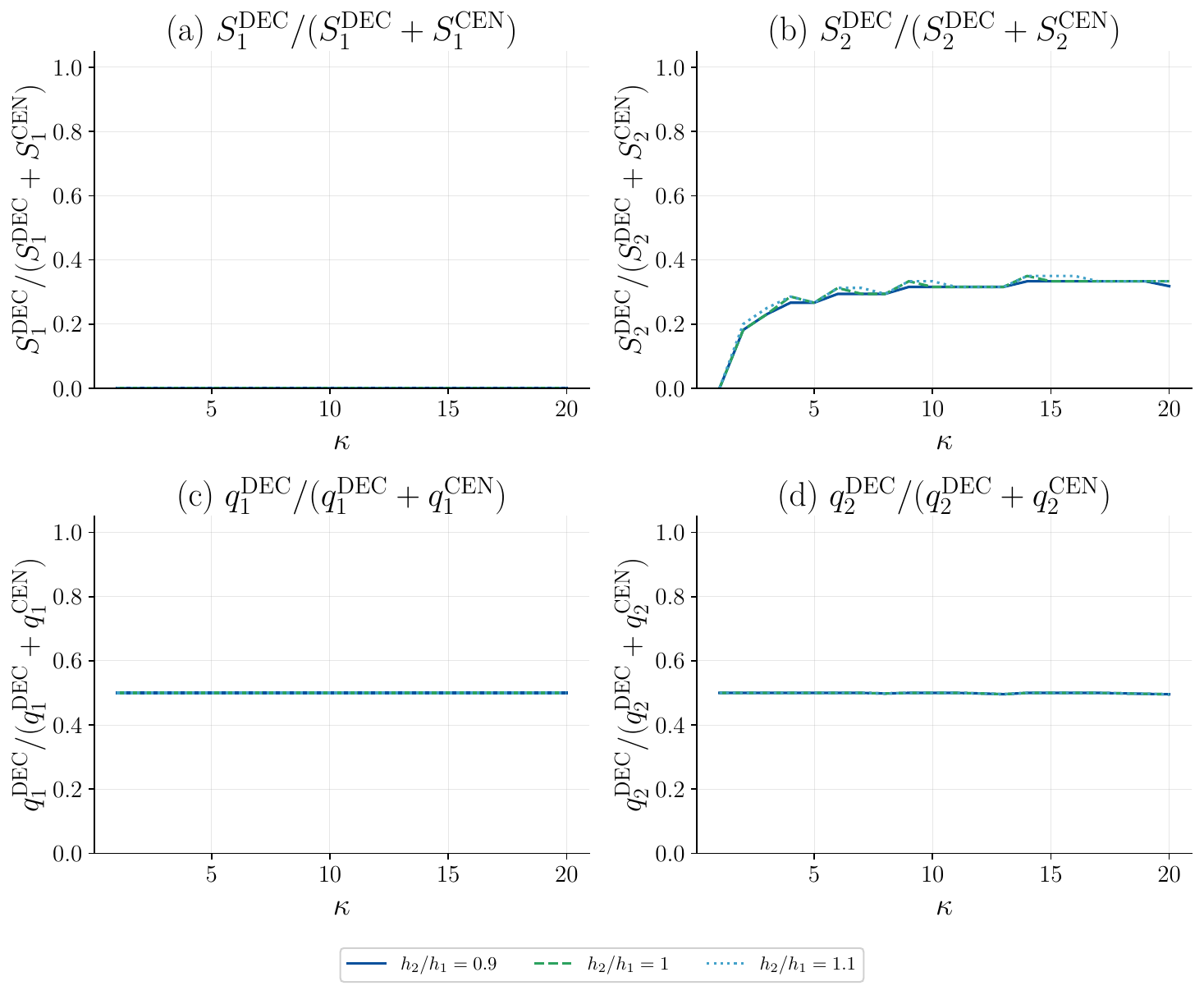}
\caption{Cross-section at $\rho = 0.8$: decentralized shares $X^{\textsc{dec}}/(X^{\textsc{dec}}+X^{\textsc{cen}})$ as functions of $\kappa$ for different $h_2/h_1$ ratios (reduced $h_1$). Top row: inventory shares. Bottom row: joining probability shares.}
\label{fig:cross-app2-dec-shares}
\end{figure}

\begin{figure}[!htbp]
\centering
\includegraphics[width=\textwidth]{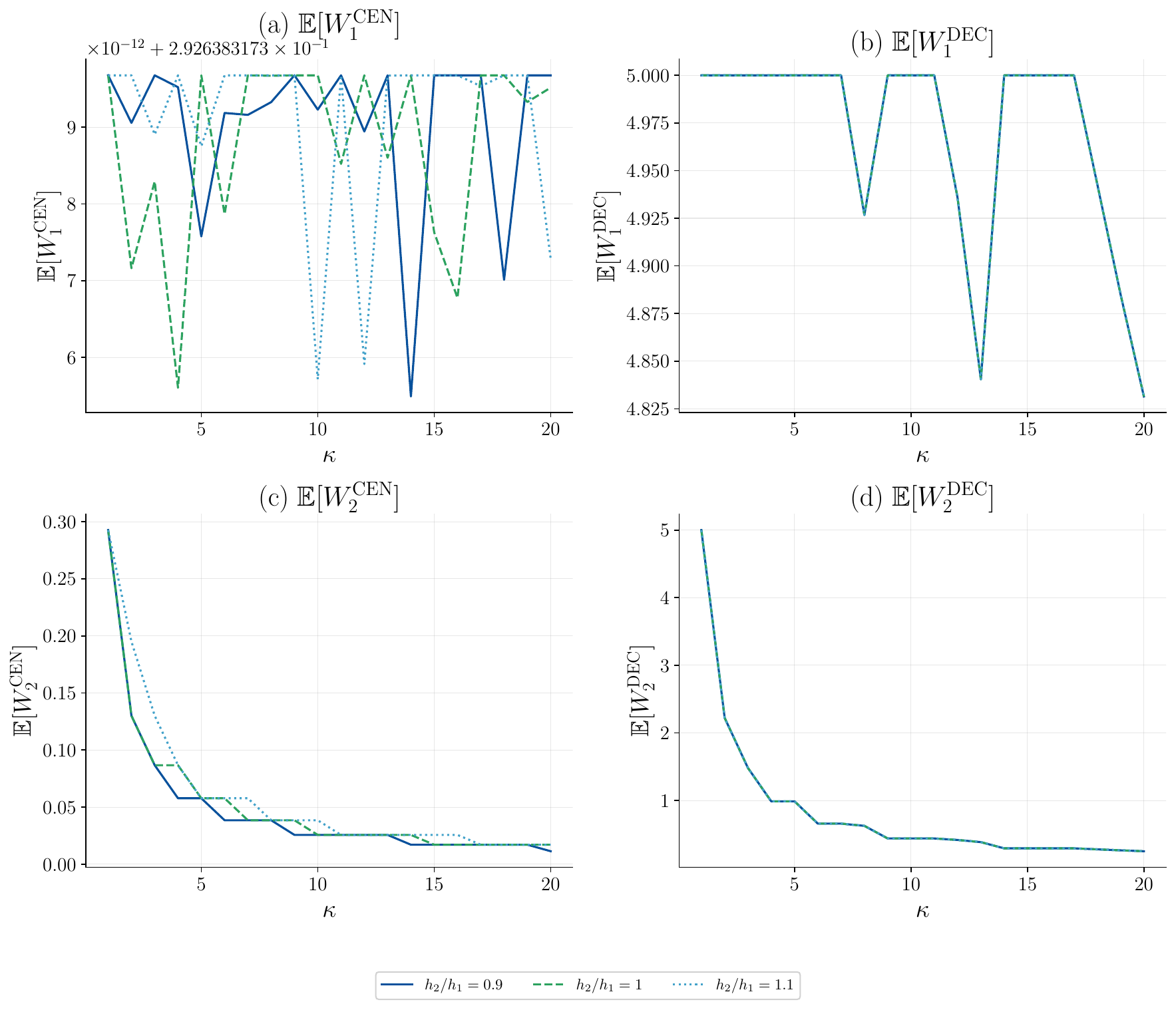}
\caption{Cross-section at $\rho = 0.8$: expected waiting times $E[W_i]$ as functions of $\kappa$ for different $h_2/h_1$ ratios (reduced $h_1$). Top row: type-1. Bottom row: type-2. Red segments indicate no joining.}
\label{fig:cross-app2-waiting}
\end{figure}

\begin{figure}[!htbp]
\centering
\includegraphics[width=\textwidth]{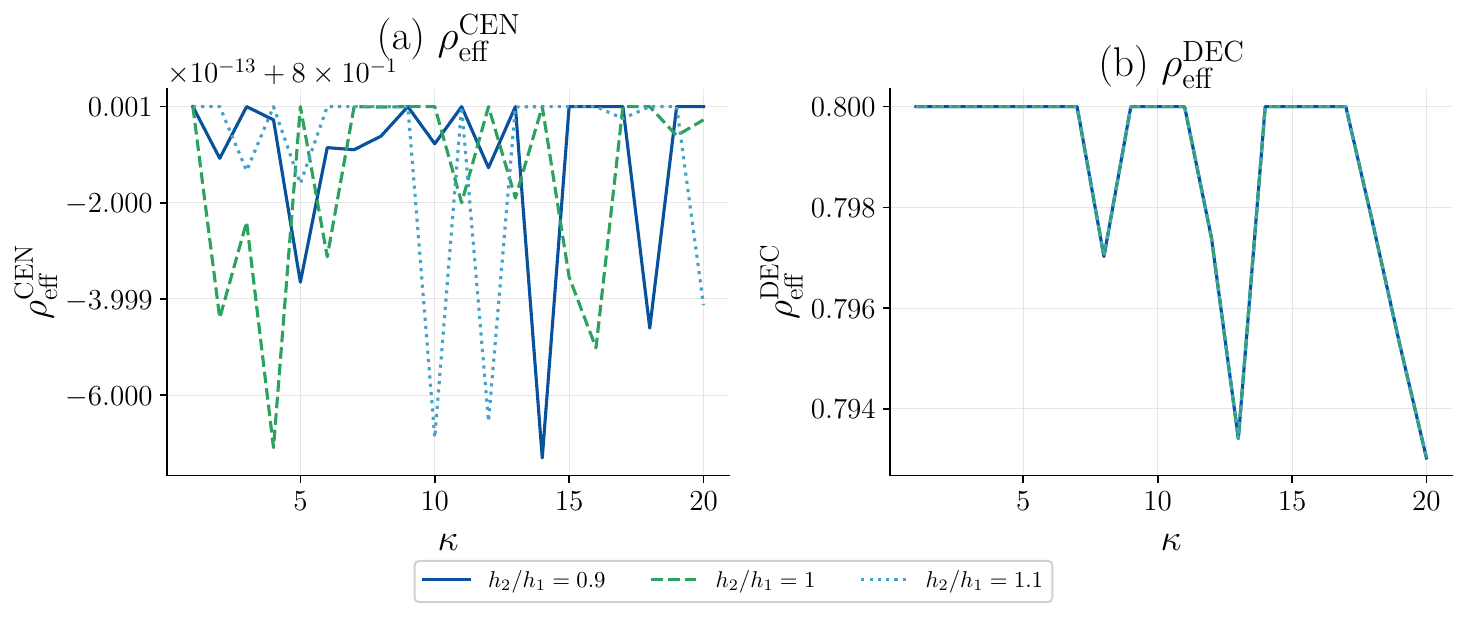}
\caption{Cross-section at $\rho = 0.8$: effective utilization $\rho_{\mathrm{eff}} = (q_1\Lambda_1 + q_2\Lambda_2)/\mu$ as function of $\kappa$ for different $h_2/h_1$ ratios (reduced $h_1$). Left: centralized. Right: decentralized.}
\label{fig:cross-app2-rho-eff}
\end{figure}
\FloatBarrier

\end{appendices}

\end{document}